\documentclass[11pt]{amsart}
\usepackage{amsmath,amsthm, amscd, amssymb, amsfonts, mathrsfs,amsaddr}
\usepackage[all]{xy}
\usepackage{color}
\usepackage{enumitem}
\usepackage{mathdots}

\usepackage{hhline} 
\usepackage[active]{srcltx}

\newcommand{\PSL}{\mathbf{PSL}}
\newcommand{\Gb}{\boldsymbol{G}}
\newcommand{\PSp}{\mathbf{PSp}}

\newcommand{\Pom}{\mathbf{P\Omega}}
\newcommand{\PSU}{\mathbf{PSU}}

\newcommand{\SU}{\mathbf{SU}}

\def\pibar{\overline{\pi}}
\newcommand{\X}{{\tt{X}}}
\newcommand{\Jf}{{\tt J}}

\newcommand{\bq}{\mathbf{q}}

\newcommand{\sa}{\dot{s}_\alpha}
\newcommand{\sbe}{\dot{s}_{\beta_2}}
\newcommand{\x}{\mathtt{x}}

\newcommand{\ta}{\mathtt{a}}
\newcommand{\tb}{\mathtt{b}}

\newcommand{\diag}{\operatorname{diag}}

\newcommand{\Irr}{\operatorname{Irr}}

\newcommand{\ord}{\operatorname{ord}}

\newcommand\w{\dot{w}}
\newcommand\sigmad{\dot{\sigma}}

\newcommand\toba{\mathfrak B }

\newcommand{\trid}{\triangleright}

\newcommand{\Ha}{{\mathbb H}}

\newcommand{\Y}{{\mathcal Y}}

\newcommand{\kk}{\Bbbk}

\newcommand{\ku}{\mathbb C}

\newcommand{\K}{{\mathcal K}}

\newcommand{\Z}{\mathbb Z}
\newcommand{\N}{{\mathbb N}}
\newcommand{\I}{{\mathbb I}}
\newcommand{\J}{{\mathbb J}}

\newcommand{\G}{{\mathbb G}}
\newcommand{\B}{{\mathbb{B}}}
\newcommand{\T}{{\mathbb{T}}}
\newcommand{\U}{\mathbb{U}}

\newcommand{\Pa}{\mathbb{P}}

\newcommand{\MM}{{\mathbb M}}

\newcommand{\F}{{\mathbb F}}

\newcommand{\GL}{\mathbf{GL}}
\newcommand{\GU}{\mathbf{GU}}
\newcommand{\PGL}{\mathbf{PGL}}
\newcommand{\SL}{\mathbf{SL}}
\newcommand{\Sp}{\mathbf{Sp}}
\newcommand{\SO}{\mathbf{SO}}

\newcommand{\Fr}{\operatorname{Fr}}

\newcommand{\Hc}{{\mathcal H}}

\newcommand{\Le}{{\mathbb L}}
\newcommand{\Vu}{{\mathbb V}}

\newcommand{\Oc}{{\mathcal O}}
\newcommand{\oc}{{\mathcal O}}

\newcommand{\ydgb}{{}^{\ku \Gb}_{\ku \Gb}\mathcal{YD}}
\newcommand{\ydg}{{}^{\ku G}_{\ku G}\mathcal{YD}}

\newcommand\ad{\operatorname{ad}}
\newcommand\Ad{\operatorname{Ad}}
\newcommand\Gsc{\G_{\operatorname{sc}}}

\numberwithin{equation}{section}

\theoremstyle{plain}

\newtheorem{lema}{Lemma}[section]
\newtheorem{theorem}[lema]{Theorem}

\newtheorem{prop}[lema]{Proposition}

\newtheorem{question-app}{Question}

\theoremstyle{definition}
\newtheorem{definition}[lema]{Definition}

\theoremstyle{remark}
\newtheorem{obs}[lema]{Remark}
\newtheorem{setting}[lema]{Setting}

\newcommand\wdot{\dot{w}}
\newcommand\id{\operatorname{id}}

\newcommand\A{\mathbb A}

\def\pf{\begin{proof}}
\def\epf{\end{proof}}

\theoremstyle{remark}

\newcounter{tabla}\stepcounter{tabla}

\begin{document}

\renewcommand{\baselinestretch}{1.2}

\thispagestyle{empty}
\title[Nichols algebras over mixed classes]
{Finite-dimensional pointed Hopf algebras\newline over finite simple groups of Lie type V. \newline 
Mixed classes in Chevalley and Steinberg groups}

\author[N. Andruskiewitsch, G. Carnovale, G. A. Garc\'ia]
{Nicol\'as Andruskiewitsch, Giovanna Carnovale$^*$ and\newline Gast\'on Andr\'es Garc\'ia}
 
\thanks{
	The work of N. A. was partially supported by CONICET, Secyt (UNC) and the Alexander von Humboldt Foundation
	through the Research Group Linkage Programme.
	The work of G. A. G. 
	was partially supported by CONICET, Secyt (UNLP) and ANPCyT-Foncyt 2014-0507.
	The work of G. C. was partially supported by Progetto BIRD179758/17 of 
	the University of Padova. The results were obtained during visits of N. A. 
	and G. A. G. to the University of Padova, and of G. C. to the University 
	of C\'ordoba,  partially supported by: the bilateral agreement between 
	these  Universities, the ICTP-INdAM Research in Pairs fellowship programme 
	and the Coimbra group  Scholarship Programme for Young Professors and 
	Researchers from Latin American Universities. The authors wish to thank the referee for his/her careful reading and suggestions.}

\address[A1]{FaMAF, 
	Universidad Nacional de C\'ordoba. CIEM -- CONICET\\ 
	Medina Allende s/n (5000) Ciudad Universitaria, C\'ordoba,\\
	Argentina\\
	andrus@famaf.unc.edu.ar}
\address[A2]{
	Dipartimento di Matematica Tullio Levi-Civita,\\
	Universit\`a degli Studi di Padova,\\
	via Trieste 63, 35121 Padova, Italia\\
carnoval@math.unipd.it, +39-049-8271354}
\address[A3]{
Departamento de Matem\'atica, Facultad de Ciencias Exactas,\\
	Universidad Nacional de La Plata. CONICET. C. C. 172, (1900)\\
	La Plata, Argentina\\
ggarcia@mate.unlp.edu.ar}

\begin{abstract}
We show that all classes that are neither semisimple nor unipotent in finite simple 
Chevalley or Steinberg groups different from $\PSL_n(q)$  
collapse (i.e. are never the support of a finite-dimensional Nichols algebra). 
As a consequence, we prove that  the only finite-dimensional pointed Hopf algebra whose group 
of group-like elements is $\PSp_{2n}(q)$, $\Pom^+_{4n}(q)$, $\Pom^-_{4n}(q)$, $^3D_4(q)$, 
$E_7(q)$, $E_8(q)$, $F_4(q)$, or $G_2(q)$ with $q$ even is the group algebra.
\end{abstract}
\maketitle
\setcounter{tocdepth}{1}

\tableofcontents

\section{Introduction}
This is the fifth of a series of papers where we study a set of group-theoretical questions on conjugacy classes of finite simple groups of Lie 
type because of their consequences on the classification of finite-dimensional Hopf algebras (over an algebraically closed field of characteristic 0, say $\ku$). 
Let us start with a brief discussion of the reduction from the latter problem to the former. See the Introductions of \cite{ACG-I,ACG-IV} for a more detailed exposition.

\begin{enumerate}[leftmargin=*]
\item  Let $G$ be a group. There is a braided tensor category $\ydg$ of so called Yetter-Drinfeld modules over 
(the group algebra of) $G$. Each $V \in \ydg$ gives rise to a graded Hopf algebra $\toba(V)$ in $\ydg$, 
called the Nichols algebra of $V$. See e.g. \cite{A-leyva} for details.
	
\item Let $\Hc$ be a pointed Hopf algebra whose group of group-like elements is isomorphic to $G$. 
There is a Nichols algebra $\toba(V)$ attached as a fundamental invariant to $\Hc$. 
For instance $\dim \Hc < \infty$ implies $\dim \toba(V) < \infty$; in fact $\toba(V)$ controls much of the structure of $\Hc$.
	
\item Assume that $G$ is finite. Then $\ydg$ is semisimple and its simple objects are parametrized by pairs 
$(\oc, \rho)$, where $\oc$ is a conjugacy class of $G$ and $\rho \in \Irr C_{G}(x)$ is an irreducible representation 
of the centralizer of a fixed element $x$ in $\Oc$; let $M(\Oc, \rho)$ be the simple object corresponding to $(\Oc, \rho)$. 
By the previous discussion, we need to determine the pairs $(\Oc, \rho)$ with $\dim \toba(M(\Oc, \rho)) < \infty$  
as a necessary initial step to classify finite-dimensional pointed Hopf algebras with group of group-likes $G$.
	
\item The second reduction goes as follows. It turns out that the Nichols algebra $\toba(M(\Oc, \rho))$ depends 
only (in  an appropriate sense) on the rack structure of $\oc$ (with rack operation given by conjugation: $x \trid y = xyx^{-1}$)
and a suitable 2-cocycle $\bq$ arising from $\rho$. In this way, it is more efficient to deal with Nichols algebras
$\toba(\Oc, \bq)$ of
pairs $(\oc, \bq)$ where $\oc$ is a rack and $\bq$ is a 2-cocycle, as the same such pair may arise from different groups.

\item Let us say that a rack $\oc$ \emph{collapses} if $\toba(\Oc, \bq)$ has infinite dimension for any
suitable cocycle $\bq$, cf. \cite[2.2]{AFGV-ampa}.
Remarkably there are rack-theoretical criteria that imply that a rack collapses, without computing 
neither cocycles nor Nichols algebras. These criteria can be spelled out in group-theoretical terms, 
once a realization of the rack in consideration as a conjugacy class is fixed. See \S \ref{subsec:racks}.
\end{enumerate}
	
In this series we investigate the applicability of these criteria for finite simple groups of Lie type,
excluding the Suzuki and Ree groups treated in \cite{CC-VI}. For 
alternating and sporadic groups, see \cite{AFGV-ampa,F,AFGV-sporadic,FV}. 
Let $\Gb$ be a finite simple Chevalley or Steinberg group defined over a finite field $\F_{q}$ where $q = p^m$,
with $m \in \N$ and $p$ prime.
An element of $\Gb$ is semisimple, respectively unipotent, if its order is coprime to $p$, respectively a power of $p$.
These notions are unambiguous except when $\Gb$ is one of the following:
\begin{align}\label{eq:exceptions}
\PSL_2(4) &\simeq \PSL_2(5),& \PSL_3(2) &\simeq \PSL_2(7),& \PSU_4(2) &\simeq \PSp_4(3).
\end{align}
A class is semisimple or unipotent if any element of it is so.
This paper deals with \emph{mixed} conjugacy classes  i.e., which are neither semisimple nor unipotent; 
evidently there are no mixed classes in $\PSL_2(q)$.
 For the groups in \eqref{eq:exceptions},
unipotent, semisimple or mixed will mean so in one of the realizations (with some ambiguity unless the realization
is clear from the context).
The mixed conjugacy classes of $\PSL_n(q)$ 
were treated in \cite{ACG-I}. In the previous four papers of the series \cite{ACG-I, ACG-II, ACG-III,ACG-IV}
we concluded the analysis of the unipotent conjugacy classes of $\Gb$
which is summarized in Table \ref{tab:unip-kthulhu}.
Elaborating on these results we obtain the first main theorem of this paper:

\begin{theorem}\label{thm:mixed-Chev-Steinberg}
A mixed conjugacy class in a simple Chevalley or Steinberg group collapses. 
\end{theorem} 

This Theorem and its proof propagate in many directions. For instance, if $\varpi: H \to \Gb$
is a projection from a finite group $H$ to our $\Gb$ and $x\in H$ is such that $\varpi(x)$ is neither semisimple nor unipotent, then $\Oc^H_{x}$ collapses. 
Here and below we use the notation $\Oc^H_{x}$ for the conjugacy class of $x$ in $H$ and we omit $H$ if clear from the context. 
For example let $H = \Gb_1 \times \Gb_2 \dots \times \Gb_t$, where all $\Gb_j$'s are simple Chevalley or Steinberg groups and 
$x = (x_1, \dots, x_t) \in H$. If $x_j$ is neither semisimple nor unipotent for at least one $j$, then 
$\Oc^H_{x}$ collapses.
Another application of Theorem \ref{thm:mixed-Chev-Steinberg} is to Nichols algebras of Yetter-Drinfeld modules over
 our $\Gb$.

\begin{theorem}\label{thm:collapse}
Let $\Gb$ be a simple Chevalley or Steinberg group and $V \in \ydgb$. Assume that $\Gb$ is not listed in \eqref{eq:exceptions}.
If $\dim\toba(V)<\infty$, then $V \simeq M(\Oc,\rho)$ is simple and $\Oc$ is semisimple.
\end{theorem}

Notice that we are not claiming the converse in Theorem \ref{thm:collapse}. Indeed Theorem \ref{thm:mixed-Chev-Steinberg}
is the culmination of our analysis of unipotent classes where we rely on the combinatorial description of such classes in the algebraic group behind $\Gb$. 
We started the consideration of the semisimple classes in $\Gb = \PSL_n(q)$ in \cite{ACG-III}; already in this case the situation is
much more involved and requires a different type of arguments.

\medbreak
As for the groups listed in \eqref{eq:exceptions}, we have the following. Let $V \in \ydgb$.

\begin{itemize}[leftmargin=*]\renewcommand{\labelitemi}{$\circ$}
\item $\Gb \simeq \PSL_2(4) \simeq \PSL_2(5) \simeq \A_5$: $\dim\toba(V) = \infty$ by \cite{FGV1,FGV2}.

\item $\Gb \simeq \PSL_3(2) \simeq \PSL_2(7)$: If $\dim\toba(V)<\infty$, then $V$ is simple, 
hence isomorphic to $M(\Oc_x,\rho)$
for some $x$ and $\rho$ by  \cite[Corollary 8.3]{HS}. We conclude that the order of $x$ is 4 and $\rho(x) = -1$ by \cite{FGV2}.

\item $\Gb \simeq \PSU_4(2) \simeq \PSp_4(3)$: $\dim\toba(V) = \infty$. Indeed, 
assume otherwise. Arguing as in the proof of Theorem \ref{thm:collapse}, necessarily $V \simeq M(\Oc_x,\rho)$ 
and the order of $x$ is odd and not divisible by 3; by inspection, it is 5.
Now there is only one class of elements of such order, hence this class is real. But this contradicts \cite{AZ}. 
\end{itemize}

\medbreak
For some Chevalley or Steinberg groups the arguments for Theorem \ref{thm:collapse} can be pushed even further. 
A folklore conjecture  claims that there is no finite-dimensional pointed Hopf algebra whose group 
of group-like elements is simple non-abelian, except the group algebra. 
The conjecture is known to hold for the simple alternating groups \cite{AFGV-ampa}, 
the sporadic groups (including the Tits group) except $Fi_{22}$, $B$ and $M$ \cite{AFGV-sporadic} and for some families of 
$\PSL_{n}(q)$'s \cite{FGV1,ACG-III}. As we just saw, it also holds for $\PSU_4(2) \simeq \PSp_4(3)$; 
otherwise the Conjecture is open.
We add to the list of confirmations of the Conjecture the families in \eqref{eq:gps-1inW}, except from $\PSL_2(q)$ 
for $q$ even, that had already been treated.

\begin{theorem}\label{thm:collapse-gp}
Let $q$ be even and let $\Gb$ be a simple group in one of the following families:
\begin{align}\label{eq:gps-1inW}
\begin{aligned}
&\PSL_2(q),& &\PSp_{2n}(q),& &\Pom^+_{4n}(q),& &\Pom^-_{4n}(q), && \\
&^{3}D_4(q),& &E_7(q), & & E_8(q),& & F_4(q),& & G_2(q).
\end{aligned}
\end{align}
Let  $\Hc$ be a finite-dimensional pointed Hopf algebra whose group of group-like elements is isomorphic to $\Gb$.
Then $\Hc \simeq \ku \Gb$.
\end{theorem}
 
\medbreak
Here is the structure of the article. In Section \ref{sec:preliminaries} we set 
the notation and collect basic preliminary results used in the paper. 
In order to prove Theorem \ref{thm:mixed-Chev-Steinberg} we adopt the following strategy. 
Let $x =x_sx_u$ be the Chevalley-Jordan decomposition of $x\in \Gb$ 
(which amounts to the $p$-decomposition in the theory of finite groups)
where we assume that neither $x_s$ nor $x_u$ are trivial. 
Recall that there is a surjective group homomorphism $\pi: \G^F \to \Gb$, 
where $\G$ is a simply connected simple algebraic group and $F$ is a Steinberg endomorphism--see Section \ref{sec:preliminaries} for details.
We pick $\x \in \G^F$ such that $\pi(\x) = x$. If $\x = \x_u\x_s$ is its Chevalley-Jordan decomposition, then 
$x_u = \pi(\x_u)$.
There is a natural inclusion of racks \eqref{eq:reduction-unipotent} 
of $\Oc_{\x_u}^{\K}$ in $\Oc_{x}^{\Gb}$, where $\K$ is
the centraliser of $\x_s$ in $\G^F$; we restrict our attention to $\Oc_{\x_u}^{\K}$. 
The structure of $\K$ is known and 
$\Oc_{\x_u}^{\K}$ contains a product of unipotent conjugacy classes of smaller 
groups of Lie type.  
Thus the unipotent classes in the factors of the centraliser must occur in 
Table \ref{tab:unip-kthulhu}. 
In this way  we are able to provide strong restrictions to potential kthulhu classes, 
see \eqref{cond:Jv=1} and \eqref{cond:Jv=2}.
Then we focus on the resulting reduced list of classes and we bring the action of the 
Weyl group into the picture, providing further techniques to deal 
with most of the potential kthulhu classes. 
These methods are still quite general and allow to 
cover a considerable list of groups, see Proposition \ref{prop:center}, 
and leave out a few exceptions in the remaining ones. 
The substance of the preceding sketch is contained in Section \ref{sec:mixed}.
The classes for which all above techniques 
fail can be described with sufficient precision, 
allowing us to carry out the final analysis. 
However, this part of the proof is laborious and it needs  several ad-hoc 
considerations and a separation of the treatment for Chevalley (Section \ref{sec:chevalley}) and Steinberg groups
(Section \ref{sec:steinberg}), where a case-by-case analysis is performed.
Finally we prove Theorems \ref{thm:collapse} and \ref{thm:collapse-gp} in Section \ref{sec:nichols}.

\section{Notation and Preliminaries}\label{sec:preliminaries}
In this Section we establish notation, recall some known facts about racks and simple algebraic groups over finite fields, and collect 
some results on unipotent conjugacy classes 
that are used throughout the paper.

For $k\leq l$, $k,l\in\N$ we set $\I_{k,l}:=\{k,\,k+1,\,\ldots,\,l\}$ and $\I_{l}=\I_{1,l}$.  

Let $p$ be a prime, $q=p^m$ for $m\in\N$, $\kk=\overline{\F_q}$. 

Let $\G$ be a simple algebraic group over $\kk$; we fix a maximal torus $\T$ and a 
Borel subgroup $\B\supset\T$. We denote by $\U$ the unipotent radical of $\B$. 
Let $\Phi$ be the root system of $\G$, let $\Phi^+$ be the set of positive roots 
associated with $(\T,\B)$ and let $\Delta=\{\alpha_i,\,i\in\I_\ell\}$ 
be the corresponding base, with 
$\omega_i$, for $i\in \I_\ell$ the set of fundamental weights.
The root subgroup corresponding to $\alpha$ will be denoted by $\U_\alpha$. It is the image of 
a monomorphism of groups $x_\alpha: \kk\to \G$.
Since we will consider conjugacy classes as racks (see Subsection \ref{subsec:racks}), 
the conjugation action of $g\in \G$ on $h\in \G$  will be indicated by $g\trid h$. 
Adopting the notation from \cite[8.1.4]{springer}, we have the commutation rule:
\begin{align}\label{eq:comm-rule}
t\trid x_{\alpha}(a) &= t x_{\alpha}(a)t^{-1}=x_{\alpha}(\alpha(t)a),& t\in \T, \ \alpha&\in\Phi,\,\ a\in \kk.
\end{align}
We shall also make use of  Chevalley's commutator
formula \cite[Lemma 15, p. 22 and Corollary, p. 24]{yale}. Namely, for $\alpha, \beta \in\Phi^+$ such that $\alpha+\beta\in\Phi^+$ and any order on the set
$\Gamma$ of positive roots of the form $i\alpha+j\beta$ for $i,j\in\N$,  there exist integers $c_{ij}^{\alpha\beta}$ such that
\begin{align}\label{eq:chev}x_{\alpha}(\xi)x_\beta(\eta)x_\alpha(\xi)^{-1}x_\beta(\eta)^{-1}
&= \prod_{i\alpha+j\beta\in \Gamma}  x_{i\alpha + j\beta}(c_{ij}^{\alpha\beta}\xi^i\eta^j), &
\forall \xi, \eta&\in \kk.
\end{align}
For further unexplained notation and commutation rules we refer to \cite[Section 3]{ACG-II}.
By $\alpha_0$ we will denote the highest root in $\Phi$, and we set $\widetilde {\Delta}:=\Delta\cup\{-\alpha_0\}$. 
As usual, $W=N_{\G}(\T)/\T$ denotes the Weyl group, 
acting naturally on $\T$ and $s_\alpha$ is the reflection with respect to the root $\alpha$. 
Let $\sigma\in W$ be represented by $\sigmad\in N_\G(\T)$ and let $\MM$ be a subgroup of $\G$ normalized by  $\T$. Then the subgroup $\sigmad \MM\sigmad^{-1}$  is independent of the choice of $\sigmad$ and we shall simply write $\sigma\MM$.  

If we insist that $\G$ is simply-connected, we shall write $\G_{sc}$. 
For an algebraic group $\MM$, we shall denote by $\MM^\circ$ the connected component containing the identity. 
If $\phi$ is an automorphism of a group $H$, then $H^\phi$ will denote the set of elements of $H$ that are fixed by $\phi$. 

\subsection{The groups $\Gb$}\label{subsec:F}
Let $F$ be a Steinberg endomorphism of $\G$, and assume $\T$ and $\B$ are $F$-stable.
Let $\Gb:= \G_{sc}^F/Z(\G_{sc}^F)$ and let $\pi\colon \G^{F}_{sc}\to \Gb$ be the natural projection.
In most cases $\Gb$ is simple and, conversely, every finite simple group of Lie type is obtained this way
(except the Tits group that was treated in the paper \cite{AFGV-sporadic} with the sporadic groups). 
Here we deal with Chevalley or Steinberg groups; thus, $F$ is the composition of a standard graph automorphism $\vartheta$ of $\G$ with the Frobenius map $\Fr_q$. Then $\vartheta$ determines an 
automorphism $\theta$ of the Dynkin diagram of $\G$.
If $\theta=\id$, then $\Gb$ is a Chevalley group, whereas if $\theta\neq\id$, then $\Gb$ is a Steinberg group. 
In the latter case, $\T$ and $\B$ are chosen to be both $F$-stable and $\Fr_q$-stable and in both cases $\Fr_q(t)=t^q$ for any $t\in\T$.
Since $\T$ is $F$-stable, $F$ acts on $W$ and we denote by $W^{F}$ the subgroup fixed by $F$.
Thus $W^{F} = W$ for Chevalley groups.
For each $w\in W^{F}$, there exists a representative $\w \in N_{\G^{F}}(\T)$, see \cite[Prop. 23.2]{MT}.
Also, for any representative $\w \in N_{\G}(\T)$  of $w\in W$ and any $\alpha\in\Phi$ one has that $\w \trid \U_{\alpha} = \U_{w(\alpha)}$.

\subsection{$F$-stable tori}\label{subsec:tori}We recall that a maximal torus $g\T g^{-1}$ in $\G$ is $F$-stable if and only if 
$\dot{w}=g^{-1}F(g)\in N_{\G}(\T)$ and this relation induces a map $\phi$ from the set of $F$-stable 
maximal tori to $W$.  
In addition, two $F$-stable maximal tori $\T_1$ and $\T_2$  are conjugate by an element in $\G^F$  if and only if their images $\phi(\T_1)$ and $\phi(\T_2)$
lie in the same $F$-twisted conjugacy class  in $W$ \cite[Prop. 25.1]{MT}. If this is the case, their fixed points subgroups  $\T_1^F$ and $\T_2^F$  are also conjugate in $\G^F$.
We will indicate by $\T_w$ a maximal torus $g\T g^{-1}$ such that $g^{-1}F(g)=\dot{w}\in N_{\G^{\Fr_q}}(\T)$.  
Observe that every $w$ in 
$W^{\Fr_q}=W$ has a representative 
$\w$ in $N_{\G^{\Fr_q}}(\T)$ by \cite[\S 23.1]{MT}  and existence of a $g\in\G$ such that 
$g^{-1}F(g)=\dot{w}$ is guaranteed by the Lang-Steinberg theorem. Hence, such a torus always exists. 

\subsection{The elements $w_{\theta}$}\label{subsec:steinberg_w_0}
Let $\vartheta$ be as above and let $\theta$ be the associated automorphism of the Dynkin diagram; it
extends to an automorphism of the root system of $\G$ that we denote again by $\theta$. 
Assume that $\theta=\id$ or $\theta=-w_0$, 
where $w_0$ is the longest element in $W$. 
The latter occurs for Steinberg groups with root system of type: $A_\ell$ for $\ell\geq2$; $D_{\ell}$ 
for $\ell\geq 5$ and odd; and $E_6$. 
Let $w \in W$ and take $\T_w$, $\w\in N_{\G^{\Fr_q}}(\T)$ and $g$ as in Subsection \ref{subsec:tori}. 
We set $w_\theta:=ww_0^{|\theta|-1}$. 
Since $C_W(w_\theta)=C_W(w\theta)$,
\cite[Prop. 25.3]{MT} guarantees that there exists a representative 
$\w_{\theta}$ of $w_\theta$ in $N_\G(\T)$ such that $g\w_\theta g^{-1}\in N_{\G^F}(\T_w)$. 
If $\theta=-w_0$, then we set $\w_0:=\w^{-1}\w_{\theta}$.

\subsection{Semisimple classes}\label{subsec:semisimple}
Every semisimple element of $\G^F$ lies in an $F$-stable torus 
$\widetilde \T$ \cite[Prop. 26.6]{MT}. In particular, if $\phi(\widetilde \T)=w$ and 
$s\in\G^F\cap \widetilde \T$, then $\T_w^F$ is $\G^F$-conjugate to $\widetilde \T^F$, 
so $\Oc_s^{\G^F}\cap \T_w\neq\emptyset$. 

\smallbreak
Let $s\in \G^F\cap \T_w = \T_w^F$, with $g,\,\w$ as in Subsection \ref{subsec:tori}. 
We set $F_w:=\Ad(\w)\circ F$. Then $s=gtg^{-1}$ for some $t\in \T^{F_w}$, 
so 
\begin{align}\label{eq:Ft}
F(t)=w^{-1}(t) = \Ad(\dot{w}^{-1})(t).
\end{align} 

Sometimes it is convenient to consider other $F$-stable tori intersecting $\Oc_s^{\G^F}$.
Assume that $\G=\G_{sc}$ and pick $\sigma\in W$ such that $\sigma^{-1}(t)=F(t)$, i.~e.
$\sigma$ and $w$ lie in the same left coset for the stabiliser of $t$ in $W$.
Then for every representative $\sigmad\in N_{\G}(\T)$
and every $g_\sigma$ such that $g_\sigma^{-1}F(g_\sigma)=\sigmad$ we have 
\begin{align*}
s':=g_\sigma t g_\sigma^{-1}\in \G^F\cap\Oc_t^{\G}\overset{\star}{=}\Oc_s^{\G^F}
\end{align*} 
so  $\Oc_s^{\G^F}$ has a representative in $\T_\sigma^F$ for any such $\sigma$.	
Here $\star$ holds because $C_\G(s)$ is connected, since $\G=\Gsc$ is simply-connected, \cite[Theorem 2.11]{hu-cc}. 

\subsection{Commuting elements in a semisimple class}\label{subsec:semisimple-comm}
Let $s, w, g$ and $t$ be as in Subsection \ref{subsec:semisimple}.
We look for powers of $s$ lying in $\Oc_s^{\G^F}$.
If $\theta=\id$ or $\theta=-w_0$, we have:
\begin{align*}
 (g\w_\theta^{-1} g^{-1})\trid s&=(g\w_{\theta}^{-1})\trid t= g\trid(w_0^{|\theta|-1}w^{-1}t)\\
 &=g\trid (w_0^{|\theta|-1}\theta t^q)
 =\begin{cases}\begin{matrix}
g\trid t^q=s^q,& \text{ if }\theta=\id,\\
g\trid t^{-q}=s^{-q},& \text{ if }\theta=-w_0,\end{matrix}
\end{cases}
\end{align*} 
so either $\Oc_s^{\G^F}=\Oc_{s^q}^{\G^F}$  (when $\theta=\id$) or 
$\Oc_s^{\G^F}=\Oc_{s^{-q}}^{\G^F}$ ($\theta=-w_0$). 

Assume that $\theta=\id$. Even when $s^q=s$ we get some information.
In this case, $t\in \G^F$, so $t\in \Oc_s^{\G} \cap \G^F$. 
If, in addition, $\G=\Gsc$ then $C_\G(s)$ is 
connected, \cite[Theorem 2.11]{hu-cc}, and $\Oc_s^{\G^F}=\Oc^{\G^F}_t$. Hence $\Oc_s^{\G^F}\cap\T^F\neq\emptyset$,
that is, the class intersects $\T$.

\subsection{Real classes}\label{subsec:w0}
Assume $w_0=-\id$ and let $s$, $\T_w$, $\w$ and $g$ as in Subsection \ref{subsec:tori}. Then $w_0\in C_W(w\theta)$ 
so there exists 
a representative $\w_0$ of $w_0$ in $N_\G(\T)$ such that $g\w_0g^{-1}\in N_{\G^F}(\T_w)$. Then
\begin{align*}
 (g\w_0 g^{-1})\trid s&=(g\w_0)\trid t= g\trid t^{-1}=s^{-1}\in \Oc_s^{\G^F}.
\end{align*}

\subsection{Non-central elements in the torus}\label{subsec:alpha} We need the following fact on algebraic groups.
Assume $\G=\G_{sc}$. 
We recall that, for $t\in \T-Z(\G)$, there always exists $\alpha\in \Delta$ such that $s_\alpha(t)\neq t$. Indeed, suppose that $s_\alpha(t)=t$ for every $\alpha\in\Delta$. 
Then $w(t)=t$ for every $w\in W$. 
By \cite[\S 2.2]{hu-cc}, $W$ is then the Weyl group of the reductive group 
$C_{\G}(t)=\langle \T, \U_\alpha~|~\alpha(t)=1\rangle$. 
Now the Dynkin diagram of $C_{\G}(t)$ can be read off from the Coxeter graph of $W$, except when $\G$ is
of type $B_\ell$ or $C_\ell$. 
This forces $C_{\G}(t)=\G$ (because the two groups have the same Dynkin diagram, 
or by comparison of the dimensions), a contradiction as $t$ is not central.

Such an $\alpha$ also satisfies $\alpha(t) \neq 1$.
Indeed, $t$ commutes with $\U_{\alpha}$ if and only if $\alpha(t)=1$ if and only if $t$ commutes
with $\U_{-\alpha}$. In addition,
$s_\alpha (t)=\dot{s}_\alpha \trid t$ for any representative $\dot{s}_\alpha$ of $s_\alpha$ in $N_{\G}(T)$.
Since there is a representative of $s_\alpha$ lying in $\U_\alpha\U_{-\alpha}\U_\alpha$ \cite[8.1.4]{springer}, 
we get a contradiction and whence the claim follows.

\subsection{Racks}\label{subsec:racks} 
We recall the notion of rack and some
definitions that will be needed later. In the present paper all racks are (unions of) conjugacy classes. 
See \cite{ACG-III} for details and more information.

A \emph{rack} is a set $X \neq \emptyset$  with a self-distributive
operation $\trid: X \times X \to X$ such that ${\mathtt x}\trid \underline{\quad}$ is
bijective for every ${\mathtt x} \in X$. The standard example is the conjugacy class $\Oc^M_{\mathtt z}$
of an element ${\mathtt z}$ in a group $M$ with the operation ${\mathtt x}\trid {\mathtt y} = {\mathtt x}{\mathtt y}{\mathtt x}^{-1}$, 
for ${\mathtt x},{\mathtt  y} \in \Oc^M_{\mathtt z}$. 
A subrack of a rack $X$ is a non-empty subset $Y$ such that ${\mathtt x}\trid {\mathtt y} \in Y$ for all ${\mathtt x},\,{\mathtt y} \in Y$.  
A rack $X$ is \emph{abelian} if ${\mathtt x}\trid {\mathtt y} = {\mathtt y}$, for all ${\mathtt x}, {\mathtt y} \in X$.
In \cite{AFGV-ampa,ACG-I,ACG-III} the notions of racks of type D, F and C were introduced.  
Since we will  deal with subracks of conjugacy classes, 
we recall the translation of these notions in this case. 

\begin{definition} A conjugacy class $\Oc$ in a finite group $M$ is of type
\begin{enumerate}[leftmargin=*]
\item[\textbf{C}]  if there exist $H\leq M$ and $r,\,s\in \Oc\cap H$ 
such that
\begin{enumerate}
\item $rs\neq sr$,
\item $\Oc_r^H\neq \Oc_s^H$,
\item $H=\langle \Oc_r^H,\, \Oc_s^H\rangle$,
\item either $\min(| \Oc_r^H|,\, |\Oc_s^H|)>2$ or $\max(| \Oc_r^H|,\, |\Oc_s^H|)>4$;
\end{enumerate}
\item[\textbf{D}] if there exists $r,\,s\in\Oc$ such that
\begin{enumerate}
\item $(rs)^2\neq(sr)^2$,
\item $\Oc_r^{\langle r,\,s\rangle}\neq \Oc_s^{\langle r,\,s\rangle}$;
\end{enumerate}
\item[\textbf{F}]
 if there exists $r_i\in\Oc$, for $i\in\I_4$ such that
\begin{enumerate}
\item $r_ir_j\neq r_jr_i$, for $i\neq j$,
\item $\Oc_{r_i}^{\langle r_i,\,i\in\I_4\rangle }\neq \Oc_{r_j}^{\langle r_i,\,i\in\I_4\rangle }$, for $i\neq j$.
\end{enumerate}
\end{enumerate}
A rack is \textit{kthulhu} if it is neither of type C, D nor F.
\end{definition}

\begin{obs} \label{obs:types-quotient}
The properties of type C, D, F are well-behaved with respect to projections and inclusions, \cite[Section 3.2]{AFGV-ampa}, \cite[Remark 2.9 (a)]{ACG-I},\cite[Lemma 2.10]{ACG-III}. 
In other words, if a rack $X$ has either a quotient or a subrack which is not kthulhu, then $X$ is not kthulhu. 
\end{obs}

The main motivation for seeking conjugacy classes that are not kthulhu stems from the classification of finite-dimensional Nichols algebras.  

\begin{theorem} \cite[Theorem  3.6]{AFGV-ampa},  
\cite[Theorem 2.8]{ACG-I},  \cite[Theorem 2.9]{ACG-III}.
A conjugacy class $\Oc$ of type D,  F or C collapses. \qed
\end{theorem}

\begin{obs}\label{obs:subgroup}Let $M$ be a finite group whose order is divisible by $p$. 
Let ${\mathtt y},{\mathtt z} \in M$ and 
write ${\mathtt y}={\mathtt y}_{s}{\mathtt y}_{u}$ and ${\mathtt z}={\mathtt z}_{s}{\mathtt z}_{u}$ as
products of their $p$-regular and $p$-parts. Then  $\langle {\mathtt y},\,{\mathtt z}\rangle=\langle {\mathtt y}_s, {\mathtt y}_u,{\mathtt z}_s,{\mathtt z}_u\rangle$ and if 
${\mathtt y}_s {\mathtt z}_u\neq {\mathtt z}_u{\mathtt y}_s$ or ${\mathtt y}_s {\mathtt z}_s\neq {\mathtt z}_s{\mathtt y}_s$ then 
${\mathtt y} {\mathtt z}\neq {\mathtt z} {\mathtt y}$. In particular, if $M$ is a finite group of Lie type, then 
$\langle {\mathtt y},\, {\mathtt z}\rangle$ is generated by the semisimple and unipotent parts of ${\mathtt y}$ and ${\mathtt z}$ and non-commutation of 
${\mathtt y}$ and ${\mathtt z}$ can be checked on their semisimple or unipotent parts.
\end{obs}

The next lemma is instrumental to detect classes of type C in finite groups. 

\begin{lema}\label{lem:tecnico0} 
Let $M$ be a finite group whose order is divisible by $p$.
Let ${\mathtt y}, {\mathtt z} \in M$ and write ${\mathtt y}={\mathtt y}_{s}{\mathtt y}_{u}$ and ${\mathtt z}= {\mathtt z}_{s} {\mathtt z}_{u}$ as
products of their $p$-regular and $p$-parts. 
If
\begin{equation}\label{eq:ineq-tecnico0}
{\mathtt y}_{s} {\mathtt z}_{u}\neq {\mathtt z}_{u}{\mathtt y}_{s} \quad \text{ and }  \quad {\mathtt z}_{s}{\mathtt y}_{u}\neq {\mathtt y}_{u} {\mathtt z}_{s},
\end{equation}
then ${\mathtt y} {\mathtt z} \neq {\mathtt z} {\mathtt y}$, $|\Oc_{{\mathtt y}}^{\langle {\mathtt y}, {\mathtt z}\rangle}|\geq 3$ and 
$|\Oc_{ {\mathtt z}}^{\langle {\mathtt y}, {\mathtt z}\rangle}|\geq 3$.
If in addition ${\mathtt y}$ and ${\mathtt z}$ are conjugate in $M$
and $\Oc_{{\mathtt y}}^{\langle {\mathtt y}, {\mathtt z}\rangle}\cap \Oc_{ {\mathtt z}}^{\langle {\mathtt y}, {\mathtt z}\rangle} = \emptyset$, then $\Oc_{{\mathtt y}}^{M}$
is of type C.
\end{lema}
\pf By Remark \ref{obs:subgroup} any of the inequalities in \eqref{eq:ineq-tecnico0} implies that 
${\mathtt y} {\mathtt z}\neq {\mathtt z} {\mathtt y}$. Moreover, $\Oc_{{\mathtt y}}^{\langle {\mathtt y}, {\mathtt z}\rangle}$ contains 
the orbits $\langle {\mathtt z}_{s} \rangle \trid {\mathtt y}$ 
and $\langle {\mathtt z}_{u} \rangle \trid {\mathtt y}$. 
Since by \eqref{eq:ineq-tecnico0} we know that ${\mathtt z}_{u}$ and ${\mathtt z}_{s}$ are nontrivial, 
one of these two elements has odd order bigger than $1$. This implies that 
$|\Oc_{{\mathtt y}}^{\langle {\mathtt y}, {\mathtt z}\rangle}|\geq 3$. The inequality $|\Oc_{ {\mathtt z}}^{\langle {\mathtt y}, {\mathtt z}\rangle}|\geq 3$ 
follows by symmetry. The last assertion follows by taking the subgroup $H=\langle {\mathtt y}, {\mathtt z}\rangle$.
\epf

\subsection{Known results on unipotent classes}\label{subsec:unipotent}

In this paper we analyze conjugacy classes that are neither semisimple nor 
unipotent in finite simple groups of Lie type (mixed classes), in order to determine when such a class is non-kthulhu. 
By  Remark \ref{obs:types-quotient}, it is enough to find subracks or 
quotients which are not kthulhu.
With this aim, we will first look for non-kthulhu subracks that are isomorphic to unipotent classes in 
simple groups, relying on the results in \cite{ACG-I,ACG-II, ACG-III,ACG-IV}   
on unipotent conjugacy classes. These are summarized in the following theorem.

\begin{theorem}\label{thm:slspsu}
	Let $\Gb=\G_{sc}^F/Z(\G_{sc}^F)$ be a Chevalley or a Steinberg group and let $\Oc \neq \{e\}$ be a 
	unipotent conjugacy class in $\Gb$, not listed in Table \ref{tab:unip-kthulhu}. Then $\Oc$ is not kthulhu.
	\qed
\end{theorem}
	\begin{table}[ht]
		\caption{Kthulhu unipotent classes}\label{tab:unip-kthulhu}
		\begin{center}
			\begin{tabular}{|c|c|c|}
				\hline $\Gb$   & class & $q$ \\
				\hline
				\hline
				$\PSL_2(q)$   & $(2)$ & even, or $9$, or odd not a square   \\
				\hline
				$\PSL_3(2)$  &$(3)$  & $2$   \\
				\hline
				$\PSp_{2n}(q)$, $n\geq 2$  & $W(1)^{n-1}\oplus V(2)$  & even \\
				\hline
				$\PSp_{2n}(q)$, $n\geq2$  & $(2, 1^{2n-2})$  & $9$, or odd not a square\\
				\hline
				$\PSp_{4}(q)$  &  $W(2)$  & even  \\
				\hline
				$\PSU_{n}(q)$, $n\geq 3$ &  $(2, 1^{n-2})$ & even \\
				\hline
			\end{tabular}
		\end{center}
	\end{table}

\begin{obs}\label{obs:properties-kthulhu}
We list some properties of the unipotent classes from Table \ref{tab:unip-kthulhu}, see 
\cite{ACG-I, ACG-II, ACG-IV}.
\begin{enumerate}[leftmargin=*]
 \item\label{item:una-o-due} For every group in Table \ref{tab:unip-kthulhu} there are at most  two unipotent conjugacy classes  
 labeled as in the second column of the table and they are isomorphic as racks.  There is only one class labeled by $(3)$ in $\PSL_3(2)$. 
 Also, there is only one class labeled by $(2)$ in $\PSL_2(q)$ for $q$ even.
 \item\label{item:rappresentanti} Assume $\Oc$ is a class different from the one labeled by $(3)$ in $\PSL_3(2)$. Then
 $\Oc$ is represented by an element in $\U_\beta$  where $\beta$ is either the highest root in 
 $\Phi^+$ or the highest short root in $\Phi^+$, so $F(\U_\beta)=\U_\beta$.  
If there are two such classes, then one of these is represented by $x_\beta(1)$. 
\item\label{item:invarianza} All classes in  Table \ref{tab:unip-kthulhu} are $\Fr_p$-stable. Indeed,  $\Fr_p$ being an automorphism of 
 $\Gb$, stabilizes the set of kthulhu classes and preserves the label and the order of the elements. 
 If there are two classes with the same label, the one represented by $x_\beta(1)$ is clearly $\Fr_p$-stable, 
 therefore $\Fr_p$ must stabilize the second one, too.
 \item\label{item:iso} Let $q$ be even. When $n=2$, the non-standard graph automorphism interchanging long and 
 short roots in $\PSp_{4}(q)$  maps the class labeled by $W(2)$ to the one labeled by $W(1)\oplus V(2)$. For any $n\geq2$, 
 there is only one class labeled by $W(1)^{n-1}\oplus V(2)$ in $\PSp_{2n}(q)$. 
 \item\label{item:una-sola} There is only one class labeled by $(2, 1^{n-2})$ in $\PSU_{n}(q)$ for $q$ even. 
 Indeed, let $u\in \SU_n(q)$ have label $(2, 1^{n-2})$. Then $u$ is conjugate in $\GU_n(q)$ to 
 $x_{\beta}(1)$. Since every element in $\GU_n(q)$ is the product of a diagonal matrix $d$ fixed by $F$ and a matrix in $\SU_n(q)$, 
  $u$ is conjugate in $\SU_n(q)$ to an element in  $\U_{\beta}^F$. A direct computation shows that $\U_\beta^F=\U_{\beta}^{\Fr_q}$ 
 and that all nontrivial elements in $\U_{\beta}^{\Fr_q}$ are  conjugate to $x_{\beta}(1)$ by a matrix in $\SU_n(q)$ of the form 
 $\diag(\eta,1\,\,\ldots,\,1,\eta^{-1})$ for $\eta\in\F_q^\times$. 
 \item Some of the groups occurring in the table are not simple. We need to consider them to implement an induction procedure.
 \end{enumerate}
 \end{obs}

We also recall, for the purpose of application of Theorem \ref{thm:slspsu}, that unipotent conjugacy classes in a simple group of Lie 
type $\Gb=\G_{sc}^F/Z(\G_{sc}^F)$ are isomorphic as racks to the corresponding unipotent conjugacy class
in $\G^F$ for any $\G$ isogenous to $\G_{sc}$.

\subsection{Results on products of racks}\label{subsec:product}
Another approach for detecting non-kthulhu conjugacy classes
is to find subracks isomorphic to 
a product of conjugacy classes in a smaller group and apply results from \cite{ACG-I,ACG-II}.
With this in mind, we will make use of the following Lemma.
\begin{lema} \cite[Lemma 2.10]{ACG-I}\label{lem:productD}
Let $\Oc$ be a conjugacy class in a finite group $M$ containing a subrack 
of the form $X_1\times X_2$ such that
\begin{align}
\label{eq:x1x2}\mbox{there are }  & x_1,x_2\in X_1 \mbox{ such that } (x_1x_2)^2\neq (x_2x_1)^2,\\
\label{eq:y1y2}\mbox{there are }  & y_1\neq y_2\in X_2\mbox{ such that } y_1y_2= y_2y_1.
\end{align}
Then $\Oc$ is of type D.\qed
\end{lema}

In order to apply Lemma \ref{lem:productD}, we determine which unipotent conjugacy classes 
satisfy \eqref{eq:x1x2} or \eqref{eq:y1y2}.
\begin{lema}\label{lem:product}
If $X_1$ is a conjugacy class occurring in Table \ref{tab:unip-kthulhu}, then it satisfies \eqref{eq:x1x2}. If $X_2$ is a  
class occurring in Table \ref{tab:unip-kthulhu} but not in 
Table \ref{tab:unip-not-y1y2}, then it satisfies
\eqref{eq:y1y2}. In particular, if $\oc$ is a conjugacy class of a finite group $M$ containing a 
	 subrack isomorphic to a product of unipotent conjugacy 
classes in Chevalley or Steinberg groups, not both listed in  
	 Table \ref{tab:unip-not-y1y2}, then $\Oc$ is not kthulhu.
\end{lema}
	\begin{table}[ht]
		\caption{Kthulhu unipotent classes not satisfying \eqref{eq:y1y2} }\label{tab:unip-not-y1y2}
		\begin{center}
			\begin{tabular}{|c|c|c|}
				\hline $\Gb$   & class & $q$ \\
				\hline
				\hline
				$\PSL_2(q)$   & $(2)$ & $2,\,3$   \\
				\hline
				$\PSU_{3}(q)$  &  $(2, 1)$ & $2$ \\
				\hline
			\end{tabular}
		\end{center}
	\end{table}
\pf Let $\overline{\Oc}$ be the class in $\PSU_n(q)$ occurring in Table \ref{tab:unip-kthulhu}. 
For all conjugacy classes occurring in Table \ref{tab:unip-kthulhu} different from $\overline{\Oc}$, 
conditions \eqref{eq:x1x2} and   \eqref{eq:y1y2} are proved in \cite[Lemma 3.5]{ACG-IV}.
Assume $X_1=\overline{\Oc}$. We set $x_1=\id_n+e_{1,n}$, $\sigma:=e_{1,n}+e_{n,1}+\sum_{i\neq1,n}e_{i,i}$, $x_2:=\sigma\trid x_1$ 
and \eqref{eq:x1x2} holds. We now prove \eqref{eq:y1y2} for $X_2=\overline{\Oc}$ and $(n,q)\neq(3,2)$. Let $\zeta$ be a generator of $\F_{q^2}^\times$. 
If $n=3$ and $q>2$ we take
\begin{align*}
d={\rm diag}(\zeta,\zeta^{q-1},\zeta^{-q}),&&y_1=\id_3+e_{1,3},&&y_2=d\trid y_1.
\end{align*}
If $n>3$, we take 
    $\sigma:=e_{1,2}+e_{2,1}+e_{n-1,n}+e_{n,n-1}+\sum_{i\in\I_{3,n-2}}e_{i,i}$, $y_1=\id_n+e_{1,n}$ and $y_2:=\sigma\trid y_1$. This settles the first statement.

Let now $X=X_1\times X_2\subset\Oc$ with $X_i$ for $i=1,2$ isomorphic to a product of unipotent classes in Chevalley or Steinberg groups. If $X_1$ or $X_2$ is not in 
Table \ref{tab:unip-kthulhu}, then $X$ is not kthulhu. If  $X_1$ nd $X_2$ are in Table \ref{tab:unip-kthulhu} and $X_2$ is not in 
Table   \ref{tab:unip-not-y1y2}, then the statement follows from the preceding and Lemma \ref{lem:productD}.
    \epf
    
    \begin{lema}\label{lem:product_different}Let $\Oc\neq \Oc'\subset \G^F$ be unipotent conjugacy classes in 
    Table \ref{tab:unip-kthulhu} corresponding to the same label. If $q\neq 3$, then there exist  $x_1\in\Oc$, $x_2\in \Oc'$  
    such that $ (x_1x_2)^2\neq (x_2x_1)^2$. If $q=3$,  then there exist  $x_1\in\Oc$, $x_2\in \Oc'$  such that $ x_1x_2\neq x_2x_1$. 
    \end{lema}
    
    \pf  By Remark \ref{obs:properties-kthulhu} \eqref{item:una-o-due} and  \eqref{item:una-sola} the classes $\Oc$ 
    and $\Oc'$ occur either in $\PSL_2(q)$ or in $\PSp_{2n}(q)$. By  Remark \ref{obs:properties-kthulhu} \eqref{item:rappresentanti}  
    we may assume $x_{\beta}(1)\in\Oc$ for $\beta$ the highest root or the highest short root in $\Phi^+$ and  $y_1\in\Oc'\cap\U_\beta^F$. 
    Let $\dot{s}_{\beta}$ be a representative of $s_\beta$ in $N_{\G^F}(\T)$ and let $\dot{s}_{\beta}\trid y_1=x_{-\beta}(\xi) \in\Oc'\cap\U_{-\beta}^F$. 
    A computation in $\langle \U_\beta,\U_{-\beta}\rangle$ shows that 
    $x_{\beta}(1)x_{-\beta}(\xi)\neq x_{-\beta}(\xi)x_{\beta}(1)$ always and that
    $(x_{\beta}(1)x_{-\beta}(\xi))^2\neq(x_{-\beta}(\xi)x_{\beta}(1))^2$ if and only if $\xi(2+\xi)\neq0$. 
    A direct verification in $\SL_2(q)$ and $\Sp_{2n}(q)$ shows that if $q\neq3$, then $\dot{s}_{\beta}$ can always be chosen in such a way that 
    $\xi$ satisfies the latter condition. We take $x_1=x_{\beta}(1)$ and $x_2=\dot{s}_{\beta}\trid y_1$.\epf
    
    \begin{obs}
    The inequalities in Lemmata \ref{lem:product} and \ref{lem:product_different} can be written as inequalities between unipotent elements.  
    Since the restriction of an isogeny  to unipotent elements is an isomorphism, such inequalities hold independently of the  isogeny type of $\G$.\end{obs}

\section{Mixed classes}\label{sec:mixed}

In this Section we introduce some reduction techniques in order to deal with mixed classes. 
From now on  $\G=\G_{sc}$ is simply-connected and $\Gb \neq \PSL_{n}(q)$, since 
these last groups have been treated in \cite{ACG-I}. 

Let $x\in \Gb$ and let $\x\in \G^F$ such that $\pi(\x) = x$, with Chevalley-Jordan decomposition $\x = \x_s\x_u$.
Then $x$ has Chevalley-Jordan decomposition $x=x_sx_u$ with $x_s = \pi(\x_s)$, $x_u = \pi(\x_u)$. 
We assume that $\Oc_x^{\Gb}$ is mixed, that is $x_s,\,x_u\neq1$.
By construction $\x_u$ belongs to $\K := C_{\G}(\x_s)\cap\G^F$, thus $x_u\in K := \pi(\K)$. 

We recall the morphisms of racks from \cite[Lemma 1.2]{ACG-I}
\begin{equation}\label{eq:reduction-unipotent}
\oc^{\K}_{\x_u}\simeq  \oc^{K}_{x_u} \hookrightarrow \Oc_x^{\Gb},
\end{equation}
where the isomorphism follows from injectivity of the restriction of the isogeny 
$\G\to \G/Z(\G)$ to unipotent elements. 
This motivates the quest for classes $\Oc_\x^{\G^F}$ such that $\Oc_{\x_u}^{\K}$ is 
not kthulhu and a better understanding of the group $\K$. For the latter we will use \cite{carter-classical}.

\subsection{Classes in centralisers of semisimple elements} \label{sec:unipotent_in_K}

By the discussion in Subsection \ref{subsec:semisimple}, we may assume that  
$\x_s=gtg^{-1}$ for some $t\in\T$, and some $g\in \G$ such that 
$g^{-1}F(g)=\dot{w}\in N_{\G^{\Fr_q}}(\T)$. Let $\Ha:=C_\G(t)$. By \cite[2.2]{hu-cc}, 
$\Ha$ is connected and reductive, so $\Ha=Z(\Ha)^\circ [\Ha,\Ha]$, it is given by
\begin{align*}
\Ha=\left\langle \T,\,\U_\alpha~|~\alpha\in\Phi,\,\alpha(t)=1\right\rangle.
\end{align*}
and $\Fr_q \Ha=C_{\G}(t^q)=\Ha$. Since $\Phi_t:=\{\alpha\in\Phi~|~\alpha(t)=1\}$ is 
$W$-conjugate to a root system with base a subset of 
$\widetilde \Delta$, 
 \cite[Proposition 30]{McS}, 
up to replacing $t$ by $\sigmad\trid t$,  $g$ by $g\sigmad^{-1}$  and $\w$ by $\sigmad\w F(\sigmad^{-1})$ for $\sigmad$ a 
suitable element in $N_{\G}(\T)$, we may always assume that $\Phi_t$ has a base $\Pi$ contained in $\widetilde \Delta$. 
The subset $\Pi$ is unique up to  $W$-action. 

If $\Ad(\w)^{-1}(t)=F(t)$, then $\Ad(h\w)^{-1}(t)=F(t)$ for every $h\in \Ha$, 
so $\w$ could be replaced by $\sigmad \w\in\G^{\Fr_q}$ for any $\sigmad\in N_{\Ha}(\T)$ and $w$ could be replaced by 
$\sigma w$, for $\sigma\in N_\Ha(\T)/\T=W_\Pi$, the group generated by the reflections with respect to roots in $\Pi$.

Let $F_w:=\Ad(\wdot)\circ F$. Then $\Ha$ is $F_w$-stable and
\begin{align*}\K&=(g \Ha g^{-1})^F=g\Ha^{F_w}g^{-1},\end{align*}
 \cite[Remark 2.5(c)]{ACG-II}.
By uniqueness of the Chevalley-Jordan decomposition, $\x_u$ lies in 
$(g[\Ha,\Ha]g^{-1})^F= g[\Ha,\Ha]^{F_w}g^{-1}$ . 

For the rest of the paper, we assume that 
$\x = \x_s\x_{u}= gtvg^{-1} $ with $t\in \T$, $v\in [\Ha,\Ha]^{F_w}$ and 
$\x_s= gtg^{-1}$, $\x_u=gvg^{-1}$.

This leads us to the following statement.
\begin{lema}\label{lem:reduction_to_HH}With notation as above, 
if $\Oc_v^{[\Ha,\Ha]^{F_w}}$ is not kthulhu, then $\Oc_x^{\Gb}$ is again so.
\end{lema}
\pf 
This follows from \eqref{eq:reduction-unipotent} and the inclusion $g\trid \Oc_v^{[\Ha,\Ha]^{F_w}}\subset \Oc_{\x_u}^{\K}$.
\epf

We explore now several conditions ensuring that the hypothesis of 
Lemma \ref{lem:reduction_to_HH} is satisfied. 
In order to do so, we describe the structure of the $F_w$-stable and $\Fr_q$-stable, 
hence $\Ad(\wdot)\vartheta$-stable,  
semisimple  group  $[\Ha,\Ha]$.  There exist uniquely determined 
$\Fr_q$-stable simple algebraic subgroups $\G_j\leq [\Ha,\Ha]$ for $j\in\I_r$  
satisfying:
\begin{align}\label{eq:deco1}
[\Ha,\Ha]=\G_1\cdots \G_r,&& [\G_i,\G_j]=1\mbox{ if }i\neq j,&&\G_i\cap \prod_{j\neq i}\G_j\subseteq Z([\Ha,\Ha]). 
\end{align}
We will denote by $\Phi_{t}^i$ the root system of $\G_i$ with base $\Delta_i\subseteq\widetilde \Delta$.
The automorphism $\Ad(\wdot)\vartheta$ permutes the factors $\G_i$ and the systems $\Phi_{t}^i$, 
inducing a permutation of the indices $i\in\I_r$ 
which we denote by $\omega$. 
By suitably rearranging the indices, we may assume that $\omega$ is a product of $l$ disjoint cycles of the form 
$c_j=(i_j+1,\ldots,i_j+a_j)$ with $\sum_ja_j=r$ and $i_j=\sum_{b=1}^{j-1}a_b$, 
so $\omega$ is completely determined by the
$a_j$, $j\in\I_l$. Let $C_j:=\{i_j+1,\ldots,i_j+a_j\}$ and $\Ha_j:=\prod_{l\in C_j}\G_l$. 
Then each $\Ha_j$ is 
semisimple and $F_w$-stable, since it is $\Fr_q$-stable and 
$\Ad(\wdot)\vartheta$-stable. As a consequence of  \eqref{eq:deco1} we have:
\begin{align}\label{eq:deco2}
[\Ha,\Ha]=\Ha_1\cdots \Ha_l,&& [\Ha_i,\Ha_j]=1\mbox{ if }i\neq j,&&
\Ha_i\cap \prod_{j\neq i}\Ha_j\subseteq Z([\Ha,\Ha]). 
\end{align}
The element $v$ decomposes accordingly as $v=\prod_{j=1}^l v_j$ with $v_j\in \Ha_j$ unipotent. The equation
\begin{align*}\prod_{j=1}^l F_w(v_j)=F_w(v)=v=\prod_{j=1}^l v_j\end{align*} implies that
$v_j F_w(v_j)^{-1}=z_j \in Z([\Ha,\Ha])$.  Hence, $v_j =z_j F_w(v_j)$ with $v_j$, $F_w(v_j)$ 
unipotent, forcing $z_j=1$ for every $j\in \I_l$. 
Thus, $v\in \prod_{j=1}^l\Ha_j^{F_w}$. 

We analyse now the structure of $\Ha_j^{F_w}$, for $j\in\I_l$. 
For simplicity, we assume for the moment that  $l=1$, so $v=v_1$. We set $a:=a_1$. 
Let $u_k$ be the (unipotent) component of $v$ in $\G_k$. Since $v$ is $F_w$-invariant, we have
\begin{align*}
\prod_{j=1}^{a}u_j =v=F_w(v)=\prod_{j=1}^{a}F_w(u_j),& &F_w(u_j)\in\G_{j+1}, j\in\I_{a-1},&&F_w(u_a)\in\G_{1}.
\end{align*}
By \eqref{eq:deco1} and the unipotency of the $u_j$ we have
\begin{align}\label{eq:simple1}
u_j=F_w^{j-1}(u_1),\,j\in\I_{1,a}&&F_w^a(u_1)=u_1.
\end{align}
In other words, $v$ lies in the subgroup $G_1$ of $\Ha_1$ consisting of elements whose components satisfy \eqref{eq:simple1}. 
Projection onto the first component induces an isomorphism $G_1\simeq\G_1^{F_w^a}$. 
Since, by construction, $\Fr_q(\dot{w})=\dot{w}$, we have 
\begin{align*}
F_w^a&=(\Ad(\w)\vartheta)^a\circ \Fr_q^a
\end{align*}
where $(\Ad(\w)\vartheta)^a$ is an automorphism of the simple algebraic group $\G_1$ and $F_w^a$ 
is a Steinberg endomorphism of $\G_1$. 
The corresponding finite simple group will be  Chevalley or Steinberg 
according to whether $(\Ad(\w)\vartheta)^a$ is an inner automorphism of $\G_1$ or not \cite[10.9]{St}. 

If $l>1$, a similar analysis shows that $v\in\prod_{j=1}^lG_j$ where $G_j\simeq \G_{i_j+1}^{F_w^{a_j}}$. 
Thus, the rack $\Oc_v^{[\Ha,\Ha]^{F_w}}$ 
contains the subrack $\Oc_v^{\prod_jG_j}$ which is isomorphic to $\prod_j\Oc_{v_j}^{G_j}$. 
Each component is a unipotent conjugacy class among those studied in \cite{ACG-I,ACG-II,ACG-III, ACG-IV}.
By abuse of notation we will identify $G_j$ and $\G_{i_j+1}^{F_w^{a_j}}$ and $v_j$ 
with its component in $\G_{i_j+1}^{F_w^{a_j}}$. 

\medbreak
Let $\J_v:=\{j\in \I_l~|~v_j\neq1\}$. 

\begin{lema}\label{lem:reduction}Assume that one of the following conditions holds for  $\Oc_{x}^{\Gb}$:
\begin{enumerate}
\item\label{item:uno} For some $j\in \J_v$, the class  $\Oc_{v_j}^{G_j}$  is not in Table \ref{tab:unip-kthulhu},
\item\label{item:due} $|\J_v|\geq2$ and for some $j\in \J_v$,  $\Oc_{v_j}^{G_j}$ 
is not in Table \ref{tab:unip-not-y1y2}.
\end{enumerate}
Then $\Oc_{x}^{\Gb}$ is not kthulhu. 
\end{lema}
\pf \eqref{item:uno} The rack $\Oc_{\x_u}^{\K}$ contains $\Oc_{v_j}^{G_j}$ which is not kthulhu. 
\eqref{item:due} The rack $\Oc_{\x_u}^{\K}$ contains a subrack  satisfying 
the hypothesis of Lemma \ref{lem:product}.  
\epf

\subsection{Subracks obtained by the action of the Weyl group}\label{subsec:weyl}

From now on we assume that we are in  the following situation: either

\begin{eqnarray}
|\J_v|=1\text{ and }\Oc_{v_1}^{G_1}\text{ occurs in Table \ref{tab:unip-kthulhu}, or}\qquad\qquad\qquad\qquad\quad\ \ \label{cond:Jv=1}\\   
\qquad \qquad |\J_v|\geq2\text{ and for every }j\in \J_v\text{ the rack }\Oc_{v_{j}}^{G_j}\text{ occurs in Table \ref{tab:unip-not-y1y2}}\label{cond:Jv=2}.  
 \end{eqnarray}

The latter case occurs only when $q^{a_j}\in\{2,3\}$ for every $j\in \J_v$, 
whence $\omega(j)=j$ for every $j\in \J_v$. 
Observe that if $(w\theta)^{a_j}(\Phi_t^j\cap\Phi^+)=\Phi_t^j\cap\Phi^+$ 
for some $j\in\J_v$, then we may always assume $v_j\in\U$. In particular, if  $q\neq2$ and $w\theta$ acts trivially on $\Phi_t^j$, 
then by Remark \ref{obs:properties-kthulhu} we may take $v_j\in\U^F_{\beta}$ for $\beta$ the highest root or the highest short root in $\Phi_t^j$.

\begin{obs}\label{obs:rackiso-conjg}
In order to deal with mixed  classes satisfying   \eqref{cond:Jv=1} or \eqref{cond:Jv=2}, we need to look at subracks of $\Oc_\x^{\G^F}$ 
that are different from those in  Lemma \ref{lem:reduction}.
Since for every $F_w$-stable subgroup $\MM$ of $\G$, we have
$(g\MM g^{-1})^F=g \MM^{F_w}g^{-1}$, we have the rack isomorphism 
\begin{align}\label{eq:isomorphism}g\trid\underline{\,\,\,}: \Oc_{tv}^{\G^{F_w}} \overset{\sim}{\to} \Oc_\x^{\G^F}, \end{align} 
and we will mainly perform calculations in $\G^{F_w}$.
We will often make use of the fact that any 
$\sigma$ in the centraliser $C_W(w\theta)$ has a representative in $N_{\G}(\T)\cap \G^{F_w}$, 
\cite[Proposition 25.3]{MT}. In addition $Z(\G^F)=Z(\G^{F_w})$ 
so the isomorphism of racks \eqref{eq:isomorphism} is compatible with the isogeny $\G_{sc}\to \G_{\ad}$.
\end{obs}

\begin{lema}\label{lem:no-centro}
Assume that there exists $\sigma\in W$ such that: 
\begin{enumerate}[leftmargin=*]
\item\label{item:comm} $\sigma \in C_W(w\theta)$;
\item\label{item:preserve} $\sigma(\Ha_j)=\Ha_j$ for some $j\in \J_v$;
\item\label{item:nocentro} $\sigma (t)\not\in Z(\G^F)t$.
\end{enumerate}
Then  $\Oc_x^{\Gb}$ is not kthulhu.
\end{lema}
\pf 
Observe that the subgroup $\sigma\Ha_j$ is well-defined because  $\Ha_j$ is normalized by $\T$. Let $\sigmad\in  N_{\G}(\T)\cap \G^{F_w}$ be a 
 representative of $\sigma$. Let
\begin{align*}X:=t \Oc^{\Ha^{F_w}}_v,&&\text{ and }&&Y:=\sigmad\trid X=\sigma(t)\Oc_{\sigmad\trid v}^{\Ha^{F_w}}.
\end{align*} 
Condition \eqref{item:nocentro} ensures that  $X\cap Y=\emptyset$ and that
  the restriction of $\pi$ to $X\coprod Y$ is injective. 
 By Condition \eqref{item:comm},  $\sigma(t)$ commutes with $\sigma(\Ha_j)=\Ha_j$, whence with $G_j$. In addition,  
 $\Ha_j$ commutes with all $\Ha_i$ and $\sigma(\Ha_i)$ for $i\neq j$, hence it commutes with the  subgroup they generate.
Thus, in order to verify the non-commutation of elements in $X$ and 
 $Y$ it is enough to look at the components in  $\Oc_{v_j}^{G_j}$ and $\Oc_{\sigmad\trid v_j}^{G_j}$. 
 For this reason we assume for simplicity that $\J_v=\{1\}$, so $v=v_j=v_1$. Notice that  $\sigmad\trid v\in G_1$ has the same label as $v$.  
\begin{enumerate}[leftmargin=*]
\item[(a)] Assume $q$ even. Since $v$ and $\sigmad \trid v$ have the same label, $\sigmad\trid v\in\Oc_{v}^{G_1}$ by Remark \ref{obs:properties-kthulhu}.  
By  Lemma \ref{lem:product} applied to $\Oc_{v}^{G_1}$,
 there exist $x_{1}, x_{2} \in \Oc_{v}^{G_1}$ such that 
$(x_{1}x_{2})^{2}\neq (x_{2}x_{1})^{2}$. Then, for $r:=tx_1\in X$ and  $s:=\sigma(t)x_2\in Y$  
we have $\Oc_{r}^{\langle r,s\rangle}\subseteq X$, $\Oc_{s}^{\langle r,s\rangle}\subseteq Y$ and $(rs)^2\neq (sr)^2$. 
This shows that $\Oc_\x^{\G^F}$ is of type $D$. 

\item[(b)] Assume $q$ odd. Either by Lemma \ref{lem:product} or by Lemma \ref{lem:product_different}, there exists $v'\in\Oc_{\sigmad\trid v}^{G_1}$ such that $vv'\neq v'v$. 
Let 
\begin{align*}
&r=tv\in X,&&s=\sigma(t)v'\in Y,&&H:=\langle r,\,s\rangle=\langle t,\sigma(t),v,\,v'\rangle.
\end{align*} 
By construction, $t,\sigma(t)\in Z(H)$,  $rs\neq sr$, $\Oc_r^H\subset X$,  $\Oc_s^H\subset Y$, so 
$\Oc_r^H\cap \Oc_s^H=\emptyset$. Also,
$H\leq\langle \Oc_r^H,\,\Oc_s^H\rangle\leq H$. 
Finally, $v$ and $v'$ are $p$-elements,  so $v, v'\trid v$ and $(v')^2\trid v$ are all distinct, hence
\begin{align*}
\left|\Oc_r^H\right|=\left|t\Oc_v^H\right|=\left|\Oc_v^H\right|=\left|\Oc_v^{\langle v,v'\rangle}\right|\geq 3
\end{align*}  and similarly for $\left|\Oc_s^H\right|$.
%
\end{enumerate}
Since $\Oc_x^{\Gb}$ contains the subrack $\pi(X \amalg Y) \simeq X\amalg Y $, from the discussion 
above it follows that $\Oc_x^{\Gb}$ is of type $D$ or $C$.
\epf

\begin{obs}\label{obs:conditions}
Conditions  \eqref{item:comm} and \eqref{item:preserve} from 
Lemma \ref{lem:no-centro} are verified in the following situations:
\begin{enumerate}
\item[(a)]\label{item:q+1} $\theta=\id$ or $\theta=-w_0$ and $\sigma=w_\theta^{-1}$. Here   $\sigma(t)=t^{\pm q}$.
\item[(b)]\label{item:square} $w_0=-1$ and $\sigma=w_0$. Here $\sigma(t)=t^{-1}$.
\end{enumerate}
In the above situations, $C_\G(\sigma(t))=\Ha$, so \eqref{item:preserve} holds for any $j\in\J_v$.
\end{obs}

\begin{lema}\label{lem:specifico}Assume that one of the following conditions holds
\begin{enumerate}
\item $\Gb$ is Chevalley and $t^q\not\in Z(\G^F)t$.
\item $\Gb$ is Steinberg with $\theta=-w_0$ and $t^{-q}\not\in Z(\G^F)t$.
\item $w_0=-1$ and $t^2\not\in Z(\G^F)$.
\end{enumerate}
Then $\Oc$ is not kthulhu.
\end{lema}
\pf 
It is a direct consequence of Lemma \ref{lem:no-centro} and Remark \ref{obs:conditions} because $w_\theta^{-1}(t)=t^q$ 
when $\theta=\id$ and $w_\theta^{-1}(t)=t^{-q}$  when $\theta=-w_0$.
\epf

The following remark will be useful to locate the cases in which the hypotheses of Lemma \ref{lem:no-centro}  or, 
more specifically,  \ref{lem:specifico} do not hold.

\begin{obs}\label{obs:zeta}Let $\tau\in W$ and $s=\prod_{j\in\I_\ell}\alpha_j^\vee(\xi_j)\in\T$. We set 
\begin{align*}
\I_{(\tau)}&=\{i\in\I_\ell~|~s_i\textrm{ occurs in a reduced decomposition of } \tau\}, \textrm{ and }\\
\I^{(s)}&=\{i\in\I_\ell~|~\xi_j\neq1\}.
\end{align*}
We recall that $\I_{(\tau)}$ does not depend on the chosen reduced decomposition and that the coefficients $\xi_j$ are uniquely determined. 
\begin{enumerate}
\item For $\tau$ and $s$ as above, $\tau(s)\in s\prod_{i\in\I_{(\tau)}}\alpha_i^\vee(\kk)$.
\item If for some $\sigma\in W$, $z\in Z(\G^F)$ and $t\in \T$ we have $\sigma(t)=z t$ then $\I^{(z)}\subseteq \I_{(\sigma)}$. 
\item In the special case in which $\theta=\id$ or $\theta=-w_0$, $z\in Z(\G_{sc}^F)- 1$ and 
$t\in\T^{F_w}$ satisfy $w_\theta^{-1}(t)=zt$, we necessarily have $\I^{(z)}\subseteq \I_{(w_\theta)}$.
\end{enumerate}
\end{obs}

\begin{prop}\label{prop:center}
Assume $q$ is even and $\Phi$ is of type 
$A_1$, $B_n$, $C_n$, $D_{2m}$,  $E_7$, $E_8$,  $F_4$ or $G_2$. Then $\Oc_x^{\Gb}$ is not kthulhu.
\end{prop}
\pf
In this case $w_0=-1$ and $\ord t$ is odd. By Lemma \ref{lem:specifico} the statement could fail only if $t^2\in Z(\G^F)$, 
but this would force $t\in Z(\G^F)$ contradicting our assumption on $\x_s$.
\epf

We end this Subsection
with a lemma 
that will be useful to discuss some  of the cases not covered by Lemmata \ref{lem:reduction}, \ref{lem:specifico} or Proposition \ref{prop:center},
both in Chevalley and Steinberg groups.

\begin{lema}\label{lem:A3}
Let $\{\beta_i,\,i\in\I_3\}\subset\Phi$ be the base of a  root subsystem of type $A_3$ such that
\begin{enumerate}
\item $\{\beta_1,\,\beta_3\}\subset \Phi_t$ and $\{\beta_1,\,\beta_3\}$ is $w\theta$-stable;
\item $\U_{\beta_2}$ is $F_w$-stable and $\beta_2\not\in \Phi_t$;
\item $v\in((\U_{\beta_1}-1)\U_{\beta_3})^{F_w}$. 
\end{enumerate}
Then $\Oc_{x}^{\Gb}$ is of type C. 
\end{lema}

\pf Let $\widetilde {\U}=\langle \U_{\beta_i},\,i\in\I_3\rangle$ 
and $\widetilde {\B}=\T\widetilde {\U}$. Then $\widetilde {\U}$ and $\widetilde {\B}$ are $F_w$-stable 
as well as $\Fr_q$-stable. 
Let $r:=tv=tx_{\beta_1}(\xi)x_{\beta_3}(\xi')\in t\widetilde {\U}$ 
for $\xi\in{\kk}^\times$, $\xi'\in{\kk}$. 
Since $w\theta(\beta_2)=\beta_2$ we have $w\theta s_{\beta_2}=s_{\beta_2}w\theta$ 
so we may find a representative $\sbe$ of $s_{\beta_{2}}$ in $\G^{F_w}\cap N_{\G}(\T)$.
We set: $t'=\sbe\trid t$, $v':=\sbe\trid v=x_{\beta_1+\beta_2}(\zeta)
x_{\beta_2+\beta_3}(\zeta')$ for some $\zeta\in {\kk}^\times$ and $\zeta'\in\kk$, and 
$s:=\sbe\trid r=t'v'\in t'\widetilde {\U}$. Since $\langle r,s\rangle\subset \widetilde {\B}$, then 
$\Oc_r^{\langle r,s\rangle}\subset t\widetilde {\U}$ and  
$\Oc_s^{\langle r,s\rangle}\subset  t'\widetilde {\U}$. 
In addition, since $s_{\beta_{2}}(\beta_{1}) \notin \Phi_{t}$ by $(2)$, 
it follows that $C_{\G}(t')\neq \Ha$ so $t'\not\in Z(\G)t$. 
Therefore,  $\pi(\Oc_r^{\langle r,s\rangle})\cap\pi(\Oc_s^{\langle r,s\rangle})=\emptyset$. 
Since $(\beta_{1}+\beta_{2})(\pi(t)) = (\beta_{1}+\beta_{2})(t) \neq 1$ and 
$\beta_{1}(\pi(t')) = \beta_{1}(t') \neq 1$, by \eqref{eq:comm-rule} the inequalities \eqref{eq:ineq-tecnico0} 
in  Lemma \ref{lem:tecnico0} hold for ${\mathtt y} =\pi(r)$ and ${\mathtt z} = \pi (s)$.
Hence, $\pi(\Oc_r^{\G^F})$ is of type C. 
\epf

\section{Mixed classes in Chevalley groups}\label{sec:chevalley}

In this Section, $F=\Fr_q$ and $\Phi$ is not of type $A_n$. 
We keep the notation introduced in 
Section \ref{sec:mixed}.

By Lemma \ref{lem:specifico} and Proposition \ref{prop:center} it remains to deal with the cases
$\x_s^2\in  Z(\G^F)$ if $w_0=-1$ and $q$ is odd, and
$\x_s^{q-1}\in Z(\G^F)$ if $w_0\neq-1$.  In the latter situation $q>2$ because 
$\x_s\not\in Z(\G^F)$. In other words, we have to deal with the following cases:
\begin{enumerate}[leftmargin=*]
\item $\Phi$ of type $B_\ell$, $\ell\geq 3$; $C_\ell$, $\ell\geq 2$; $D_{2n}$, $n\geq 2$;  
$E_7$; $E_8$;  $F_4$; $G_2$; $q$ is odd and $\x_s^2\in Z(\G^F)$ (every $z \in Z(\G^F)$ has order $\leq 2$) 
\cite[Table 24.2]{MT}; 
\item $\Phi$ of type $D_{2n+1}$, $n\geq 2$; $E_6$ and $\x_s^{q-1}\in Z(\G^F)$ (here $|Z(\G^F)|\leq 4$);
\end{enumerate}
and  $v$ satisfies either \eqref{cond:Jv=1} or \eqref{cond:Jv=2}.

\subsection{The case $\x_s\in \T^F$}
\label{subsec:split}

We recall from Subsection \ref{subsec:semisimple-comm} that 
if $\x_s^{q-1}=1$, then $t\in\T^F$ and we may assume $\x_s=t\in \T^F$. 
We now deal with this situation.  

\begin{lema}\label{lem:alpha-injective}Let $\x_s\in \T^F - Z(\G^F)$.
By  Subsection \ref{subsec:alpha} there always exists $\alpha\in\Delta$ such that 
$s_\alpha(\x_s)\neq \x_s$. 	 Assume that we are not in the situation: 
\begin{align}\label{eq:special-situation}
&\G \text{ is of type }B_\ell, & &q \text{ is odd,} && t \text{ satisfies } s_{\alpha_j}(t) \neq t \iff j = \ell.
\end{align}
Then the root $\alpha$  can be chosen so that 
$s_\alpha(\x_s)\not\in \x_s Z(\G^F)$. 
\end{lema}
\pf  The statement follows when $\G$ is of type $E_8, F_4, G_2$, 
because $Z(\G)=Z(\G^F)$ is trivial in these cases. In the general case, assume that 
$s_\alpha(\x_s)\overset{\star}{=}z\x_s$ for some $z\in Z(\G^F) - 1$. 
Applying $s_\alpha$, we get $z^2 = 1$. Thus $\G$ is not of type $E_6$  
(here $Z(\G) \simeq \Z/3$) and $q$ should be odd. By $\star$, we have
\begin{align}
\omega_i(\x_s) &= \omega_i(zs_\alpha(\x_s)) = \omega_i(z) s_\alpha(\omega_i)(\x_s) , &\forall i &\in \I_\ell.
\end{align}
Say $\alpha = \alpha_j$, $j \in \I_\ell$. If $i\neq j$ then $s_\alpha(\omega_i) = \omega_i$, hence 
$\omega_i(z) = 1$. Now such $z$ exists only if we are in \eqref{eq:special-situation}, see
the description of $Z(\G^F)$ in Table \ref{tab:center}.
\epf

\begin{table}[ht]
\caption{Center of $\G_{sc}^F$, $q$ odd; $\zeta, \omega \in \F_q^{\times}$, $\vert \zeta\vert = 4$,
 $\vert \omega\vert  = 3$}
\label{tab:center}
\begin{tabular}{|c|c|}
\hline  type  &  $Z(\G_{sc}^F)$   \\
\hline
$B_\ell$ & $\langle \alpha^{\vee}_\ell(-1)\rangle$\\
\hline
  $C_\ell$, $\ell > 2$ & 
  \footnotesize{$\left\langle \displaystyle\prod_{i\text{ odd}}\alpha^{\vee}_i(-1)\right\rangle$}\\
  \hline
  $D_{\ell}$,  $\ell\in 2\Z$  & 
  \footnotesize{$\left\langle \displaystyle\prod_{i \text{ odd}} \alpha^{\vee}_i(-1),\,
  \alpha^{\vee}_{\ell-1}(-1)\alpha^{\vee}_{\ell}(-1)\right\rangle$}\\
  \hline
$D_{\ell}$,  $\ell\in 2\Z + 1$,  $q \equiv 1 (4)$
 & \footnotesize{$\left\langle \displaystyle \prod_{i\text{ odd}
 \leq \ell - 2}\alpha^{\vee}_i (-1)\alpha^{\vee}_{\ell -1}(\zeta)\alpha^{\vee}_{\ell}(\zeta^3)\right\rangle$}
\\
\hline
 $D_{\ell}$,  $\ell\in 2\Z + 1$, $q \equiv 3 (4)$
& \footnotesize{$\left\langle \displaystyle \alpha^{\vee}_{\ell -1}(-1)\alpha^{\vee}_{\ell}(-1)\right\rangle$}
\\
\hline
$E_6$, $q \equiv 1(3)$ &
$\langle \alpha^{\vee}_1(\omega)\alpha^{\vee}_3(\omega^2)
\alpha^{\vee}_5(\omega)\alpha^{\vee}_6(\omega^2)\rangle$\\
\hline
$E_6$, $q \not\equiv 1(3)$ &1\\
\hline
$E_7$ & $\langle  \alpha^{\vee}_2(-1)\alpha^{\vee}_5(-1)\alpha^{\vee}_7(-1) \rangle$
\\ 
\hline 
\end{tabular}
\end{table}
 
\begin{lema}\label{lem:toro-split-C} Let $\x_s\in\T^F -Z(\G^F)$, let $\alpha\in\Delta$ be as in Subsection 
\ref{subsec:alpha}, and let $\Pa$ be the minimal standard $F$-stable parabolic subgroup with standard 
Levi complement $\Le$ associated with $\alpha$. Let $L=\Le^F$. Then $\Oc_{\x_s}^{L}$ is of type C. 
\end{lema}
\pf Let $\X_1 = \x_s \U_\alpha^F$ and $\X_2 = s_{\alpha}(\x_s) \U_\alpha^F$. 
Then $\X_1=\U_{\alpha}^F \trid \x_s=\U_\alpha^F \x_s$ and 
$\X_2=\U_{\alpha}^F \trid s_{\alpha}(\x_s) =\U_\alpha^F s_{\alpha}(\x_s)$, using \eqref{eq:comm-rule} and 
$\alpha(s_\alpha(\x_s))\neq1$. Hence, $\X_i\trid \X_j=\X_j$ for $i,j\in\I_2$.
Set $\Y = \X_1 \cup \X_2$.
Note that $\langle \Y \rangle = \langle \x_s, s_\alpha(\x_s),  \U_{\alpha}^F \rangle \subseteq L$.
Then $\X_1 \cap \X_2 = \emptyset$, $\x_s \in \X_1 = \Oc_{\x_s}^{\langle \Y \rangle}$, 
$s_\alpha(\x_s) \in \X_2 = \Oc_{s_\alpha(\x_s)}^{\langle \Y \rangle}$ and 
$\vert X_1\vert = \vert X_2\vert = q > 2$. Taking
$r = \x_s$ and $s = s_\alpha(\x_s)x_\alpha(1)$, we conclude that $\Y$ is of type C.
\epf

\begin{lema}\label{lem:mixed-split}
If $\x_s\in\T^F -Z(\G^F)$ and 
we are not in the situation \eqref{eq:special-situation}, then $\Oc_x^{\Gb}$ is not kthulhu.
\end{lema}
\pf In this case we take $g,\w=1$, so $\Ha_i=\G_i$ for all $i\in\I_r$ 
(cf. Subsections \ref{subsec:semisimple}, \ref{sec:unipotent_in_K}). 
Recall that $v$ satisfies either \eqref{cond:Jv=1} or \eqref{cond:Jv=2}.
Since $\w=1$, the class in $\PSU_n(q)$ does not occur
because $\PSU_n(q)$ does not occur in the decomposition of $[\Ha,\Ha]^{F}$, see Subsection \ref{sec:unipotent_in_K}.
Also, we may assume that $q>2$, as the classes in groups over $\F_2$ occur only for 
$\x_s=1$
Thus, each component in $v=\x_u$ can be chosen to lie in a 
$\U^F_{\beta}$ for some $\beta\in\Phi^+$ such that $\beta(\x_s)=1$.  
By the hypothesis on $\x_s$, there is $\alpha\in\Delta$ as in 
Lemma \ref{lem:alpha-injective}. Let $\Pa$ be the standard $F$-stable parabolic subgroup with 
standard Levi complement $\Le$ associated with $\alpha$, and let $P=\Pa^F$, $L=\Le^F$.  
Then for $\pi_L\colon P\to L$ we have that $\pi_L(\Oc^{P}_{\x_s\x_u})=\Oc_{\x_s}^L$  is of type C by 
Lemma \ref{lem:toro-split-C}. 
Let $\Y$ as in the proof of Lemma \ref{lem:toro-split-C} and let 
\begin{align*}
\X_1'&:=\pi_L^{-1}(\X_1)\cap\Oc^P_{\x}=\x_s\U^F\cap \Oc_{\x}^P,\\
\X_2'&:=\pi_L^{-1}(\X_2)\cap\Oc_{\x}^P=s_\alpha(\x_s)\U^F\cap \Oc_{\x}^P,\\
\Y'&:=\pi_L^{-1}(\Y)\cap\Oc_{\x}^P=\X_1'\cup \X_2'.
\end{align*}
Let $\pi\colon \G^F\to\Gb$ be the isogeny and let $\ta \neq \tb \in \Y'$ such that $\pi(\ta) = \pi(\tb)$, 
i.e. there is $z \in Z(\G^F)$, $z\neq 1$, such
$\ta = z\tb$. Hence either $\ta = \x_su $ and $\tb = s_\alpha(\x_s) v$ for some $u,v\in\U^F$, or vice versa.
In any case, $\x_s  = zs_\alpha(\x_s)$, impossible by our choice of $\alpha$. 
Hence, the restriction of $\pi$ to $\Y'$ is injective and $\Oc_x^{\Gb}$ is of type C.\epf

\bigskip

In the next lemma we deal with the case in which $\x_s\in\T^F -Z(\G^F)$ is as in \eqref{eq:special-situation}.  
Then $\x_s$ has the form
\begin{align}\label{eq:special-situation-explicit}
 & \x_{s}=\left(\prod_{i\in \I_{\ell -1}} \alpha^{\vee}_i ((-1)^{i})\right) \alpha^{\vee}_{\ell} (\eta),&
\text{where }  \eta^2&= (-1)^{\ell}.
\end{align}
Notice that if $\ell$ is odd, then such an $\x_s$ belongs to $\G^F$ if and only if $q \equiv 1 (4)$.

\smallbreak
We shall also need a realization of the symplectic group $\Sp_{2n}(\kk)$ for $q$ odd. We choose it to be the subgroup of $\GL_{2n}(\kk)$ 
consisting of matrices preserving the skew-symmetric bilinear form with associated matrix \begin{align}\label{eq:forma}
\widetilde {\Jf}_{2n}=\left(\begin{smallmatrix}
 0&\Jf_n\\                                                                                                                                                                                                                               
-\Jf_n&0                                                                                                                                                                                                                              
\end{smallmatrix}\right), \text{ where }
\Jf_n=\left(\begin{smallmatrix}
&&1\\ 
&\iddots&\\
1&&\end{smallmatrix}\right).\end{align}
When clear from the context, we shall omit the index $n$ from $\Jf_n$ and $\widetilde {\Jf}_{2n}$.

\begin{lema}\label{lem:special}Assume $\Phi$ is of type $B_\ell$ and $q$ is odd. 
 Let $\x=\x_s\x_u$ with $\x_u\neq1$ and  $\x_s$ is as in \eqref{eq:special-situation}. Then $\Oc_x^{\Gb}$ is not kthulhu.
 \end{lema}
\pf Here $\Ha=C_{\G}(\x_s)=C_{\G}(t)$ is simple with $\Phi_t$ of type 
$D_\ell$. As $\ell>2$ is not compatible with \eqref{cond:Jv=1} nor \eqref{cond:Jv=2}, we have
$\ell=2$, $\G\simeq \Sp_4(\kk)$ and $\Ha\simeq \SL_2(\kk)\times\SL_2(\kk)$ is semisimple. 
Our assumptions on $\x_u=v$ give either: 
 \begin{enumerate}[leftmargin=*]
 \item\label{item:qodd} $q=9$ or not a square and only one component of 
 $x_u$ in $\SL_2(\kk)\times\SL_2(\kk)$ is non-trivial; or else
 \item\label{item:q=3} $q=3$ and the two components of $x_u$ in $\SL_2(\kk)\times\SL_2(\kk)$ are non-trivial.
 
 \end{enumerate}
To deal with \eqref{item:qodd}, we work in $\G^F=\Sp_4(9)$ for simplicity. In this case, 
$\B^F$ is the subgroup of $\Sp_4(9)$ consisting of upper triangular matrices. Let $\pi\colon \G^F\to \Gb$. 
We take
\begin{align*}
 \x= \x_sx_{\alpha_{2}}(a)&=\left(\begin{smallmatrix}
            -1&0&0&0\\
            0&1&a&0\\
            0& 0& 1&0\\
            0& 0& 0 &-1\\
           \end{smallmatrix}\right),\quad a\in\F_9^\times.
\end{align*}
Let 
\begin{align*}
 \sigma&=\left(\begin{smallmatrix}
            0&1&0&0\\
            1&0&0&0\\
            0& 0& 0&1\\
            0& 0& 1 &0
           \end{smallmatrix}\right),&r=  \sigma\trid \x=\left(\begin{smallmatrix}
            1&0&0&a\\
            0&-1&0&0\\
            0& 0& -1&0\\
            0& 0& 0 &1
           \end{smallmatrix}\right),\\   
  x_{\alpha_1}(1)&=\left(\begin{smallmatrix}
            1&1&0&0\\
            0&1&0&0\\
            0& 0& 1&-1\\
            0& 0& 0 &1
           \end{smallmatrix}\right),&s= x_{\alpha_1}(1)\trid \x=\left(\begin{smallmatrix}
            -1&2&a&a\\
            0&1&a&a\\
            0& 0& 1&2\\
            0& 0& 0 &-1
           \end{smallmatrix}\right).         
\end{align*}

A direct computation shows that $(\pi(r)\pi(s))^2\neq (\pi(s)\pi(r))^2$. Also, 
$r,s\in\B^F=\T^F\U^F$ and, for $V=\langle \U_{\beta}^F~|~\beta\in\Phi^+-\alpha_2\rangle$ 
by Chevalley's commutator formula we have 
$\B^F\trid r\in \T^F V$ whereas $\B^F\trid s\in \T^F x_{\alpha_2}(\F_9^\times) V$.  
Hence, $\Oc_{\pi(r)}^{\langle \pi(r),\pi(s)\rangle}\neq \Oc_{\pi(s)}^{\langle \pi(r),\pi(s)\rangle}$ 
and $\Oc_x^{\Gb}$ is of type D.

\smallbreak
There are $2$ classes in case \eqref{item:q=3}, namely the classes represented by  
\begin{align*} \x_{a}&=\left(\begin{smallmatrix}
            1&0&0&1\\
            0&2&a&0\\
            0& 0& 2&0\\
            0& 0& 0 &1\\
           \end{smallmatrix}\right), \quad a\in\F_3^\times.
           \end{align*}
They lie in the same orbit for the conjugation action of the group of diagonal matrices ${\rm diag}(\xi,\eta,\lambda \xi^{-1},\lambda\eta^{-1})$, 
for $\xi,\,\eta\in\F_9^\times$, $\lambda\in\F_3^\times$ on  $\Sp_4(q)$. Hence $\Oc_{\x_{1}}\simeq \Oc_{\x_{2}}$. Let $r:=\x_{1}$.            
With notation as in \eqref{item:qodd}, let
\begin{align*}
r_1=  \sigma\trid r=\left(\begin{smallmatrix}
            2&0&0&1\\ 
            0&1&1&0\\
            0& 0& 1&0\\
            0& 0& 0 &2
           \end{smallmatrix}\right),& &s= x_{\alpha_1}(1)\trid r_1=\left(\begin{smallmatrix}
            2&2&1&2\\
            0&1&1&1\\
            0& 0& 1&2\\
            0& 0& 0 &2
           \end{smallmatrix}\right).
  \end{align*}
Again, $(\pi(r)\pi(s))^2\neq (\pi(s)\pi(r))^2$. Also, $r,s\in\B^F$ and 
$\B^F\trid r\in \T^F x_{\alpha_2}(2) V$ whereas $\B^F\trid s\in \T^F x_{\alpha_2}(1)V$.  Hence, 
$\Oc_{\pi(r)}^{\langle \pi(r),\pi(s)\rangle}\neq \Oc_{\pi(s)}^{\langle \pi(r),\pi(s)\rangle}$ so this class is  of type D.    
\epf
 
\begin{prop}\label{prop:e8f4g2-qeven}
Assume $\Phi$ and $q$ fall into one of the following cases: 
\begin{enumerate}
 \item $E_8$, $F_4$, $G_2$, $q$ arbitrary;
 \item $E_6$ and  
$q\not\equiv 1(3)$; 
\item $D_{2m+1}$ and $q$ is even.
\end{enumerate}
Then $\Oc_x^{\Gb}$ is not kthulhu. 
\end{prop}

\pf In these cases, we have that $Z(\G^F)=1$, see \cite[Table 24.2]{MT}. Thus,
by Subsection \ref{subsec:semisimple-comm}, the condition 
$\x_{s}^{q-1}\in Z(\G^F)$ reads $\x_{s} \in \T^{F} $. This is dealt with in Lemma \ref{lem:mixed-split}.
\epf

\subsection{The remaining cases}\label{subsec:remaining}
According to the list at the beginning of Section \ref{sec:chevalley}, 
we are left with the following cases:
\begin{setting}\label{set:chev}$\x_s^{q-1}=z\in Z(\G^F)-\{1\}$; either \eqref{cond:Jv=1} or \eqref{cond:Jv=2} holds for $v$; 
and $\Phi$, $q$ and $\x_s$ satisfy: 
\begin{enumerate}[leftmargin=*]
\item $C_\ell$, $\ell\geq 2$;  $\x_s^2 = z$,  so $q\equiv 3(4)$;
\item $B_\ell$, $\ell\geq 3$;  $\x_s^2 = z$, so $q\equiv 3(4)$;
\item $D_{2n}$, $n\geq 2$;  $\x_s^2 = z$,  so $q\equiv 3(4)$;
\item $E_7$;  $\x_s^2 = z$,  so $q\equiv 3(4)$;
\item $D_{2n+1}$,  $n\geq 2$, $q$ is odd;
\item $E_6$, $q\equiv 1(3)$. 
\end{enumerate}
\end{setting}

Assume that we are in either of the situations $(1), \ldots, (4)$. By Table \ref{tab:center}, we must have that
$q\equiv 3(4)$, since otherwise $\x_s^{q-1}=1$. Analogously, $q\equiv 1(3)$ if $\Phi$ is of type $E_{6}$.

\smallbreak
In what follows, we provide several technical remarks and lemmata to deal with these cases.


\begin{obs}\label{obs:A1}If $\Phi_{t}^{j}$ is of type $A_{1}$ with base $\Delta_j=\{\alpha_{j_1}\}$ and $\Ad(\w)\G_j=\G_j$, 
then,  by replacing $w$ 
by $s_{\alpha_{j_1}}w$, we can always ensure that $w$ acts trivially on $\Phi_{t}^j$.
Observe that this replacement does not affect the action of $w$  
on other subsets $\Delta_l$, for $l\in\J_{v}-\{j\}$.
\end{obs}

\begin{lema}\label{lem:tecnico}Assume that \eqref{cond:Jv=2} holds for $v$,
that $\Ad(\w)^{-1}(t)=t^q=zt\neq t$ with $z\in Z(\G^{F})$, and that $w$ acts trivially on $\Phi_{t}^{i}$ for all
$i \in \J_{v}$.
If there exists $\alpha\in\Delta$ such that:
\begin{enumerate}
\item\label{item:commu} $s_\alpha w=w s_\alpha$;
\item\label{item:nonPsi} $s_\alpha(\beta)\not\in\Phi_t$ for some $\beta \in (\Phi_t^{i_{0}})^{+}$ 
and some $i_{0} \in \J_{v}$ with $\beta \neq \alpha_{0}$;
\end{enumerate} 
then $\Oc_x^{\Gb}$ is of type C.
\end{lema}

\pf Note that in this situation $q\neq 2$. As $|\J_{v}|>1$, we necessarily have $q=3$
and $\Phi_{t}^{i}$ is of type $A_{1}$ for all $i \in \J_{v}$. 

We set for simplicity $i_0=1$. 
Then  $v_{1}\in\Ha_{1}=\G_{1}$, $G_{1}=\G_{1}^{\Fr_q}$, and since $w$ acts trivially on $\Phi_{t}^{1}$
we may assume that $v_{1}=x_{\beta}(\xi)$ for some $\xi\in\F_q^\times$.  

Let  $\Delta_1\subset\Delta$ be the base of $\Phi_{t}^{1}$; let $\widetilde \Pi=\Delta_1\cup\{\alpha\}$; let $\Pa$ be the $F_w$-stable standard parabolic subgroup of $\G$ 
associated with $\widetilde {\Pi}$;  $\Vu$ be its unipotent radical and $\Le$ be the $F_w$-stable standard Levi subgroup. 
Up to conjugation in $\Ha^{F_w}$, we have that $v=x_{\beta}(\xi)v'$ for some $v'\in \Vu$.
Let
\begin{align*}
r:=tv=tx_{\beta}(\xi)v'\in tx_{\beta}(\xi)\Vu.
\end{align*} 
Condition \eqref{item:commu}  
and \cite[Proposition 25.3]{MT} applied to $F_w$ ensure that there exists a representative $\sa\in N_{\G}(\T)\cap\G^{F_w}$.
We set 
\begin{align*}
s:=\sa\trid r=s_\alpha(t)x_{s_\alpha\beta}(\zeta)v''\in s_\alpha(t)x_{s_\alpha\beta}(\zeta)\Vu
\end{align*} 
for some $\zeta\in \kk^\times$ and some $v''\in\Vu$, since $\sa \in \Le$. 
Condition \eqref{item:nonPsi} 
ensures that $\alpha\not\in\Phi_t$ so $s_\alpha\gamma\in\Phi^+$ for every $\gamma\in\Phi^{+}_t$.  
Hence, $\langle r,\,s\rangle\subset\T\U$. Thus:
\begin{align*}
\Oc_r^{\langle r,s\rangle}\subset \Oc_r^{\T\U}\subset t\U,& & \Oc_s^{\langle r,s\rangle}\subset 
\Oc_s^{\T\U}\subset s_\alpha(t)\U.
\end{align*}
Condition \eqref{item:nonPsi}  also gives $s_\alpha\Phi_t\not=\Phi_t$ whence $\sa\not\in N_{\G}(\Ha)$.
In particular, $\sa\not\in\Ha$, so $s_\alpha(t)\neq t$ and therefore 
$\Oc_r^{\langle r,s\rangle}\cap\Oc_s^{\langle r,s\rangle}=\emptyset$. 
Since $\sa\not\in N_{\G}(\Ha)$, we have $s_\alpha(t)\not\in Z(\G) t$. 
Thus, the restriction of the isogeny  
to $\Oc_{r}^{\langle r,s\rangle}\coprod  \Oc_s^{\langle r,s\rangle}$ 
is injective and $\pi\left(\Oc_r^{\langle r,s\rangle}\right)
\cap\pi\left(\Oc_s^{\langle r,s\rangle}\right)=\emptyset$. 
As $\alpha(\pi(t)) = \alpha(t) \neq 1$ and 
$\beta(\pi(t')) = \beta(t') = \beta(s_\alpha(t)) =  (s_\alpha\beta)(t) \neq 1$, 
by \eqref{eq:comm-rule} the inequalities in 
Lemma \ref{lem:tecnico0} hold for ${\mathtt y} = \pi(r)$ and ${\mathtt z} = \pi (s)$.
Hence, $\pi\left( \Oc_r^{\G^F} \right)$ is of type C, and consequently, 
$\Oc_x^{\Gb}$ is also of type C, by Remark \ref{obs:rackiso-conjg}.

\epf


\begin{obs}\label{obs:representativev}
If \eqref{cond:Jv=1} holds for $\Oc_x^{\Gb}$, 
$\Ad(\w)^{-1}(t)=t^q=zt\neq t$ with $z\in Z(\G^{F})$, and $w$ acts trivially on $\Phi_{t}^{1}$, then 
$q\neq 2$. By the discussion at the beginning of Subsection \ref{subsec:weyl} we can always 
make sure that 
$v=v_{1}=x_{\beta}(\xi)$ for some $\beta \in (\Phi_{t}^{1})^{+}$ and $\xi \in \F_{q}^{\times}$.
\end{obs}

\begin{lema}\label{lem:tecnico2}
Assume $\Ad(\w)^{-1}(t)=t^q=zt\neq t$ with $z\in Z(\G^{F})$, 
$v=v_{1}=x_{\beta}(\xi) \in G_{1}$ for some $\beta \in (\Phi_{t}^{1})^{+}$ and $\xi \in \F_{q}^{\times}$.  
If there exists $\alpha\in\Phi^+$ such that
$s_\alpha w=w s_\alpha$ and 
$s_\alpha(\beta)\in\Phi^+-\Phi_t$, 
then $\Oc_x^{\Gb}$ is of type C.
\end{lema}

\pf
Similar to that of Lemma \ref{lem:tecnico}. 
Since $|\J_v|=1$, the computations for $r=tx_{\beta}(\xi)$ and $s=\sa\trid r$ 
follow as before without introducing  $\Pa$ and $\Vu$. 
The condition $s_\alpha(\beta)\in\Phi^+-\Phi_t$ ensures that $s\in\T\U$.
\epf

We are now ready to state the main result of this Subsection.

\begin{prop}\label{prop:tq=tz}
Assume that we are in the Setting \ref{set:chev}. Then $\Oc_x^{\Gb}$ is not kthulhu.
\end{prop}

We prove this proposition by analysing the different groups separately.

\begin{lema}\label{lem:cn}
Assume that we are in Setting \ref{set:chev} with $\Phi$ of type $C_\ell$, so that $\G=\Sp_{2\ell}(\kk)$, $\ell\geq 2$ and $q\equiv 3(4)$. 
Then $\Oc_x^{\Gb}$ is not kthulhu.
\end{lema}

\pf Here $\T$ is the subgroup of diagonal matrices of the form
\begin{align}\label{eq:shape}
{\rm diag}(x_1,\ldots,\,x_\ell,x_\ell^{-1},\,\ldots,\,x_1^{-1}).
\end{align} 
The assumption on $\x_s$ gives $\x_s^2=-\id$ so 
$$t={\rm diag}(\xi_1,\ldots,\,\xi_\ell,\xi_\ell^{-1},\,\ldots,\,\xi_1^{-1})\text{ with }\xi_i^2=-1 \text{ for }i\in\I_\ell.$$

Since permuting eigenvalues compatibly with \eqref{eq:shape} gives a new element 
in $\Oc_t^{\G}$ and $\Oc_t^{\G}\cap \G^F=\Oc_{\x_s}^{\G^F}$, 
we can reorder the eigenvalues and assume that $t={\rm diag}(\xi\id_\ell,-\xi\id_\ell)$ for $\xi^2=-1$. 

The matrices in 
$\Ha$ have a diagonal block form of shape  ${\rm diag}(A, \Jf\, ^t\!A^{-1}\Jf)$ for 
$A\in \GL_\ell(\kk)$, i.e., $[\Ha,\Ha]\simeq\SL_\ell(\kk)$ is simple. 
Since $v=v_1$ occurs in Table \ref{tab:unip-kthulhu}, this is possible only if $\ell=2$.

In this case $t$ is conjugate to ${S}=\left(\begin{smallmatrix}
0&1&0&0\\
-1&0&0&0\\
0&0&0&-1\\
0&0&1&0
\end{smallmatrix}\right)$.

By a direct computation of unipotent matrices in $C_{\G}(S)$ we see that 
$\Oc_\x^{\G^F}$ is represented by a matrix of this form. 
\begin{align}
r:=S \left(\begin{smallmatrix}
   1&0&0&b\\
   0&1&b&0\\
   0&0&1&0\\
   0&0&0&1
       \end{smallmatrix}\right)=\left(\begin{smallmatrix}
  0 &1&b&0\\
   -1&0&0&-b\\
   0&0&0&-1\\
   0&0&1&0
       \end{smallmatrix}\right),\quad\text{ for some $b\in \F_q^\times$}.
\end{align}

Assume first $q\neq 3,7$. Let $\Pa$ be the standard $F$-stable parabolic subgroup with 
standard Levi factor $\Le$ associated with $\alpha_{1}$, and let $\mathbf{P}=\pi(\Pa^F)$, 
$\mathbf{L}=\pi(\Le^F)$, with $\pi_{L}: \mathbf{P} \to \mathbf{L}$ the corresponding projection.
Then $\Oc_{x}^{\Gb}$ contains the subrack $\Oc_{x}^{\mathbf{P}}$, which in turn projects onto
$\Oc_{\pi_{L}(x)}^{\mathbf{L}}$. The latter is isomorphic as rack to $\Oc^{\PGL_{2}(q)}_{y}$, with 
$y$ the class of $\left(\begin{smallmatrix}
  0 &1\\
   -1&0
       \end{smallmatrix}\right)$.
Since by \cite[Theorem 1.1]{ACG-III} this class is not kthulhu, the same holds
for $\Oc_{x}^{\Gb}$.

Assume now that $q=3$ or $7$.
Let $u:= \left(\begin{smallmatrix}
   0&1&b&0\\
   1&1&b&0\\
   0&0&0&1\\
   0&0&1&-1
       \end{smallmatrix}\right)$ and $s:=u\trid r= \left(\begin{smallmatrix}
   1&-1&-2b&0\\
   2&-1&-2b&0\\
   0&0&1&1\\
   0&0&-2&-1
       \end{smallmatrix}\right)$. 
A direct computation shows that $u\in \G^{F}$ and 
$(rs)^{2}\neq (sr)^{2}$, since  
$(rs)^{2}:= \left(\begin{smallmatrix}
   5&-3&-7b&b\\
   -3&2&14b&5b\\
   0&0&5&3\\
   0&0&3&2
       \end{smallmatrix}\right)$ and $(sr)^{2}= \left(\begin{smallmatrix}
   2&3&b&11b\\
   3&5&4b&13b\\
   0&0&2&-3\\
   0&0&-3&5
       \end{smallmatrix}\right)$. Moreover, this inequality also holds in the projection to $\Gb$ as
 $(rs)^{2}(sr)^{-2}:= \left(\begin{smallmatrix}
  1&0&2b&b\\
   0&1&2b&2b\\
   0&0&1&0\\
   0&0&0&1
       \end{smallmatrix}\right)$ if $q=3$ and 
$(rs)^{2}(sr)^{-2}:= \left(\begin{smallmatrix}
   -1&0&-b&-b\\
   0&-1&4b&-b\\
   0&0&-1&0\\
   0&0&0&-1
       \end{smallmatrix}\right)$ if $q=7$.   
In addition, if $q=3$ we have 
\begin{align*} 
\Oc^{\langle r,s\rangle}_{r} & \subseteq \left\{ \left(\begin{smallmatrix}
A&B \\
 0&\Jf_{2}A^{-t}\Jf_{2}
\end{smallmatrix}\right) \text{ with }A = \pm
\left(\begin{smallmatrix}
       0&1\\
       -1&0
      \end{smallmatrix}\right) 
\right\},\\
 \Oc^{\langle r,s\rangle}_{s} & \subseteq \left\{\left(\begin{smallmatrix}
A&B \\
 0&\Jf_{2}A^{-t}\Jf_{2}
\end{smallmatrix}\right) \text{ with } A = \pm
\left(\begin{smallmatrix}
       1&2\\
       2&2
      \end{smallmatrix}\right)\right\};
       \end{align*} 
whereas if $q=7$ we have 
\begin{align*} 
\Oc^{\langle r,s\rangle}_{r} & \subseteq \left\{ \left(\begin{smallmatrix}
A&B \\
 0&\Jf_{2}A^{-t}\Jf_{2}
\end{smallmatrix}\right) \text{ with }A \in \left\{
\pm \left(\begin{smallmatrix}
       0&1\\
       -1&0
      \end{smallmatrix}\right), 
\pm \left(\begin{smallmatrix}
       4&2\\
       2&3
      \end{smallmatrix}\right) 
      \right\} \right\},\\
 \Oc^{\langle r,s\rangle}_{s} & \subseteq \left\{\left(\begin{smallmatrix}
A&B \\
 0&\Jf_{2}A^{-t}\Jf_{2}
\end{smallmatrix}\right) \text{ with } A \in \left\{
\pm \left(\begin{smallmatrix}
       1&-1\\
       2&-1
      \end{smallmatrix}\right), 
\pm \left(\begin{smallmatrix}
       -1&-2\\
       1&1
      \end{smallmatrix}\right)
      \right\} \right\}.
       \end{align*} 
       In both cases $\pi(\Oc^{\langle r,s\rangle}_{r})\neq \pi(\Oc^{\langle r,s\rangle}_{s})$, whence
 $\Oc^{\Gb}_{x}$ is of type D.
       \epf

\begin{lema}\label{lem:bn}
Assume that we are in Setting \ref{set:chev} with $\Phi$ of type $B_\ell$, $\ell\geq 3$, so that $q\equiv 3(4)$. Then $\Oc_{x}^{\Gb}$ is not kthulhu.
\end{lema}
\pf Let $\widetilde {\pi}\colon \G\to \G_{\ad}=\SO_{2\ell+1}(\kk)$ be the natural projection. 
We realize $\SO_{2\ell+1}(\kk)$ as the subgroup of $\SL_{2\ell+1}(\kk)$ 
consisting of matrices preserving the symmetric bilinear form with 
associated matrix ${\Jf}_{2\ell+1}$ as in \eqref{eq:forma}, 
and $\widetilde {\pi}(\T)$ to be the group of diagonal matrices of the form
\begin{align}\label{eq:shapebn}
{\rm diag}(x_1,\ldots,\,x_\ell,1,x_\ell^{-1},\,\ldots,\,x_1^{-1}).
\end{align} 
Since $\widetilde {\pi}(t)^2=1$, its eigenvalues are $\pm1$.
Reordering eigenvalues (i.e., acting via the Weyl group) we assume 
that $\widetilde {\pi}(t)={\rm diag}(-\id_k,\id_{2\ell-2k+1},-\id_k)$, i.e., 
\begin{align*}
t=\left(\prod_{i=1}^{k}\alpha_i^\vee((-1)^i)\right)
\left(\prod_{j=k+1}^{ \ell-1}\alpha_j^\vee((-1)^k)\right)\alpha_\ell^\vee(\xi), \text{with $\xi^2=(-1)^k$}.
\end{align*} 
The condition $t^{2}\in Z(\G^F)-1$ forces $k$ to be odd. 
If $k\geq 3$ the root system  $\Phi_t$ is the union of two orthogonal subsystems, 
with base $\Delta_1=\{\alpha_i,\,i\in\I_{0,k-1}\}$ of type $D_k$ and 
$\Delta_2=\{\alpha_j,\,j\in\I_{\ell-k+1,\ell}\}$ of type $B_{\ell-k}$, whereas if $k=1$, 
$\Phi_t$ has base $\Delta_2=\{\alpha_j,\,j\in\I_{2,\ell}\}$ of type $B_{\ell-1}$. 
Let $\gamma=\alpha_1+\cdots+\alpha_\ell=\varepsilon_1$. A direct calculation 
shows that we can take $\w=s_{\gamma}$. 
Indeed, $s_\gamma(\alpha_1)=-\alpha_0$ and $s_\gamma(t)=\alpha_\ell^\vee(-1)t$. So from 
Table \ref{tab:center},  $s_\gamma(t)=zt=t^{q-1}$ as desired. 
Observe that when $k=3$ we have $D_3=A_3$. 
Then \eqref{cond:Jv=2} holds with $\J_v=\{2\}$ and $\ell-k\in\{1,2\}$. 
Remark \ref{obs:representativev} applies with  $\beta=\alpha_\ell$ if $\ell=k+1$ and $\beta=\alpha_{\ell-1}$ if $\ell=k+2$. 
Lemma \ref{lem:tecnico} applies with $\alpha=\alpha_{\ell-1}$, or $\alpha=\alpha_{\ell-2}$, 
respectively with the exception of the case $\ell=3,\,k=1$.
In this case, $t=\alpha_{1}^{\vee}(-1)\alpha_{2}^{\vee}(-1)\alpha_{3}^{\vee}(\xi)$, 
$w(\alpha_{1})=-\alpha_{0}$, $w(\alpha_{2})=\alpha_{2}$, $w(\alpha_{3})=\alpha_{3}$,
$\Phi_{t}$ is
of type $B_{2}$ and  $tv = \alpha_{1}^{\vee}(-1)\alpha_{2}^{\vee}(-1)\alpha_{3}^{\vee}(\xi)x_{\alpha_{2}}(\zeta)$ for 
some $\zeta\in \F_{q}^{\times}$.
Let us set $\alpha_1=\beta_1$, $\alpha_2=\beta_2$ and $\beta_3=-\alpha_0$. 
This is a base for a root system of type $A_3$. Let $\widetilde {\U}=\langle \U_{\beta_i},\,i\in\I_3\rangle$, 
and let ${\mathbb K}=\langle  \widetilde {\U}, \w_0(\widetilde {\U})\rangle$ be the corresponding $F_w$-stable 
algebraic subgroup of $\G$. Let $\Pa$ be the $F_w$-stable standard parabolic subgroup of ${\mathbb K}$ 
associated with $\{\beta_2\}$,  let $\Le$ be the corresponding standard Levi subgroup and let $\Vu$ be 
the unipotent radical. Thus
\begin{align*}
tv=\beta_1^\vee(\xi)\beta_3^\vee(\xi^{-1})x_{\beta_2}(\zeta)\in t\Vu.
\end{align*} 
By the discussion in Subsection \ref{subsec:F} the element $tv$ is conjugate in $[\Ha,\Ha]^{F_w}$ to an 
element $r$ of the form:  
\begin{align*}
r=tx_{\alpha_2+2\alpha_3}(\zeta')=
\beta_1^\vee(\xi)\beta_3^\vee(\xi^{-1})x_{\beta_1+\beta_2+\beta_3}(\zeta')\in t \Vu.
\end{align*} 

Let $y=x_{\beta_1}(\eta)x_{\beta_3}(\eta')\in \widetilde {\U}^{F_w}$, with $\eta\eta'\neq0$. The existence of 
such an element is guaranteed by \cite[Proposition 23.8]{MT}. Let
$s:=y\trid tv$. 
Then \begin{align*}s\in\beta_1^\vee(\xi)\beta_3^\vee(\xi^{-1})x_{\beta_2}(\zeta)x_{\beta_1}(-2\eta)
x_{\beta_3}(-2\eta')
\U_{\beta_1+\beta_2}\U_{\beta_2+\beta_3}\U_{\beta_1+\beta_2+\beta_3}\end{align*} lies in $t x_{\beta_2}(\zeta)\Vu$ and its semisimple part equals $y\trid t=tx_{\beta_1}(-2\eta)
x_{\beta_3}(-2\eta')$.
Thus, $\langle r,s\rangle \subset \langle t,\widetilde {\U}\rangle$. 
Since $t$ commutes with $\U_{\beta_2}$, we have the inclusions
$\Oc_{r}^{\langle r,s\rangle} \subseteq t\Vu$
and $\Oc_{s}^{\langle r,s\rangle} \subseteq tx_{\beta_{2}}(\zeta)\Vu$.
As a consequence,
$\Oc_{r}^{\langle r,s\rangle}\cap \Oc_{s}^{\langle r,s\rangle} = \emptyset$.  
Moreover, $t$ does not commute with $\U_{\beta_1}$, hence it does not commute with $y\trid t$ and so $rs\neq sr$ by Remark \ref{obs:subgroup}.
%
To finish
the proof we estimate $|\Oc_{r}^{\langle r,s\rangle}|,|\Oc_{s}^{\langle r,s\rangle}|$.
A direct computation shows that
\begin{align*}
&t\trid s\in\beta_1^\vee(\xi)\beta_3^\vee(\xi^{-1})x_{\beta_2}(\zeta)x_{\beta_1}(2\eta)
x_{\beta_3}(2\eta')
\U_{\beta_1+\beta_2}\U_{\beta_2+\beta_3}\U_{\beta_1+\beta_2+\beta_3},\\
&(y\trid t)\trid(t\trid s)\in\beta_1^\vee(\xi)\beta_3^\vee(\xi^{-1})x_{\beta_2}(\zeta)x_{\beta_1}(-6\eta)
x_{\beta_3}(-6\eta')
\U_{\beta_1+\beta_2}\U_{\beta_2+\beta_3}\U_{\beta_1+\beta_2+\beta_3},\\
&(y\trid t)\trid r\in \beta_1^\vee(\xi)\beta_3^\vee(\xi^{-1})x_{\beta_1}(-4\eta)x_{\beta_3}(-4\eta')
\U_{\beta_1+\beta_2+\beta_3},\\
&t\trid((y\trid t)\trid r)\in \beta_1^\vee(\xi)\beta_3^\vee(\xi^{-1})x_{\beta_1}(4\eta)x_{\beta_3}(4\eta')
\U_{\beta_1+\beta_2+\beta_3}.
\end{align*}
Therefore $|\pi(\Oc_r^{\langle r,s\rangle})|>2$, $|\pi(\Oc_s^{\langle r,s\rangle})|>2$, 
so $\pi(\Oc_r^{\G^F})$ is of type C. 
By Remark \ref{obs:rackiso-conjg} the class $\Oc_x^{\Gb}$ is also of type C.
\epf

\begin{lema}\label{lem:d2k}Assume that we are in Setting \ref{set:chev} with $\Phi$ of type $D_\ell$, $\ell=2n\geq4$ so that $q\equiv 3(4)$. Then $\Oc_{x}^{\Gb}$ is not kthulhu.
\end{lema}

\pf Recall that $t^{q-1}=t^2=z\in Z(\G^F)-1$. By looking at 
Table \ref{tab:center} we see that $z$ can be either 
$z_1=\prod_{i\text{ odd }}\alpha^\vee_i(-1)$, or $z_2=\alpha_{\ell-1}^\vee(-1)\alpha_{\ell}^\vee(-1)$ 
or $z_3=z_1z_2$. 
Let $\overline{\pi}\colon \G\to \SO_{2\ell}(\kk)$ be the isogeny with kernel $\langle z_2\rangle$. 
We realize $\SO_{2\ell}(\kk)$ as the group of matrices in $\SL_{2\ell}(\kk)$ preserving the form $\Jf_{2l}$ from 
\eqref{eq:forma}.  Then the group $\pibar(\T)$ is given by diagonal elements of shape \eqref{eq:shape}. 

Assume first that $z=z_1$. Then $\pibar(t)^2=-\id_{2\ell}$, so 
\begin{align*}\pibar(t)={\rm diag}(\xi_1,\ldots,\,\xi_\ell,-\xi_\ell,\ldots,\,-\xi_1),\text{ with }
\xi_i^2=-1,\text{ for }i\in\I_{\ell}.\end{align*}
Up to reordering eigenvalues by acting with the Weyl group we see that $\pibar(t)$ has  the  form:
\begin{align*}\pibar(t)={\rm diag}(\xi\id_\ell,-\xi\id_\ell)\text{ or  }
{\rm diag}(-\xi,\xi\id_{\ell-1},-\xi\id_{\ell-1},\xi),\text{ with } \xi^2=-1.\end{align*}
In both cases $C_{\SO_{2\ell}(\kk)}(\pibar(t))^\circ=\pibar(\Ha)$ has semisimple part of type $A_\ell$. 
This implies that $[\Ha,\,\Ha]$ is simple so $|\J_v|=1$ but \eqref{cond:Jv=1} cannot hold. The case $z=z_3$ is dealt with in the same way.

Assume now that $z=z_2$. Then $\pibar(t)^2=\id_{2\ell}$, so 
\begin{align*}\pibar(t)={\rm diag}(\xi_1,\ldots,\,\xi_\ell,\xi_\ell,\ldots,\,\xi_1),\text{ with }\xi_i=\pm1,
\text{ for }i\in\I_{\ell}.\end{align*}
Up to reordering eigenvalues by acting with the Weyl group we see that $\pibar(t)$ has the form:
\begin{align}\label{eq:t-inv}
\pibar(t)={\rm diag}(-\id_b,\id_{2\ell-2b},-\id_b), \text{ with } b\in\I_\ell.
\end{align}

Then $\pibar(\Ha)=C_{\SO_{2\ell}(\kk)}(\pibar(t))^\circ$ is semisimple with root system of type 
$D_b\times D_{\ell-b}$, with the understanding that $D_1$ is a torus, $D_2$ is $A_1\times A_1$ 
and $D_3=A_3$. If either \eqref{cond:Jv=1} or \eqref{cond:Jv=2} holds, then necessarily 
$b=2$ and/or $\ell-b=2$. Since these two cases are obtained from one another by multiplying $\pibar(t)$ 
by $-\id_{2\ell}$ we focus on $b=2$. Up to a central element, $t=\alpha_1^\vee(-1)$. 
A direct calculation shows that we can take 
$w=s_{\varepsilon_1-\varepsilon_\ell}s_{\varepsilon_1+\varepsilon_\ell}$. 

Let us first assume $\ell>4$. Here, $[\Ha,\Ha]=\G_1\G_2\G_3$ with $\Delta_1=\{\alpha_1\}$, 
$\Delta_2=\{-\alpha_0\}$, $\Delta_3=\{\alpha_i,\,i\in\I_{3,\ell}\}$. Then $\Ha_1=\G_1\G_2$, $\Ha_2=\G_3$ 
and $G_1=\Ha_1^{F_w}\simeq\SL_2(q^2)$. Hence, $v$ satisfies \eqref{cond:Jv=1}  and  $q=3$.

If, instead, $\ell=4$, then $[\Ha,\Ha]=\G_1\G_2\G_3\G_4$ with $\Delta_1=\{\alpha_1\}$, 
$\Delta_2=\{-\alpha_0\}$, $\Delta_3=\{\alpha_3\}$ and $\Delta_4=\{\alpha_4\}$. Here, 
 $\Ha_1=\G_1\G_2$, $\Ha_2=\G_3\G_4$ and $G_i=\Ha_i^{F_w}\simeq\SL_2(q^2)$ for $i=1,2$.
 So condition \eqref{cond:Jv=1} holds and $v$ has only a non-trivial component 
 in $G_1$ or $G_2$ and  $q=3$. A diagram automorphism interchanges the roles of $G_1$ and 
 $G_2$, hence it is enough to look at $v=v_1\in G_1$ for $\ell\ge4$ and $q=3$.

We set $\alpha_1=\beta_1$, $\alpha_2=\beta_2$ and $\beta_3=-\alpha_0$. 
This is a base for a root system of type $A_3$. Let $\widetilde {\U}=\langle \U_{\beta_i},\,i\in\I_3\rangle$ 
and $\widetilde {\B}=\T\widetilde {\U}$. Then $\widetilde {\U}$ and $\widetilde {\B}$ are $F_w$-stable as well 
as $\Fr_q$-stable. Up to conjugation in $[\Ha,\Ha]^{F_w}$ we may assume that 
$r:=tv=\beta_1^\vee(-1)x_{\beta_1}(\xi)x_{\beta_3}(\xi')\in \beta_1^\vee(-1)\widetilde {\U}$ 
for $\xi,\xi'\in\F_{q^2}^\times$. In such a case the lemma follows,
since Lemma \ref{lem:A3} applies. 
\epf

\begin{lema}\label{lem:E7}Assume that we are in Setting \ref{set:chev} with $\Phi$ of type $E_7$ and  $q\equiv 3(4)$. Then $\Oc_{x}^{\Gb}$ is not kthulhu.
\end{lema}
\pf Recall that $\x_s^{q-1}=\x_s^2 = z \in Z(\G^F) -1$. By looking at Table \ref{tab:center} we see that $t$ is necessarily of the form:
\begin{align*}
 t=\alpha_1^\vee(\epsilon_1)\alpha_2^\vee(\eta\epsilon_2)\alpha_3^\vee(\epsilon_3)
 \alpha_4^\vee(\epsilon_4)\alpha_5^\vee(\eta\epsilon_5)\alpha_6^\vee(\epsilon_6)
 \alpha_7^\vee(\eta\epsilon_7)
\end{align*}
for $\eta$ a primitive fourth root of unity and $\epsilon_i\in\{\pm1\}$.
Since $E_7$ is simply-laced, the components of $v$ 
can only be non-trivial unipotent classes in $\PSL_2(q^{a_i})$. 
Also, by the discussion at the end of Subsection \ref{subsec:semisimple}, we may take  $w=w_0$, 
so the action of $\Ad(\w)$ does not permute the simple factors $\G_i$ of $\Ha$.  
Hence, for every $j\in \J_v$ we have $\Phi_t^j=\{\pm \alpha_{i_j}\}$ and 
$\Ha_j=\langle \U_{\pm\alpha_{i_j}}\rangle$. Assume $\alpha_k\in \Phi_t^1$. 
For every possible choice of $k\in \I_{0,7}$ we will provide a $\sigma\in W$ satisfying the 
hypotheses of Lemma \ref{lem:no-centro}. Note that Condition (1) of Lemma \ref{lem:no-centro} is always satisfied because 
$w_0=-1$ and $\theta=\id$. 

The construction of $\sigma$ relies on the fact that if $\alpha_k\in\Phi_t^1$, 
then all roots in $\widetilde {\Delta}$ that are adjacent to $\alpha_k$ do not lie in 
$\Phi_t$. These conditions pose a series of constraints on some of the $\epsilon_j$'s 
ensuring that there exists $\beta_k\in \widetilde {\Delta}$ such that $\beta_k\perp\alpha_k$ and  
$\beta_k\not\in\Phi_t$. The first property guarantees condition (2) in Lemma \ref{lem:no-centro} 
for $\sigma=s_{\beta_k}$. The second property implies $s_{\beta_k}(t)\neq t$. 
Condition (3) from Lemma \ref{lem:no-centro}  follows from Remark \ref{obs:zeta} (2)  if $\beta_k\in \Delta$, whereas 
if $\beta_k=-\alpha_0$, it follows because the equation 
$s_{\beta_k}(t)=t \alpha_2^\vee(\epsilon_1)\alpha_5^\vee(\epsilon_1)\alpha_7^\vee(\epsilon_1)$ cannot be satisfied. 

If $k=0$, then we have $\epsilon_1=1$, $\epsilon_3=-1$. Also, $\alpha_3\in\Phi_t$ $\Leftrightarrow$
$\epsilon_4=1$ $\Leftrightarrow$ $\alpha_2\not\in\Phi_t$. We take either $\beta_0=\alpha_2$ or $\alpha_3$.

If $k=1$, then we have $\epsilon_1=-1$, $\epsilon_3=\epsilon_4=1$ so 
$\alpha_2\not\in\Phi_t$ and we take $\beta_1=\alpha_2$.

If $k=2$, then we have $\epsilon_4=-1$ and $\epsilon_2\epsilon_3\epsilon_5=1$. 
Also, $\alpha_5\in\Phi_t$ $\Leftrightarrow$ $\epsilon_6=1$ $\Leftrightarrow$ $\alpha_7\not\in\Phi_t$. 
We take either $\beta_2=\alpha_5$ or $\alpha_7$.

If $k=3$, then $\epsilon_1\epsilon_4=1$, $\epsilon_3=-1$ and $\epsilon_2\epsilon_5=-1$.
Then $\alpha_0\in\Phi_t$ $\Leftrightarrow$ $\epsilon_4=1$ $\Leftrightarrow$ $\alpha_2\not\in\Phi_t$. 
We take either $\beta_3=\alpha_2$ or $\alpha_0$.

If $k=4$, then $\epsilon_2\epsilon_3\epsilon_5=-1$, $\epsilon_1=-1$, 
$\epsilon_4=\epsilon_6=1$, so $\beta_4=\alpha_7\not\in\Phi_t$.

If $k=5$, then $\epsilon_2\epsilon_3\epsilon_5=1$, $\epsilon_4\epsilon_6=-1$, 
$\epsilon_5\epsilon_7=1$. 
Then  $\alpha_7\in\Phi_t$ $\Leftrightarrow$ $\epsilon_4=1$ $\Leftrightarrow$ $\alpha_2\not\in\Phi_t$. 
We take either $\beta_5=\alpha_2$ or $\alpha_7$.

If $k=6$, then $\epsilon_4=\epsilon_6=1$ and $\epsilon_5\epsilon_7=-1$, 
so $\beta_6=\alpha_2\not\in\Phi_t$.

If $k=7$, then $\epsilon_6=-1$ and $\epsilon_5\epsilon_7=1$. 
Also, $\alpha_5\in\Phi_t$ $\Leftrightarrow$ $\epsilon_4=1$ $\Leftrightarrow$ $\alpha_2\not\in\Phi_t$. 
We take either $\beta_7=\alpha_2$ or $\alpha_5$.
\epf

\begin{lema}\label{lem:d2k+1}Assume that we are in Setting \ref{set:chev} with $\Phi$ of type $D_\ell$, $\ell=2n+1\geq 5$ and $q$ odd. 
Then $\Oc_{x}^{\Gb}$ is not kthulhu.
\end{lema}
\pf Recall that  $t^{q-1}=z\in Z(\G^F)-1$. By looking at 
Table \ref{tab:center} we see that $z$  can be either 
$z_1=(\prod_{i\leq \ell-2\text{ odd }}\alpha^\vee_i(-1))\alpha_{\ell-1}^\vee(\zeta)\alpha_\ell^\vee(\zeta^3)$, 
or $z_2=z_1^2$ 
or $z_3=z_1^3$ and the latter case is treated as the case $z=z_1$.

Let us assume $z=z_1$. This can occur only if $q\equiv1(4)$. 
We proceed as in the proof of Lemma \ref{lem:d2k}, notation as therein. We have $\pibar(z_1)=-\id_{2\ell}$, so 
$\pibar(t)^{q-1}=-\id_{2\ell}$, and the eigenvalues of $\pibar(t)$ do not lie in $\F_q$. 
Therefore, no eigenvalues of $\pibar(t)$ are equal to $\pm1$. 
By direct calculation this implies that there are no components of type $D$ in $[\Ha,\Ha]$. 
Conditions on $v$ and $q$ imply that \eqref{cond:Jv=1} holds and $v=v_{1}$ must 
live in a subgroup of type $A_1$. Such factor occurs 
when exactly two eigenvalues are repeated. Since 
$F(t)=zt$ is conjugate to $t$, it follows that $\pibar(t)$ is conjugate to $-\pibar(t)$. Hence 
if $\xi$ is an eigenvalue $\pibar(t)$, then  
$-\xi$ is again so, and $\xi^{-1}\neq -\xi$, otherwise $\xi=\pm\zeta\in\F_q$. Up to reordering, we have
\begin{align*}
\pibar(t)={\rm diag}(\xi\id_2,-\xi\id_2,\, d, -\xi^{-1}\id_2 ,\xi^{-1}\id_2)
\end{align*}
where $d$ is a diagonal matrix in $\SO_{2\ell-8}(\kk)$ and $\xi\neq\pm1$. Next we identify $w$. 
We have $\w^{-1}\trid(\pibar(t))=\pibar(t)^q=-\pibar(t)$ for some $w\in W$, so $\w$ acts on 
the first $4\times 4$-block as $C=\left(\begin{smallmatrix} 
0&\id_2\\
\id_2&0\end{smallmatrix}\right)$. This shows that $w$ interchanges the simple factors 
$\G_1$, $\G_{2}$ of type $A_1$ in $[\Ha,\Ha]$ and therefore the corresponding group $G_1$ 
is isomorphic to $\SL_{2}(q^2)$. Since $q\equiv 1(4)$, condition \eqref{cond:Jv=1} cannot be verified, concluding the case $z=z_1$. 
 
Let us assume $z=z_2=\alpha_{\ell-1}^\vee(-1)\alpha_\ell^\vee(-1)\in{\rm Ker}(\pibar)$. 
In this case, $\pibar(t)^q=\pibar(t)$ so $\pibar(\w)\in C_{\SO_{2\ell}}(\pibar(t))-\pibar(\Ha)$ and 
all  eigenvalues of $\pibar(t)$ lie  in $\F_q$. We observe that $2k$ eigenvalues equal to $\pm1$ 
give a component in $C_{\SO_{2\ell}(\kk)}(\pibar(t))$ isomorphic to ${\mathbf O}_{2k}(\kk)$, 
whereas $k$ repeated eigenvalues different from $\pm1$ give a component of type $\GL_k(\kk)$, embedded in $\SO_{2\ell}(\kk)$ 
as the group of  block diagonal matrices of the form ${\rm diag}(\id_c,A,\id_{2c'},\Jf_k\,^t\!A^{-1}\Jf_k,\id_c)$, for $c,c'\geq0$. 
Therefore, $C_{\SO_{2\ell}(\kk)}(\pibar(t))$ is the subgroup of matrices of determinant $1$ in a group 
isomorphic to ${\mathbf O}_{2a}(\kk)\times{\mathbf O}_{2b}(\kk)\times\prod_{j=1}^m\GL_{a_j}(\kk)$ with $a,\,b,\,a_j\geq0$. 
Since $C_{\SO_{2\ell}(\kk)}(\pibar(t))$ is not connected, we necessarily have $ab\neq0$, i.e., we have $2a>0$ 
eigenvalues equal to $-1$ and  $2b>0$ 
eigenvalues equal to $1$.

Up to multiplication by $z_1$ this gives 
\begin{align*}
t=\left(\prod_{j=2}^{\ell-2}\alpha_j^\vee(c_j)\right)\alpha^\vee_{\ell-1}(c)\alpha^{\vee}_\ell(-c)
\textrm{ with }c_{j}, c \in \kk^{\times},
\end{align*}
and $w=s_{\varepsilon_1-\varepsilon_\ell}s_{\varepsilon_1+\varepsilon_\ell}$.  
Thus, $w$ acts trivially on the root system generated by $\{\alpha_j,\,j\in\I_{2,\ell-2}\}$. 
If all components of $v$ are in the groups of type $A_1$ in $[\Ha,\Ha]$, then each of them can be chosen to 
be of the form $x_{\alpha_j}(\xi)$ for some $j\in \I_{3,\ell-3}$. In this case, either Lemma \ref{lem:tecnico} 
or \ref{lem:tecnico2} applies with $\beta=\alpha_j$ and $\alpha=\alpha_{j\pm1}$.  Assume a component of $v$ occurs in a factor 
$\G_1$ of type $D_a$. Then $a=2$ and it corresponds either to $\{-\alpha_0,\alpha_1\}$ or to 
$\{\alpha_{\ell-1},\alpha_{\ell}\}$. Also, $G_1\simeq \SL_2(q^2)$, so this is possible only if $q=3$ and $v=v_1$. 
In this case  we apply Lemma \ref{lem:A3} either to $\{-\alpha_0,\alpha_2, \alpha_1\}$ or to  
$\{\alpha_{\ell-1},\alpha_{\ell-2},\,\alpha_{\ell}\}$. 
\epf

\begin{lema}\label{lem:E6}Assume that we are in Setting \ref{set:chev} with $\Phi$ of type $E_6$ and $q\equiv 1(3)$. Then $\Oc_{x}^{\Gb}$ is not kthulhu.
\end{lema}
\pf Since $q\equiv 1(3)$, the element $v$ satisfies necessarily \eqref{cond:Jv=1}.
As $\Phi$ is simply laced, $\Oc_{v_1}^{G_1}$ is not a unipotent conjugacy class in $\PSp_{2k}(q^a)$.  
Thus, $\Oc_{v_1}^{G_1}$ is either isomorphic to the rack labeled by $(2)$ in $\PSL_2(q^a)$ with $q^a$ not a square 
if $q$ is odd, or the rack labeled by $(2,1^{m-2})$ in $\PSU_m(q^a)$ and the latter occurs only if $q$ is even.  
Therefore, $[\Ha,\Ha]$ must have a component  $\G_i$ of type $A_k$ and the action of $F_w^a$  on $\G_1$ should 
be twisted if $k>1$. 
Recall from Section \ref{sec:unipotent_in_K} that to each $\x_s\in\G^F$ we can associate a base $\Pi\subset\widetilde \Delta$ 
of the root system $\Phi_t$ of the connected centraliser $\Ha$ and a Weyl group element $w$ such that $w^{-1}(t)=t^q$.
The element $w$ is determined up to multiplication by elements in $W_\Pi$. Also, $w\in N_W(W_\Pi)$ because it stabilises $\Phi_t$. The pair  
$(\Pi, [w])$, where $\Pi$ is a proper subset  of $\widetilde {\Delta}$ (up to $W$-action) and $[w]= wW_\Pi$ ranges through the 
set of representatives of the conjugacy classes in $N_W(W_\Pi)/W_\Pi$ is uniquely determined up to $W$-action.

In our situation, $\w\not\in\Ha$ and  $\w^3\in\Ha$, since $\w^{-3}(t)=z^{3}t$ and $z^{3}=1$ by Table \ref{tab:center}. 
Thus, $w\in N_W(W_\Pi)- W_\Pi$ and its class in $N_W(W_\Pi)/W_\Pi$ has order $3$. 
In addition, by Remark \ref{obs:zeta} (3)  and Table \ref{tab:center} any reduced expression of $w$ must contain the reflections $s_1,\,s_3,\,s_5$ and $s_6$.
The order of a class $[w]\in N_W(W_\Pi)/W_\Pi$ can be calculated by using the package CHEVIE of GAP3 \cite{Chevie,GAP}.
For those of order $3$ we will make use of the list in \cite{FJ} of all possible pairs $(\Pi, [w])$ up to conjugation in $W$ and of the isomorphism classes of the 
corresponding $G_i^{F_w}$.  The reader should be aware that the numbering of simple roots therein differs from ours. We look through the list considering 
only the pairs $(\Pi, w)$ for which $\Pi$ is non-trivial,  $w W_\Pi=3$,  the condition from Remark \ref{obs:zeta} (3) is satisfied, and $G_1$ corresponds either to $\SL_2(q^a)$ or $\SU_m(q^a)$. 
The class $w W_\Pi$ is given by means of a representative $w\in W$. 
The following roots are used to describe $w$ as a product of 
reflections: numbering of these non-simple roots is conformal to \cite{FJ} to simplify double-checking.
\begin{align*}
&\beta_{10}=\alpha_2+\alpha_4;
&&\beta_{11}=\alpha_4+\alpha_5;\\
& \beta_{19}=\alpha_5+\alpha_6;
& & \beta_{20}=\alpha_1+\alpha_2+2\alpha_3+3\alpha_4+2\alpha_5+\alpha_6; \\
 &\beta_{21}=\alpha_1+\alpha_3;
 & &\beta_{22}=\alpha_3+\alpha_4.
 \end{align*}

We remain with the possibilities for $(\Pi,w)$ listed in Table \ref{tab:remaining-e6}. Here, we know $w$ and how it acts on each factor 
of $[\Ha,\,\Ha]$ and we proceed case-by-case.
 
 \begin{table}[ht]
		\caption{Pairs $(\Pi,wW_{\Pi})$ discussed separately in type $E_6$}\label{tab:remaining-e6}
		\begin{center}
			\begin{tabular}{|c|c|c|}
				\hline $\Pi$  & $wW_\Pi$ &  $[\Ha,\Ha]^{F_w}$ (up to isogeny) \\
			\hline
$A_1=\{-\alpha_0\}$ & $s_1s_3s_5s_6$& $\SL_2(q)$\\
				\hline
	$3A_1=\{\alpha_4,\alpha_6,-\alpha_0\}$& $s_{\beta_{11}} s_{\beta_{19}}s_{\beta_{20}} s_{\beta_{10}}$ & 	$\SL_2(q^3)$\\		
				\hline
	$4A_1=\{\alpha_1,\,\alpha_4,\,\alpha_6,\,-\alpha_0\}$ & $s_{\beta_{19}}s_{\beta_{11}}s_{\beta_{21}}s_{\beta_{22}}$  & $\SL_2(q)\times \SL_2(q^3)$\\
				\hline
			\end{tabular}
		\end{center}
	\end{table}
Let $\Pi=\{-\alpha_0\}$ with $w=s_1s_3s_5s_6$ or  $\Pi=\{\alpha_4,\alpha_6,-\alpha_0\}$ with $w=s_{\beta_{11}} s_{\beta_{19}}s_{\beta_{20}} s_{\beta_{10}}$. 
Then Lemma \ref{lem:no-centro} applies for $\sigma=s_1s_3$ or  $\sigma=s_1$.

Let $\Pi=\{\alpha_1,\,\alpha_4,\,\alpha_6,\,-\alpha_0\}$, with $w=s_{\beta_{19}}s_{\beta_{11}}s_{\beta_{21}}s_{\beta_{22}}$. Here 
$w(\alpha_0)=\alpha_0$, $w(\alpha_1)=\alpha_6$, $w(\alpha_6)=\alpha_4$ and  $w(\alpha_4)=\alpha_1$, so $v=v_1$ lies either in 
$\langle\U_{\pm\alpha_0}\rangle$ or in $\langle\U_{\pm\alpha_j},\,j=1,4,6\rangle$.
In the former case, Lemma \ref{lem:no-centro} applies with $\sigma=s_{\beta_{11}}s_{\beta_{21}}$. In the latter, the condition $-\alpha_0\in\Pi$ 
implies that $tv\in {\mathbb K}:= \langle\U_{\pm\alpha_j},\,j=1,3,4,5,6\rangle\simeq \SL_6(\kk)$ by \cite[Corollary 5.4]{sp-st}. Also, $w$ can be 
represented by an element in $N_\G(\T)\cap{\mathbb K}$, so ${\mathbb K}^{F_w}\simeq\SL_6(q)$, \cite[10.9]{St}. The claim follows from the main result in \cite{ACG-I}.
\epf

\noindent\textit{Proof of Proposition \ref{prop:tq=tz}.}
By Lemmata \ref{lem:cn} up to \ref{lem:E6}. 
\qed

\medbreak

Putting together Lemmata \ref{lem:reduction}, \ref{lem:specifico} and Propositions  \ref{prop:center}, \ref{prop:e8f4g2-qeven}, and 
\ref{prop:tq=tz} we have the main result of this Section.
\begin{theorem}
 Let $\Gb$ be a Chevalley group and $x=x_{s}x_{u} \in \Gb$ with $x_s,\,x_u\neq1$.
Then $\Oc_x^{\Gb}$ is not kthulhu, and consequently $\Oc_x^{\Gb}$ collapses.\qed
\end{theorem}

\section{Mixed classes in Steinberg groups}\label{sec:steinberg}

In this Section, $\theta\neq\id$. By Lemma \ref{lem:specifico} and  Proposition \ref{prop:center} it remains to consider the following cases:

\begin{setting}\label{set:steinberg}
$\x_s^{q+1}=z\in Z(\G^F)$; either \eqref{cond:Jv=1} or \eqref{cond:Jv=2} holds for $v$; $\Phi$, $q$ and $\x_s$ satisfy: 
\begin{enumerate}[leftmargin=*]
\item $A_\ell$, $\ell\geq 2$;
\item $D_{2n+1}$, $n\geq 2$; 
\item $E_6$, $q$ arbitrary;
\item $D_{2n}$, $n\geq 2$, $\theta^2=\id$, $\x_s^2\in Z(\G^F)$ and $q$ is odd; 
\item $D_4$, $\theta^3=\id$, $\x_s^2=1$ and $q$ is odd.
\end{enumerate}
\end{setting}
For the reader's convenience  the center $Z(\G^F)$  is recalled in Table \ref{tab:center-steinberg}.
We begin with the groups for which $w_0=-1$, i.e., cases (4) and (5).


\begin{table}[ht]
\caption{Center of $\G_{sc}^F$, Steinberg groups}\label{tab:center-steinberg}
\begin{tabular}{|p{2cm}|c|p{4cm}|c|}
\hline  type  &  $Z(\G_{sc}^F)$   & conditions  \\
\hline
$^{2}\!D_\ell$ & $\langle  \alpha^{\vee}_{\ell-1}(-1)\alpha^{\vee}_\ell(-1)\rangle$ & $q$ odd\\
\hline
$^{2}\!A_\ell$&  $\langle  \prod_{i=1}^\ell\alpha^{\vee}_{i}(\zeta^i)\rangle$  &$d=(q+1,\ell+1)$, 
$\zeta$ a primitive $d^{\rm th}$ root of $1$\\
\hline
  $^{2}E_6$ & $\langle\alpha_1^\vee(\omega)\alpha_3^\vee(\omega^2)
  \alpha_5^\vee(\omega)\alpha_6^\vee(\omega^2)\rangle$&  $q\equiv 2(3)$, $\omega^3=1$\\
\hline
$^{3}D_4$ & $1$&\\
\hline

\end{tabular}
\end{table}

\begin{lema}\label{lem:d2k-steinberg}
Assume that we are in Setting \ref{set:steinberg}, with $\Phi$ of type $D_\ell=2n$, $\theta^2=\id$ and $q$ odd. Then $\Oc$ is not kthulhu.
\end{lema}
\pf We proceed as in the proof of Lemma \ref{lem:d2k}, from which we adopt notation. 
Since $q$ is odd, $v$ has only components in type $A_1$ and  $|\J_v|>1$ only if $q=3$. By Lemma \ref{lem:specifico}
we need to consider the two possibilities 
$t^2=z_2$ or $t^2=1$. 
We consider $\pibar(t)$ as in \eqref{eq:t-inv}.  It is always an involution, its connected centralizer 
$\pibar(\Ha)$ has root system of type $D_b\times D_{\ell-b}$ and $[\Ha,\Ha]$ decomposes as in 
Lemma \ref{lem:d2k}. Here however, $\vartheta$ acts on $\T$ as $s_{\varepsilon_\ell}$ 
and we need to describe $w$. 

We observe that if $b,\ell-b>2$, the unipotent part has no components in 
Table \ref{tab:unip-kthulhu}, whereas if $b=1$ or $\ell=b+1$ the component of type $D_1$ 
is a torus and has no unipotent component.  Hence $b=2$ and/or $\ell=b+2$, so $b$ is always even. 
Up to $W$-action and multiplication by $z_2$ we have
\begin{align*}
t=\prod_{j=1}^b\alpha_j^\vee((-1)^j), 
\end{align*}
so $\vartheta(t)=t$ and $t^2=1$ necessarily. Therefore  $t\in \T^F$ and $w=1$, so $\G^{F_w}=\G^F$. 
The following cases may occur:

\noindent $\bullet$ $b=2$ and $\ell>2+b$. Then we have $\Ha_1=\G_1=\langle \U_{\pm\alpha_1}\rangle$, 
$\Ha_2=\G_2=\langle \U_{\pm\alpha_0}\rangle$, $\Ha_3=\G_3=\langle \U_{\pm\alpha_j},\,j\in\I_{3,\ell}\rangle$. 
Then up to $W$-action, either $\J_v=\{1\}$, or else $\J_v=\{1,2\}$ and $q=3$.

\noindent $\bullet$ $b>2$ and  $\ell=b+2$. Then we have $\Ha_1=\G_1=\langle \U_{\pm\alpha_j,\,j\in\I_b}\rangle$, 
$\G_2=\langle \U_{\pm\alpha_{\ell-1}}\rangle$, $\G_3=\langle \U_{\pm\alpha_{\ell}}\rangle=\vartheta(\G_2)$ so 
$\Ha_2=\G_2\G_3$, and $G_2\simeq \SL_2(q^2)$. Then $\J_v=\{2\}$ and $q=3$. 

\noindent $\bullet$ $b=2$ and $\ell=4$.  Then we have $\Ha_1=\G_1=\langle \U_{\pm\alpha_1}\rangle$, 
$\Ha_2=\G_2=\langle \U_{\pm\alpha_0}\rangle$,  $\G_3=\langle \U_{\pm\alpha_{\ell-1}}\rangle$, 
$\G_4=\langle \U_{\pm\alpha_{\ell}}\rangle=\vartheta(\G_3)$ and $\Ha_3=\G_3\G_4$ with  
$G_3\simeq \SL_2(q^2)$.
Then either $\J_v=\{1\}$, or $\J_v=\{3\}$, or $\J_v=\{1,2\}$ with $q=3$.

Summarizing, we have either $b=2$ and $v=x_{\alpha_1}(\xi)x_{-\alpha_0}(\zeta)$ for $\xi\in\kk^\times$ 
and $\zeta\in\kk$ (with $q=3$ if $\zeta\neq0$)  or 
$\ell=b+2$ and $v=x_{\alpha_{\ell-1}}(\xi)x_{\alpha_\ell}(\eta)$ (with $\xi,\,\eta\in\kk^\times$ and $q=3$). 
We set: $\beta_1=\alpha_1,\beta_2=\alpha_2, \beta_3=-\alpha_0$ in the first case and 
$\beta_1=\alpha_{\ell-1}$, $\beta_2=\alpha_{\ell-2}$, $\beta_3=\alpha_\ell$ in the second case. 
The statement then follows from Lemma \ref{lem:A3}.
\epf

\begin{lema}\label{lem:triality}Assume that we are in Setting \ref{set:steinberg} with $\Phi$ of type $D_4$, $\theta^3=\id$ and $q$ odd. Then $\Oc$ is not kthulhu.
\end{lema}

\pf 
We proceed as in the proof of Lemma \ref{lem:d2k-steinberg}, from which we adopt notation.  
We take $\vartheta$ to be a graph automorphism with associated $\theta$ given by
$\alpha_1\mapsto\alpha_3\mapsto\alpha_4$. 
Here again $q$ is odd, $t^2=1$ 
and $v$ has components only in type $A_1$.  Also, $F(t)$ is conjugate to $t$. 
These conditions imply that up to conjugation by an element in $N_\G(\T)$ 
we have either  
$t=t_1=\alpha_1^\vee(-1)\alpha_3^\vee(-1)\alpha_4^\vee(-1)\in \G^F$, or 
$t=t_2=\alpha_1^\vee(-1)$. 
We have 
$F(t_2)=\vartheta(t_2)=\alpha_3^\vee(-1)=s_{\alpha_2}s_{\alpha_1+\alpha_2+\alpha_3}(t_2)$. 
For both choices of $t$ we have $\G_1=\langle \U_{\pm\alpha_1}\rangle$, $\G_2=\langle \U_{\pm\alpha_3}\rangle$, 
$\G_3=\langle \U_{\pm\alpha_4}\rangle$ and $\G_4=\langle \U_{\pm\alpha_0}\rangle$. 
If $t=t_1$, then $w=1$, so $\Ha_1=\G_1\G_2\G_3$ and $G_1$ is either $\SL_2(q^3)$ or 
$\PSL_2(q^3)$, and $\Ha_2=\G_4$. 
If $t=t_2$, then  $w^{-1}= s_{\alpha_2}s_{\alpha_1+\alpha_2+\alpha_3}$ so 
$\Ha_1=\G_1$, $\Ha_2=\G_2\G_3\G_4$ and $G_2$ is either $\SL_2(q^3)$ or $\PSL_2(q^3)$. 
Thus, $|\J_v|=1$  and $v$ can be chosen of the form  
$x_{\beta_1}(\xi_1)x_{\beta_2}(\xi_2)x_{\beta_3}(\xi_3)$ 
with $\beta_1,\beta_2,\beta_3\in\{\alpha_0,\alpha_1,\alpha_3,\alpha_4\}$, $\xi_1\in\kk^\times$, $\xi_2,\,\xi_3\in\kk$.
We deal with the case $w=1$, the other is similar. 
We consider $s_{\alpha_2}\in W^F$ and its representative $\dot{s}_{\alpha_2}$ in $N_{\G^F}(\T)$. 
Let $r=tv$ and $s:= \dot{s}_{\alpha_2}\trid r=t'v'$. 
Since $s_{\alpha_{2}}(\beta_{1}) \notin \Phi_{t}$, 
it follows that $C_{\G}(t')\neq \Ha$ so $t'\not\in Z(\G)t$. 
Therefore,  $\pi\left(\Oc_r^{\langle r,s\rangle}\right)\cap\pi\left(\Oc_s^{\langle r,s\rangle}\right)=\emptyset$. 
As $\alpha_{2}(\pi(t)) = \alpha_{2}(t) \neq 1$ and 
$\beta_{1}(\pi(t')) = \beta_{1}(t') = \beta_{1}(s_{\alpha_{2}}(t)) =  (s_{\alpha_{2}}\beta)(t) \neq 1$, 
by \eqref{eq:comm-rule} the inequalities in 
Lemma \ref{lem:tecnico0} hold for ${\mathtt y} = \pi(r)$ and ${\mathtt z} = \pi (s)$.
Hence, $\pi(\Oc_r^{\G^F})$ and $\Oc_x^{\Gb}$ are of type C.
\epf

In the remaining groups we always have $\theta=-w_0$, so
\begin{align}\label{eq:relationE6}t=w(F(t))=w(\vartheta t^q)=ww_0(t^{-q})=z^{-1}ww_0(t).\end{align} 
Therefore, $ww_0\in N_W(W_\Pi)$ and $ww_0\in W_\Pi$ if and only if $z=1$. 
If this is the case,  since $\sigma=w_0w^{-1}\in W_\Pi$, we may replace $w$ by  $\sigma w=w_0$, so $C_W(w\theta)=C_W(-\id)=W$. 

\begin{obs}\label{obs:simple_root}Assume that $\theta=-w_0$ and $t^{q+1}=1$, $w=w_0$.
\begin{enumerate}
\item Suppose in addition that  for some irreducible component $\Delta_j$ of $\Pi$ there exists a 
simple root $\alpha_j\not\in\Pi$ such that $\alpha_j$ is orthogonal to $\Delta_j$. If 
$j\in\J_v$, then $\sigma=s_{\alpha_j}$ satisfies the hypothesis of Lemma \ref{lem:no-centro} by virtue of 
Remark \ref{obs:zeta} (2) and Table \ref{tab:center}. 
\item  If for any  component $\Delta_j$ of $\Pi$ there exists a simple root $\alpha_j\not\in\Pi$ such that 
$\alpha_j$ is orthogonal to $\Delta_j$, then the corresponding class is not kthulhu. Indeed, Lemma \ref{lem:no-centro} applies by virtue of (1).
\end{enumerate}
\end{obs}

We consider now the Steinberg groups for $\Phi$ of type $A_\ell$. Here $\G=\SL_{\ell+1}(\kk)$ and $F$ is realised as the map 
$(a_{ij})\mapsto \Jf ^t(a_{ij})^{-q}\Jf$  on any matrix of $\G$, so $\G^F=\SU_{\ell+1}(q)$ and $\Gb=\PSU_{\ell+1}(q)$.

\begin{obs}\label{rem:psu32-no-mixed}
We shall not consider  $\PSU_3(2)$  because it has no mixed classes,
since every semisimple element in $\SU_3(2)$ is either central or regular.
\end{obs}

\begin{lema}\label{lem:su3-matrices}Let $\Gb=\PSU_3(q)$ with $q>2$. 
Let $\lambda\in\F_{q^2}$ with $\lambda^{q+1}=1$, $\lambda^3\neq1$ and let 
 \begin{align*}
&\x_1=\left(\begin{smallmatrix}
\lambda&0&\lambda\\
0&\lambda^{-2}&0\\
0&0&\lambda
     \end{smallmatrix}\right).
     \end{align*}
Then $\Oc_{\pi(\x_1)}^{\Gb}$ is not kthulhu. 
\end{lema}

\pf If $q$ is odd, then a slight modification of  the proofs of \cite[Lemma 5.3]{ACG-IV} shows that 
$\Oc_{\pi(\x_1)}^{\Gb}$ is of type D for $q>3$ and type C for $q=3$. 
Let thus $q$ be even $>2$. We show that in this case,
$\Oc_{\pi(\x_1)}^{\Gb}$ is of type $D$.

Under these assumptions $\U^F$ is the group of matrices of the form
$\left(\begin{smallmatrix}
1&\xi&\eta\\
0&1&\xi^q\\
0&0&1
     \end{smallmatrix}\right)$ for $\xi,\eta\in\F_{q^2}$ such that $\eta^q+\eta=\xi^{q+1}$. The inverse of such a matrix is $\left(\begin{smallmatrix}
1&\xi&\eta^q\\
0&1&\xi^q\\
0&0&1
     \end{smallmatrix}\right)$. These elements have order $4$ if $\xi\neq0$ and order $2$ if $\xi=0$, $\eta\neq0$.

Since $\lambda^q=\lambda^{-1}$ and $\lambda^3\neq1$, it follows that 
$\lambda^6\neq1$ so $\lambda^{3q}=\lambda^{-3}\neq\lambda^3$ and thus 
$\lambda^3\not\in\F_q$.  Let $a\in \F_{q^2}^{\times}$ such that $a^{q+1}\neq1$ and $\eta\in\F_{q^2}$ such that 
$\eta+\eta^q=1$.  We set $r:=\x_1$ and
\begin{align*}y&:=\left(\begin{smallmatrix}
a&a^{q-1}&\eta a^{-q}\\
0&a^{q-1}&a^{-q}\\
0&0&a^{-q}
     \end{smallmatrix}\right)\in\SU_3(q),\\
 s&:=y\trid r=\left(\begin{smallmatrix}
\lambda&(\lambda+\lambda^{-2})&(\lambda+\lambda^{-2})+\lambda a^{1+q}\\
0&\lambda^{-2}&(\lambda+\lambda^{-2})\\
0&0&\lambda
     \end{smallmatrix}\right)\in\Oc_r^{\SU_3(q)},\\
    u &=\left(\begin{smallmatrix}
1&(1+\lambda^{-3})&(1+\lambda^{-3})+ a^{1+q}\\
0&1&(\lambda^3+1)\\
0&0&1
     \end{smallmatrix}\right)\in\U^F,\\
    t&:=\left(\begin{smallmatrix}
\lambda&0&0\\
0&\lambda^{-2}&0\\
0&0&\lambda
     \end{smallmatrix}\right)\in\T^F .
\end{align*}
so $s=tu$ and $r=t(\id+e_{13})$. 

The $(1,2)$ entry of $(rs)^2$ equals $(\lambda+\lambda^{-2})\lambda^3(1+\lambda^{-6})$ 
whereas the $(1,2)$ entry of $(sr)^2$ equals $(\lambda+\lambda^{-2})(1+\lambda^{-6})$, hence $(rs)^2\neq(sr)^2$.

These products have same diagonal part, so $\pi(rs)^2\neq\pi(sr)^2$. Let $H=\langle r,s\rangle=\langle t, (\id+e_{1,3}),\,tu\rangle=\langle t,  (\id+e_{1,3}), u\rangle$. 
Observe that $u^2$ and $(\id+e_{1,3})\in Z(H)$ so the elements in $\Oc_r^H$ and $\Oc_s^H$ have same 
diagonal part. Thus  $\Oc_r^H\neq\Oc_s^H$ implies $\pi(\Oc_r^H)\neq\pi(\Oc_s^H)$. We verify the former. We have:
\begin{align*}
H/Z(H)&=\left\langle t Z(H),u Z(H)\right\rangle\\
&=\{t^{a_1}ut^{a_2}u\cdots t^{a_{k-1}}ut^{a_k}Z(H),\,a_1\geq0 ,a_k\geq0, a_i\geq1,\,i\in\I_{2,k-1}\}.
\end{align*}
Since $t\trid r=r$ we obtain 
\begin{align*}\Oc_r^H&=\Oc_r^{H/Z(H)}\\
&=\{((t^{a_1}\trid u)(t^{a_1+a_2}\trid u)\cdots (t^{a_1+\cdots+a_k} \trid u ))\trid r,\,  a_1\geq0, \,a_i\geq1 \textrm{ for } i\in\I_{2,k}\}.
\end{align*}
 Observe that for $m_i\geq1$
 \begin{align*}
& t^{m_i}\trid u=\left(\begin{smallmatrix}
1&\lambda^{3m_i}(1+\lambda^{-3})&a^{1+q}+(1+\lambda^{-3})\\
0&1&\lambda^{-3m_i}(1+\lambda^{3})\\
0&0&1
     \end{smallmatrix}\right),\\
  &    (t^{m_1}\trid u)(t^{m_2}\trid u)\cdots (t^{m_k} \trid u )=\left(\begin{smallmatrix}
1&(\sum_i \lambda^{3m_i})(1+\lambda^{-3})&f\\
0&1&(\sum_i\lambda^{-3m_i})(1+\lambda^{3})\\
0&0&1
     \end{smallmatrix}\right)
     \end{align*}
for some $f\in\F_{q^2}$ such that $f+f^q=\left(\sum_i \lambda^{3m_i}\right)(1+\lambda^{-3})\left(\sum_i\lambda^{-3m_i}\right)(1+\lambda^{3})$ and 
 \begin{align*}
 & (t^{m_1}\trid u)(t^{m_2}\trid u)\cdots (t^{m_k} \trid u )\trid r\\
 &=t \left(\begin{smallmatrix}
1& *&f^q+f+1+\lambda^{-3}(1+\lambda^3)(1+\lambda^{-3})(\sum_i \lambda^{3m_i})(\sum_i\lambda^{-3m_i})\\
0&1&*\\
0&0&1
     \end{smallmatrix}\right).
     \end{align*}
Thus, if  $s\in\Oc_r^H$, then comparing $(1,3)$-entries we would have 
\begin{align*}
(1+\lambda^{-3})+a^{1+q}=1+(1+\lambda^{-3})\left(\sum_i \lambda^{3m_i}\right)(1+\lambda^{-3})\left(\sum_i\lambda^{-3m_i}\right)(1+\lambda^{3})
\end{align*}
for some nonnegative integers $m_i$, or, equivalently 
\begin{align*}
1+a^{1+q}=(1+\lambda^{-3})\left[\left(\sum_i \lambda^{3m_i}\right)(1+\lambda^{-3})\left(\sum_i\lambda^{-3m_i}\right)(1+\lambda^{3})+1\right].
\end{align*}
Our choice of $a$ implies that the left hand side lies in $\F_{q}^\times$, whereas the right hand side is a product of 
$(1+\lambda^{-3})\in \F_{q^2}\setminus \F_q$ with an element in $\F_q$. The latter  cannot lie in $\F_q$, unless it is zero, so this equality cannot hold. 
\epf

\begin{lema}\label{lem:su-3matrices}
Let $\Gb=\PSU_{\ell+1}(q)$ for $\ell=3,4$. 
Assume $t^{q+1}=1$ and that $tv$ is conjugate in $\G$ to one of the following matrices, for 
$\lambda_{1},\lambda_{2}\in \F_{q^{2}}^{\times}$ with $\lambda_{1}\neq \lambda_{2}$: 
 \begin{align*}
&\x_2=\left(\begin{smallmatrix}
\lambda_1&0&0&\lambda_1'\\
0&\lambda_2&\lambda_2'&0\\
0&0&\lambda_2&0\\
0&0&0&\lambda_1
     \end{smallmatrix}\right),\;(\lambda_1',\lambda_2')\neq(0,0),\,\textrm{ and }\lambda_1'\lambda_2'\neq0\textrm{ only for }q=2,3 \\
     &\x_3=\left(\begin{smallmatrix}
\lambda_1&0&0&0&\lambda_1'\\
0&\lambda_2&0&\lambda_2'&0\\
0&0&\lambda_2&0&0\\
0&0&0&\lambda_2&0\\
0&0&0&0&\lambda_1
     \end{smallmatrix}\right),\;  q=2, (\lambda_1',\lambda_2')\neq(0,0).
     \end{align*}
Then $\Oc_x^{\Gb}$ is of type D. 
\end{lema}
\pf
Let $\x=\x_2$, $\sigma={\rm diag}(\Jf_2,\Jf_2)\in\SU_{4}(q)$, $y=\left(\begin{smallmatrix}
1&1&0&0\\
0&1&0&0\\
0&0&1&-1\\
0&0&0&1
     \end{smallmatrix}\right)\in\SU_4(q)$ and 
     \begin{align*}
     &r=\sigma\trid \x=\left(\begin{smallmatrix}
\lambda_2&0&0&\lambda_2'\\
0&\lambda_1&\lambda_1'&0\\
0&0&\lambda_1&0\\
0&0&0&\lambda_2
     \end{smallmatrix}\right), & & s=y\trid\x =\left(\begin{smallmatrix}
\lambda_1&\lambda_2-\lambda_1&\lambda_2'&\lambda_1'+\lambda_2'\\
0&\lambda_2&\lambda_2'&\lambda_2'\\
0&0&\lambda_2&\lambda_2-\lambda_1\\
0&0&0&\lambda_1
     \end{smallmatrix}\right).
\end{align*}
By looking at the $(1,2)$ entry we see that $(rs)^2\not\in\kk(sr)^2$. 
In addition, 
 \begin{align*}
     &\Oc_r^{\langle r,s\rangle}\subset  \left(\begin{smallmatrix}
\lambda_2&*&*&*\\
0&\lambda_1&\lambda_1'&*\\
0&0&\lambda_1&*\\
0&0&0&\lambda_2
     \end{smallmatrix}\right), & & \Oc_s^{\langle r,s\rangle}\subset\left(\begin{smallmatrix}
\lambda_1&*&*&*\\
0&\lambda_2&\lambda_2'&*\\
0&0&\lambda_2&*\\
0&0&0&\lambda_1
     \end{smallmatrix}\right).\end{align*}
hence looking at the diagonal or at the $(2,3)$-entry we see that 
$\pi(\Oc_r^{\langle r,s\rangle})\neq\pi(\Oc_r^{\langle r,s\rangle})$, so $\Oc_x^{\Gb}$ is of type D,
with a possible exception when $\lambda_1=-\lambda_2$ and $\lambda_1'=-\lambda_2'\neq0$ and $q=3$.
Assume this is the case, so $\lambda_{1}^{2}+{\lambda_{1}'}^{2}=0$
and $(\lambda_1')^4=\lambda_1^4=1$. 
For every choice of $\lambda_1$ there are two possible choices of $\lambda_1'$ but the corresponding elements are conjugate by
${\rm diag}(\zeta,\zeta^{-1},\zeta^3,\zeta^{-3})$ where $\langle\zeta\rangle=\F_9^\times$.
In addition, the different choices of $\lambda_1$ correspond to multiplication by a central element. Therefore it remains to consider
$\x_2=\left(\begin{smallmatrix}
\zeta^2&0&0&2\\
0&\zeta^6&1&0\\
0&0&\zeta^6&0\\
0&0&0&\zeta^2
     \end{smallmatrix}\right)$. The element 
     $y:= \left(\begin{smallmatrix}
0&\zeta^7&0&\zeta^7\\
\zeta^5&0&\zeta&0\\
0&0&0&\zeta\\
0&0&\zeta^3&0
     \end{smallmatrix}\right)$ lies in $\Oc_{\x_2}^{\G^F}$ because it has the same Jordan form as $\x_2$ and 
$(\x_2 y)^2(y\x_2)^{-2}=\left(\begin{smallmatrix}
   1 &0&\zeta^5 &0\\
   0&1 &0 &\zeta^3\\
   0&0&   1  & 0\\
   0&0&0&  1
 \end{smallmatrix}\right)\not\in Z(\SU_4(q))$.   
A computation with GAP shows that $\Oc_{\x_2}^{\SU_4(q)}\neq\Oc_{zy}^{\SU_4(q)}$ for any $z\in Z(\SU_4(q))$, so $\Oc_{x}^{\Gb}$ is of type D.

\bigskip

Let now $\x=\x_3$, $\sigma'={\rm diag}(\Jf_2,1,\Jf_2)\in\SU_5(q)$, $z=\left(\begin{smallmatrix}
1&1&0&0&0\\
0&1&0&0&0\\
0&0&1&0&0\\
0&0&0&1&-1\\
0&0&0&0&1
     \end{smallmatrix}\right)\in\SU_5(q)$ and 
     \begin{align*}
     &r=\sigma'\trid \x=\left(\begin{smallmatrix}
\lambda_2&0&0&0&\lambda_2'\\
0&\lambda_1&0&\lambda_1'&0\\
0&0&\lambda_2&0&0\\
0&0&0&\lambda_1&0\\
0&0&0&0&\lambda_2
     \end{smallmatrix}\right), & & s= z\trid\x =\left(\begin{smallmatrix}
\lambda_1&\lambda_2-\lambda_1&0&\lambda_2'&\lambda_1'+\lambda_2'\\
0&\lambda_2&0&\lambda_2'&\lambda_2'\\
0&0&\lambda_2&0&0\\
0&0&0&\lambda_2&\lambda_2-\lambda_1\\
0&0&0&0&\lambda_1
     \end{smallmatrix}\right).
\end{align*}
As in the previous case we see that $(rs)^2\not\in\kk(sr)^2$ and that
 \begin{align*}
     &\Oc_r^{\langle r,s\rangle}\subset  \left(\begin{smallmatrix}
\lambda_2&*&0&*&0\\
0&\lambda_1&0&\lambda_1'&*\\
0&0&\lambda_2&0&0\\
0&0&0&\lambda_1&*\\
0&0&0&0&\lambda_2
     \end{smallmatrix}\right), & & \Oc_s^{\langle r,s\rangle}\subset\left(\begin{smallmatrix}
\lambda_1&*&0&*&*\\
0&\lambda_2&0&\lambda_2'&*\\
0&0&\lambda_2&0&0\\
0&0&0&\lambda_2&*\\
0&0&0&0&\lambda_1
     \end{smallmatrix}\right).\end{align*}
Hence looking at the diagonal we conclude that $\pi(\Oc_r^{\langle r,s\rangle})\neq\pi(\Oc_r^{\langle r,s\rangle})$ so $\Oc_x^{\Gb}$ is of type D.
\epf

\begin{lema}\label{lem:su-z=1}Assume that we are in Setting \ref{set:steinberg} with $\Phi$ of type $A_\ell$, $\ell\geq2$ 
and $\x_s^{q+1}=1$. Then $\Oc_x^{\Gb}$ is not kthulhu. 
\end{lema}

\pf Observe that \eqref{eq:relationE6} gives $w=w_0$ and $C_W(w\theta)=W$.  
Also, $F_w$ preserves each irreducible component of $\Phi_t$ mapping each root to its opposite,  hence,  
each $\Ha_i$ is simple. An eigenvalue of $t$ of multiplicity $k\geq2$ gives a component  of type $A_{k-1}$ in $\Pi$. 
If $t$ has at least $3$ different
eigenvalues, then  $\left|\widetilde \Delta- \Pi\right|\geq3$, so the statement follows from Remark \ref{obs:simple_root}. 

Assume $t$ has exactly two eigenvalues $\lambda_1,\,\lambda_2$ of multiplicity $m_1,\,m_2$. They satisfy $\lambda_i^{q+1}=1$.
According to the parity of $m_1$ and $m_2$, and up to interchanging indices, $t$ is conjugate in $\G$ to one of the following matrices:
\begin{align*}
t_1&=\textrm{diag}(\lambda_1\id_{a},\lambda_2\id_{m_2},\lambda_1\id_a),&&\textrm{ if } m_1=2a>0\\
t_2&=\textrm{diag}(\lambda_1\id_{a},\lambda_2\id_{b},\lambda_1,\lambda_2\id_{b+1},\lambda_1\id_a),&& \textrm{ if } m_1=2a+1\geq m_2=2b+1.
\end{align*}
Then $t_1\in\SU_{\ell+1}(q)$ so we take $\x_s=t_1=t$ (and we replace for convenience $w=w_0$ by $w=\id$) 
in this case, whereas $t_2\in \G^{F_w}$ for the choice $w=s_{a+b+1}$.

Let $\x_s=t_1$. Then $[C_{\G}(t_1),C_{\G}(t_1)]$ is given by matrices of block form 
$\left(\begin{smallmatrix}A_1&0&A_{2}\\
0&B&0\\
A_{3}&0&A_{4}  
 \end{smallmatrix}\right)$ with $A_j$ of size $a\times a$ for $j\in\I_4$ and $B$ of size 
 $m_2\times m_2$ and $[C_{\G}(t_1),C_{\G}(t_1)]^F\simeq \SU_{m_1}(q)\times \SU_{m_2}(q)$. We take 
 $\x_u=v\in\id+\kk e_{1,m_1+m_2}+\kk e_{a+1,a+b}\id\in \SU_{m_1+m_2}(q)$. All cases can be reduced to computations 
 in orbits represented by elements with shape $\x_j$, for $j\in 1,2,3$ as in 
 Lemmata \ref{lem:su3-matrices} and \ref{lem:su-3matrices}, under the action of subgroups isomorphic to  $\SU_j(q)$ for $j=3,4,5$.
 Such elements may have determinant different from $1$ but the proof of these Lemmata does not require this assumption.
 Hence all cases are covered, in view of Remark \ref{rem:psu32-no-mixed}.
 
 We study now the class of $t_2v\in\Oc_{t_2v}^{\G^{F_w}}$ for $w=s_{a+b+1}$. 
 The non-trivial factors in $[C_{\G}(t_2),C_{\G}(t_2)]$ are of type $A_{2a}$ and $A_{2b}$, so $q$ is necessarily even.
 We work in 
$\G^{F_w}$
 where $F_w=\Ad(\dot{s}_{a+b+1}) F$, $\dot{s}_{a+b+1}=\left(\begin{smallmatrix}
 \id_{a+b}&0&0&0\\
 0&0&1&0\\
 0&1&0&0\\
 0&0&0&\id_{a+b}\end{smallmatrix}\right)$. We may take $v\in\id+\kk e_{1,\ell+1}+\kk e_{a+1,\ell+1-a}$ where 
 two nonzero elements outside the diagonal can occur only for $a=b=1$ and $q=2$. 
It represents the class because it has the right Jordan form, centralises $t_2$ and, 
for suitable choice of the scalars, it lies in $\G^{F_w}$. We take
\begin{align*}
\x&=\textrm{diag}(\lambda_1\id_{a},\lambda_2\id_{b},\lambda_1,\lambda_2\id_{b+1},\lambda_1\id_a)(\id+\xi_1 e_{1,\ell+1}+\xi_2 e_{a+1,\ell+1-a})\\
&\textrm{with $\xi_1\neq0$ and  $\xi_2=0$ if $b=0$,}
\end{align*}
the case $ \xi_2\neq0$ is treated similarly.  Let 
\begin{align*}
\sigma_1&=\left(\begin{smallmatrix}
 \id_{a+b}&0&0&0\\
 0&0&1&0\\
 0&1&0&0\\
 0&0&0&\id_{a+b}
     \end{smallmatrix}\right)\in\G^{F_w},\qquad \sigma_2=\left(\begin{smallmatrix}
0&0 &1\\
0& \id_{\ell-2}&0\\
1& 0&0
\end{smallmatrix}\right)\in\G^{F_w}\\
r&=\sigma_1\trid \x\\
&=\textrm{diag}(\lambda_1\id_{a},\lambda_2\id_{b+1},\lambda_1,\lambda_2\id_{b},\lambda_1\id_a)(\id+\xi_1 e_{1,\ell+1}+\xi_2 e_{a+1,\ell+1-a}),\\
s&=\sigma_{2}\trid \x\\
&=\textrm{diag}(\lambda_1\id_{a},\lambda_2\id_{b},\lambda_1,\lambda_2\id_{b+1},\lambda_1\id_a)(\id+\xi_1 e_{\ell+1,1}+\xi_2 e_{a+1,\ell+1-a}).
\end{align*}
By looking at the $(1,1)$ and $(1,\ell+1)$ entries we verify
that $\pi(rs)^2\neq\pi(sr)^2$ and by looking at the entries 
$(j,j)$ for $j\in\I_{a+b,a+b+2}$ 
we verify that $\pi(\Oc_r^{\langle r,s\rangle})\neq \pi(\Oc_s^{\langle r,s\rangle})$, hence $\Oc_x^{\Gb}$ is of type D.
 \epf

\begin{lema}\label{lem:su-znot1}
Assume that we are in Setting \ref{set:steinberg} with $\Phi$ of type $A_\ell$, $\ell\geq2$ 
and $\x_s^{q+1}=z\in Z(\G^F)-1$. Then $\Oc_x^{\Gb}$ is not kthulhu.
\end{lema}

\pf Here $t^{q+1}=\eta\id$ for some $\eta\in\kk^\times$ such that $\ord \eta=b$ is a divisor of $d=(q+1,\ell+1)$, see Table \ref{tab:center-steinberg}.
Equation \eqref{eq:relationE6} shows that $t$ and $zt$ are conjugate, 
hence the spectrum of $t$ is the disjoint union of sets of the form 
$S_j=\{\lambda_j\eta^l,\, l\in\I_{0,b-1}\}$, for $j\in\I_k$, and the multiplicity is constant in $S_j$: 
we denote it by $m_j$. Hence $t$ can be chosen to have form
\begin{align*}
t={\rm diag}(\lambda_1\id_{m_1},\lambda_1\eta\id_{m_1},\ldots,\eta^{b-1}\lambda_1\id_{m_1},\lambda_2\id_{m_2},\ldots,\lambda_k\eta^{b-1}\id_{m_k}). 
\end{align*}
Then $w(F(t))=t$ for $w=\tau w_0$, where $\tau$ is the permutation represented by 
the block-diagonal monomial matrix in ${\GL}_{\ell+1}(\kk)$
\begin{align*}
&\tau={\rm diag}(\tau_1,\,\ldots,\tau_k),&&
\tau_j=\left(
\begin{smallmatrix}
0&\id_{m_j}&0&\ldots&0\\
0&0&\id_{m_j}&\ldots&0\\
\vdots&\vdots&\vdots&\vdots&\vdots\\
0&0&0&0&\id_{m_j}\\
\id_{m_j}&0&\dots&\dots&0
\end{smallmatrix}
\right).
\end{align*}
Decomposition \eqref{eq:deco2} consists of a factor $\Ha_j\simeq\prod_{j=1}^d\SL_{m_j}(\kk)$, 
for each $j\in\I_k$ such that $m_j>1$. More precisely, $\omega=\tau w_0\theta$ cyclically permutes 
the irreducible components of the root systems of $\Ha_j$ and
maps a simple root in one factor to the opposite of a simple root in the following factor. 
Hence, $G_j\simeq \SL_{m_j}(q^d)$ when $d$ is even or $m_j=2$ and 
$G_j\simeq \SU_{m_j}(q^d)$ when $d$ is odd and $m_j>2$.

Assume $k>1$ and let  $j_0\in\J_v$. By construction, $C_W(w\theta)=C_W(\tau)$, 
so it contains the permutations corresponding to each factor $\tau_j$, they preserve every $\Ha_j$, and 
\begin{align*}\tau_j(t)={\rm diag}(\id_{d(m_1+\cdots+m_{j-1})},\eta \id_{m_j},\id_{d(m_{j+1}+\cdots+m_{k})})t\not\in Z(\G^F)t.
\end{align*}
 Thus, for any $j\neq j_0$, the element $\tau_j$ satisfies the hypotheses of 
 Lemma \ref{lem:no-centro} and $\Oc_x^{\Gb}$ is not kthulhu.

Assume now $k=1$, so $|\J_v|=1$ and write for simplicity $\lambda_1=\lambda$, $m_1=m$, so $\ell+1=dm$. Then 
\begin{align*}
1=\det(t)=\lambda^{\ell+1}\eta^{\frac{d(d-1)m}{2}}.
\end{align*}
Hence, if either $m$ is even or $d$ is odd, we have $\lambda^{\ell+1}=1$ so 
$\lambda\id\in Z(\SU_{\ell+1}(q))$ and multiplying by $\lambda^{-1}\id$ we reduce 
to the case in which $t^{q+1}=1$ and invoke Lemma \ref{lem:su-z=1}. If, instead, $m$ is odd and $d$ is even, then 
$G_1\simeq \SL_{m}(q^d)$ with $m>2$ and $q^d>2$, so $\x$ cannot satisfy condition \eqref{cond:Jv=1}.  
\epf

\begin{lema}\label{lem:su}
Assume that we are in Setting \ref{set:steinberg} with $\Phi$ of type $A_\ell$, $\ell\geq2$. 
Then $\Oc_x^{\Gb}$ is not kthulhu.
\end{lema}

\pf The case $\x_s^{q+1}=1$ is covered by Lemma  \ref{lem:su-z=1} 
whereas the case $\x_s^{q+1}\neq1$ is covered by Lemma \ref{lem:su-znot1}.\epf


\begin{lema}\label{lem:d2k1-steinberg}
Assume that we are in Setting \ref{set:steinberg}, with $\Phi$ of type $D_{\ell}$ with $\ell=2n+1\geq 5$.  
Then $\Oc$ is not kthulhu.
\end{lema}

\pf  We adopt the notation
of Lemma \ref{lem:d2k}. The non-trivial element $z\in Z(\G^F)$ generates the kernel of the projection 
$\G\to\SO_{2\ell}(\kk)$, so for $t'\in\T$ we have $t\neq zt'$ if and only if $\pi(t)\neq\pi(t')$. 
In this situation $[\Ha,\,\Ha]$ has root system of type $D_a\times D_b\times \prod_jA_{k_j}$, where 
$D_1$ corresponds to a torus, $D_2$ is $A_1\times A_1$, $D_3=A_3$. 
Recall that the components of type $A_{k_j}$ correspond to equal eigenvalues $\lambda$ of 
$\pi(t)\in\SO_{2\ell}(\kk)$ distinct from $\pm1$, whereas $D_a$ and $D_b$ 
correspond to a sequence of 
equal eigenvalues $\lambda=\pm1$. Thus, if $q$ is even there is at most one component of type $D$ in $[\Ha,\Ha]$. 

Assume first that $t^{q+1}=1$ and take $w=w_0$ so $C_W(w\theta)=W$. 
Thus, $w\theta$ acts as $-\id$ on $\Phi$, so all $\G_i$ are $F_w$-stable and 
if $\G_j$ is of type $A_{k}$, then $G_j\simeq \SU_{k+1}(q)$ for $k\geq1$, where we set
$\SU_2(q):=\SL_{2}(q)$.  By hypothesis, $\J_v\neq\emptyset$ and \eqref{cond:Jv=1} or \eqref{cond:Jv=2}  impose restrictions on the subgroups $\Ha_j$ with $j\in\J_v$.
For most possibilities for $\Ha_j$ we will prove the statement by exhibiting a $\sigma\in C_W(w\theta)$ satisfying the hypotheses of Lemma \ref{lem:no-centro}. In the remaining case we will apply Lemma \ref{lem:A3}

If $\Ha_j\simeq\SL_{2l}(\kk)$ is of type $A_{2l-1}$ with base  
$\{\alpha_c,\,\ldots,\,\alpha_{c+2l-2}\}=\{\varepsilon_c-\varepsilon_{c+1},\,\ldots,\,\varepsilon_{c+2l-2}-\varepsilon_{c+2l-1}\}$, 
then  we take
$\sigma=\prod_{i=0}^{l-1}s_{\varepsilon_{c+2i}-\varepsilon_{c+2i+1}} s_{\varepsilon_{c+2i}+\varepsilon_{c+2i+1}}$. 
Condition $\sigma(\pi(t))\neq\pi(t)$ holds because  $\sigma$ interchanges the eigenvalue $\lambda$ with $\lambda^{-1}$.  
If $\Ha_j\simeq \SL_{2l+1}(\kk)$ is of type $A_{2l}$ for some $l$, and has root system with base  
$\{\alpha_c,\,\ldots,\,\alpha_{c+2l-1}\}=\{\varepsilon_c-\varepsilon_{c+1},\,\ldots,\,\varepsilon_{c+2l-1}-\varepsilon_{c+2l}\}$,
then we take 
$\sigma=\left(\prod_{i=0}^{l-1}s_{\varepsilon_{c+2i}-\varepsilon_{c+2i+1}} s_{\varepsilon_{c+2i}+\varepsilon_{c+2i+1}}\right)
s_{\varepsilon_{c+2l-1}-\varepsilon_{c+2l}} s_{\varepsilon_{c+2l-1}+\varepsilon_{c+2l}}$.

If $\Ha_j$ is of type $D_3\simeq A_3$ then the eigenvalue $\lambda=\pm1$ has multiplicity $6$ and  $q$ 
is necessarily even so $\lambda=1$ and all remaining eigenvalues of $\pi(t)$ are different from their inverses. 
We take $\sigma=\prod_{i=0}^{n-2}s_{\varepsilon_{4+2i}-\varepsilon_{4+2i+1}} s_{\varepsilon_{4+2i}+\varepsilon_{4+2i+1}}$, recalling that $\ell=2n+1$. 

Assume finally that there are no factors as above, so $\Ha_j$ is a simple factor contained in the component of type $D_2\simeq A_1\times A_1$. 
The base of the root system of type $D_2$ is either $\{\alpha_0,\,\alpha_1\}$ or $\{\alpha_{\ell-1},\,\alpha_\ell\}$.  
In the first case we take $\sigma=s_{\varepsilon_3-\varepsilon_4}s_{\varepsilon_3+\varepsilon_4}$, in the latter we take 
$\sigma=s_{\varepsilon_1-\varepsilon_2}s_{\varepsilon_1+\varepsilon_2}$. Condition 3 in Lemma \ref{lem:no-centro} holds 
unless $\pi(t)$ has only eigenvalues $\pm1$, or, equivalently, $\Phi_t$ has type $D_2\times D_{\ell-2}$. 
If this is the case, we replace $w$ by $w'=w_0w_\Pi$, where $w_\Pi$ is the longest element in $W_\Pi$ and  we replace  $F_{w}$ by $F_{w'}$. 
Then  $w'\theta$ acts trivially on the root system of type $A_3$ containing $D_2$, so Lemma \ref{lem:A3} applies with 
$\beta_1=\alpha_1,\,\beta_2=\alpha_2,\,\beta_3=-\alpha_0$. This concludes the case $t^{q+1}=1$.

Assume now $t^{q+1}=z\in Z(\G^F)-1$. Then $q$ is necessarily odd, see Table \ref{tab:center-steinberg}.  
Also, \eqref{eq:relationE6} shows that 
$ww_0\not\in W_\Pi$ whereas  $ww_0$ fixes $\pi(t)$. Hence, $C_{\SO_{2\ell}}(\pi(t))\supsetneq \pi(\Ha)$. 
Arguing as in the proof of Lemma \ref{lem:d2k+1} we deduce that this can happen only if $\pi(t)$ has both eigenvalues equal to 
$1$ and $-1$. Since $\pi\subset\widetilde {\Delta}$, we have
$\pi(t)_{jj}=\epsilon$ for $\epsilon^2=1$ and $j\in\I_{1,a}\cup \I_{2\ell-a+1,2\ell}$ and $\pi(t)_{jj}=-\epsilon$ for $j\in\I_{\ell-b+1,\ell+b}$. Also, 
\begin{align*}C_{\SO_{2\ell}}(\pi(t))/\pi(\Ha)\simeq S({\mathbf O}_{2a}(\kk)\times {\mathbf O}_{2b}(\kk))/ \SO_{2a}(\kk)\times \SO_{2b}(\kk)\simeq \Z/2\Z\end{align*}
and  by \cite[2.2]{hu-cc} it is isomorphic to $\langle W_\Pi, ww_0\rangle/W_\Pi$. 
Thus, $ww_0W_\Pi=s_{\varepsilon_1+\varepsilon_\ell}s_{\varepsilon_1-\varepsilon_\ell}W_\Pi$ and therefore we may assume 
$w=s_{\varepsilon_1+\varepsilon_\ell}s_{\varepsilon_1-\varepsilon_\ell}w_0$. Hence, 
$C_W(ww_0)=C_W(s_{\varepsilon_1+\varepsilon_\ell}s_{\varepsilon_1-\varepsilon_\ell})$. 
Since $q$ is odd, $G_j$ is of type $A_1$ for every $j\in\J_v$. Such a factor comes from an eigenvalue 
$\lambda$ of $\pi(t)$ of multiplicity $2$ if $\lambda\neq\pm1$ and of multiplicity $4$ if $\lambda=\pm1$. 
Let $j\in\J_v$. If the corresponding $\lambda\neq\pm1$ and $\pi(t)_{ll}=\pi(t)_{l+1,l+1}=\lambda$, then 
$\sigma= s_{\varepsilon_l+\varepsilon_{l+1}}s_{\varepsilon_l-\varepsilon_{l+1}}$ satisfies the hypotheses of Lemma \ref{lem:no-centro}. 
If, instead $\pi(t)_{jj}=\pm1$ for every $j\in\I_\ell$, then \eqref{cond:Jv=1} or \eqref{cond:Jv=2} holds only if either $1$ or $-1$ has multiplicity 
$4$ and the root system of $G_j$ is one of the irreducible components of the root system of type $D_2$. 
Acting possibly by $\tau=s_{\varepsilon_1-\varepsilon_\ell}s_{\varepsilon_2-\varepsilon_{\ell-2}}$ we may assume that the multiplicity $4$ 
eigenvalue occurs in the entries indexed by $1,2,2\ell,2\ell-1$. In this situation, $D_2$ has base $\{\alpha_0,\alpha_1\}$. 
As in the case $t^{q+1}=1$ we replace $w$ by $w'=ww_\Pi$ and apply Lemma \ref{lem:A3} to $\{-\alpha_0,\alpha_2,\alpha_1\}$.
\epf

Finally, we deal with the groups $^{2}\!E_6(q)$. In order to use the list in \cite{FJ}, we need to establish a dictionary between pairs 
$(\Pi, [w])$ coming from different choices for the Steinberg endomorphism $F$.
\begin{obs}\label{obs:dictionary-e6} For $j=1,2$, let $F_j=\vartheta_j \Fr_q$ be Steinberg endomorphisms of $\G$  such that 
$F_1={{\rm Ad} \sigmad}\, F_2$ for some $\sigma\in N_\G(\T)$ and let $g_0\in\G$ be such that $g_0^{-1}F_2(g_0)=\sigmad$. 
Then $\G^{F_2}=g_0\G^{F_1}g_0^{-1}$ and if  for some $z\in\G$ we have $\x\in \G^{F_1}\cap \Oc_z^{\G}$, it follows that 
$g_0\x g_0^{-1}\in \G^{F_2}\cap \Oc_z^{\G}$. Assume $s_1\in\G^{F_1}$ is a semisimple element and let $t_1\in\T$, $\w_1\in N_{\G}(\T)$ and 
$g_1\in\G$ be such that $g_1^{-1}F_1(g_1)=\w_1$;  $s_1=g_1t_1g_1^{-1}$ and  ${\rm Ad}(\w_1)\, F_1(t_1)=t_1$. 
Then $s_2:=g_0s_1g_0^{-1}$, $t_2:=t_1$, $g_2:=g_0g_1$ and $\w_2:=\w_1\sigmad$ satisfy 
$g_2^{-1}F_2(g_2)=\w_2$; $s_2=g_2t_2g_2^{-1}$ and  ${\rm Ad}(\w_2)\, F_2(t_2)=t_2$. 
Thus, if $\Pi_1$ is a base for the root system of $C_{\G}(t_1)$, then it is also a base for the root system of $C_{\G}(t_2)$ and 
comparing the actions of $w_j\theta_j$ for $j\in\I_2$ on $\Phi$ gives
$w_1 \theta_1=w_1 \sigma \theta_2= w_2\theta_2$, hence the isomorphism classes of the groups $[\Ha,\Ha]^{{\rm Ad}(\w_j) F_j}$  for $j\in \I_2$ coincide. 
\end{obs}

\begin{lema}\label{lem:e6-steinberg}
Assume that we are in Setting \ref{set:steinberg}, with $\Phi$ of type $E_6$.   Then $\Oc$ is not kthulhu.
\end{lema}

\pf 
We proceed as in the proof of Lemma \ref{lem:E6}, using the list in \cite{FJ} 
of pairs attached to a semisimple conjugacy class. The Steinberg endomorphism $F_{FJ}$ 
used therein differs from our choice of $F$ and there exists $\sigmad\in w_0\T$ such that $F_{FJ}={\rm Ad}\,\sigmad\, F$. 
By Remark \ref{obs:dictionary-e6} the pair  $(\Pi_{FJ},[w_{FJ}])$ attached to a semisimple conjugacy class in $\G^{F_{FJ}}$ 
is related to the pair attached to the corresponding semisimple conjugacy class in $\G^{F}$ by the law 
$(\Pi_{FJ},[w_{FJ}])=(\Pi,[ww_0])=(\Pi,[w_\theta])$. Also, the groups $\G_i^{F_w}$ can be extracted from the list in \cite{FJ}. 

We have to deal with the case in which $t^q=zt^{-1}$ for some $z\in Z(\G^F)$. By \eqref{eq:relationE6}
the element $w_{FJ}=ww_0$ normalises $\Ha$, whence $W_\Pi$. Also, $(ww_0)^3t=t$, so $\left|ww_0W_\Pi\right|\in\{1,\,3\}$. 
Observe that in all cases in which $\alpha_0\in\Pi$, then $t\in\prod_{i\neq2}\alpha_i^\vee(\kk)$ so it lies in the subgroup  
${\mathbb K}:=\langle \U_{\pm\alpha_j},\,j=1,3,4,5,6\rangle\simeq\SL_6(\kk)$. 
All mixed classes in ${\mathbb K}^{F_w}$ are not kthulhu by \cite{ACG-I} and Lemma \ref{lem:su}.

Assume first $\left|ww_0W_\Pi\right|=3$ so \eqref{eq:relationE6} forces $z\neq 1$. In particular,  $q\equiv 2(3)$ in this case.
All pairs as in Lemma \ref{lem:E6} with $w_{JF}=ww_0\neq1$ that have been discarded because of the order,  
can be again discarded for the same reason, as well as those that were discarded by using Remark \ref{obs:zeta} (3).   
As observed in \cite{FJ}, the groups $\G_i^{F_w}$ for $^{2}\!E_6(q)$ are obtained from those in the list for $E_6(q)$ by interchanging 
Chevalley and Steinberg's types in each factor of type $A,\,D$ or $E_6$. We are thus left with the cases in Table \ref{tab:remaining-2e6}.
There, $\beta_i$ for $i=10,11,19,20,21,22$ are as in Lemma \ref{lem:E6} and
\begin{align*}
&\beta_{13}=\alpha_1+\alpha_2+\alpha_3+\alpha_4; 
&&\beta_{15}=\alpha_2+\alpha_4+\alpha_5+\alpha_6;\\
& \beta_{24}= \alpha_1+\alpha_2+\alpha_3+2\alpha_4+2\alpha_5+\alpha_6;
& &\beta_{28}=\alpha_3+\alpha_4+\alpha_5+\alpha_6.
 \end{align*}

\begin{small} 
\begin{table}[ht]
		\caption{Pairs $(\Pi,w_{FJ}W_{\Pi})$ in type $^{2}E_6$, for $w_{FJ}\neq1$ that are dealt with separately}\label{tab:remaining-2e6}
			\begin{tabular}{|c|c|c|}
				\hline $\Pi$  & $w_{FJ}W_\Pi$ &  $[\Ha,\Ha]^{F_w}$ (up to isogeny) \\
			\hline
$A_1=\{-\alpha_0\}$ & $s_1s_3s_5s_6$& $\SL_2(q)$\\
				
	
				\hline
	$A_2=\{\alpha_2,\,-\alpha_0\}$& $s_1s_3s_5s_6,$ &  $\SU_3(q)$\\
				\hline			
	$3A_1=\{\alpha_4,\alpha_6,-\alpha_0\}$& $s_{\beta_{11}} s_{\beta_{19}}s_{\beta_{20}} s_{\beta_{10}}$ & 	$\SL_2(q^3)$\\		
				\hline
	$4A_1=\{\alpha_1,\,\alpha_4,\,\alpha_6,\,-\alpha_0\}$ & $s_{\beta_{19}}s_{\beta_{11}}s_{\beta_{21}}s_{\beta_{22}}$  & $\SL_2(q)\times \SL_2(q^3)$\\
				\hline
		$3A_2=\{\alpha_1,\alpha_3,\alpha_5,\alpha_6,\alpha_2,-\alpha_0\}$& $s_{\beta_{13}}s_{\beta_{15}}s_{\beta_{28}}s_{\beta_{24}}$ & 	$\SU_3(q^3)$\\			
				\hline		
			\end{tabular}
	\end{table}
\end{small}

If $\Pi=\{-\alpha_0\} $ or $\{\alpha_2,\,-\alpha_0\}$ or $\{\alpha_4,\alpha_6,-\alpha_0\}$, then Lemma \ref{lem:no-centro} 
applies with $\sigma=s_1s_3$ in the first two cases and  $\sigma=s_1$ in the third one. 

Let $\Pi=\{\alpha_1,\alpha_4,\alpha_6,-\alpha_0\}$. If $v\in\langle\U_{\pm\alpha_0}\rangle$, then Lemma   \ref{lem:no-centro} 
applies with $\sigma=s_{\beta_{11}}s_{\beta_{21}}$. If, instead, $v\in \langle\U_{\pm\alpha_j},\,j=1,4,6\rangle$, then $tv\in{\mathbb K}$ and it is mixed therein.

We claim that the case $\Pi=\{\alpha_1,\alpha_3,\alpha_5,\alpha_6,\alpha_2,-\alpha_0\}$ with $t^{q+1}\in Z(\G^F)-1$ cannot occur. 
Indeed,  here $t=\alpha_1^\vee(\xi)\alpha_3^\vee(\xi^2)\alpha_5^\vee(\zeta)\alpha_6^\vee(\zeta^2)$ with $\xi^3=\zeta^3=1$, 
$\xi\neq\zeta$, so $t^{q+1}=1$ because $q\equiv2(3)$.

We deal now with the case $\left|ww_0W_\Pi\right|=1$. Then $(ww_0)^{-1}\in W_\Pi$ so \eqref{eq:relationE6} implies  
$t^{q+1}=z=1$  and we take $w=w_0$, so $C_W(w\theta)=W$.  Since either \eqref{cond:Jv=1} or   \eqref{cond:Jv=2} holds, 
we discard all choices of $\Pi$ for which no $\G_i$ is isogenous to $\SU_m(q)$,  $\PSL_3(2)$ or $\SL_2(q)$. This leaves us 
with Table \ref{tab:2e6-id}.  In most remaining cases we apply Remark \ref{obs:simple_root} and we list  the simple roots we use in 
the third column of Table \ref{tab:2e6-id}.  We deal with the  remaining cases separately.   Observe that since $w=w_0$, the subgroup 
${\mathbb K}$ is $F_w$-stable and $w\theta$ acts as $-1$ on its root system, so ${\mathbb K}^{F_w}\simeq\SU_6(q)$. 

\begin{tiny}
\begin{table}[ht]
		\caption{Pairs $(\Pi, W_{\Pi})$ in type $^{2}E_6$, for $w_{FJ}=1$}\label{tab:2e6-id}
		\begin{center}
			\begin{tabular}{|c|c|c|}
				\hline $\Pi$  & $[\Ha,\Ha]^{F_w}$ (up to isogeny)  & Reason to discard \\
				\hline
			$A_1=\{-\alpha_0\}$ & $\SL_2(q)$& $s_1$\\
				\hline
$2A_1=\{\alpha_6,-\alpha_0\}$ & $\SL_2(q)^2$ &$s_1$\\ 
				\hline
$A_2=\{\alpha_2,-\alpha_0\}$ & $\SU_3(q)$ & $s_1$ \\
\hline
$3A_1=\{\alpha_4,\alpha_6,-\alpha_0\}$	&$\SL_2(q)^3$ &  $s_1$\\
	\hline
$A_2+A_1=\{\alpha_2,\alpha_6,-\alpha_0\}$	
&  $\SU_3(q)\times\SL_2(q)$  & $s_1$\\
\hline
$A_3=\{\alpha_2,\alpha_4,-\alpha_0\}$	
&  $\SU_4(q)$ & $s_1$	\\	
\hline
$4A_1=\{\alpha_1,\,\alpha_4,\,\alpha_6,\,-\alpha_0\}$ & $\SL_2(q)^4$  & $s_3$ \textrm{ or }$s_5$ \textrm{ or }$s_2$	\\
\hline
$A_2 + 2A_1=\{\alpha_1,\,\alpha_3,\,\alpha_6,\,-\alpha_0\}$ & $\SU_3(q)\times\SL_2(q)^2$ &$s_4$ \textrm{ or }$s_5$ \textrm{ or }$s_2$\\
\hline
$2A_2=\{\alpha_2,\,\alpha_5,\,\alpha_6,\,-\alpha_0\}$ & $\SU_3(q)^2$ & $s_1$\\
\hline
$A_3 + A_1=\{\alpha_4,\,\alpha_5,\,\alpha_6,\,-\alpha_0\}$ & $\SU_4(q)\times\SL_2(q)$ &$s_1$ \\
\hline
$A_4\ =\{\alpha_2,\,\alpha_4,\,\alpha_5,\,-\alpha_0\}$ & $\SU_5(q)$ & $s_1$\\
\hline
$2A_2 + A_1 =\{\alpha_1,\,\alpha_3,\,\alpha_5,\,\alpha_6,\,-\alpha_0\}$ & $\SU_3(q)^2\times\SL_2(q)$  & $s_2$ \textrm{ or }$s_4$\\
\hline
$A_3 + 2A_1=\{\alpha_1,\,\alpha_3,\,\alpha_4,\,\alpha_6,\,-\alpha_0\}$ & $\SU_4(q)\times\SL_2(q)^2$ &$s_2$ or $s_5$ or $tv\in\SU_6(q)$ \\
\hline
$A_4 + A_1=\{\alpha_1,\,\alpha_3,\,\alpha_4,\,\alpha_5,\,-\alpha_0\}$ 
&$\SU_5(q)\times\SL_2(q)$ & $s_6$ or $tv\in\SU_6(q)$\\
\hline
$A_5\ =\{\alpha_2,\,\alpha_4,\,\alpha_5,\,\alpha_6,\,-\alpha_0\}$ & $\SU_6(q)$&$s_1$\\
\hline
$3A_2=\{\alpha_1,\,\alpha_3,\,\alpha_5,\,\alpha_6,\,\alpha_2,\,-\alpha_0\}$ & $\SU_3(q)^3$  &$s_{\beta_{13}}s_{\beta_{14}}s_{\beta_{15}}$ or $tv\in\SU_6(\kk)$ \\
\hline
$A_5 + A_1=\{\alpha_1,\,\alpha_3,\,\alpha_4,\,\alpha_5,\,\alpha_6,\,-\alpha_0\}$ & $\SU_6(q)\times \SL_2(q)$ &$w_0\in W_\Pi$ and Lemma \ref{lem:A3} \\
\hline
			\end{tabular}
		\end{center}
	\end{table}
 \end{tiny}

Let $\Pi=\{\alpha_1,\,\alpha_3,\,\alpha_4,\,\alpha_6,\,-\alpha_0\}$. If $v$ has a component in either $\langle \U_{\pm\alpha_0}\rangle$ or 
$\langle \U_{\pm\alpha_6}\rangle$, then Lemma \ref{lem:no-centro} applies with either $\sigma=s_5$ or $\sigma=s_2$. Assume 
$|\J_v|=1$ and $v\in\langle \U_{\pm\alpha_j},\,j=1,3,4\rangle$. Since $-\alpha_0\in\Pi$ we have $tv\in {\mathbb K}$ and it is mixed therein.

Let $\Pi=\{\alpha_1,\,\alpha_3,\,\alpha_4,\,\alpha_5\,-\alpha_0\}$. If $v$ has a component in $\langle \U_{\pm\alpha_0}\rangle$ then
Lemma \ref{lem:no-centro} applies with  $\sigma=s_6$. If $v\in\langle \U_{\pm\alpha_j},\,j=1,3,4,5\rangle$, since $-\alpha_0\in\Pi$ we 
have $tv\in {\mathbb K}$ and it is mixed therein.

Let $\Pi=\{\alpha_1,\,\alpha_2,\,\alpha_3,\,\alpha_5\,\alpha_6,\,-\alpha_0\}$. Here 
$t=\alpha_1^\vee(\xi)\alpha_3^\vee(\xi^2)\alpha_5^\vee(\zeta)\alpha_6^\vee(\zeta^2)$ with $\xi^3=\zeta^3=1$, $\xi\neq\zeta$. 
The Weyl group involution $\sigma=s_{\beta_{13}}s_{\beta_{14}}s_{\beta_{15}}$ satisfies 
$\sigma(\alpha_1)=\alpha_5$, $\sigma(\alpha_3)=\alpha_6$ and $\sigma(\alpha_0)=\alpha_2$. Hence, if $v$ has a component in 
$\langle \U_{\pm\alpha_0}\rangle$, then Lemma \ref{lem:no-centro} applies. If, instead, $v$ has no component  in
$\langle \U_{\pm\alpha_0}\rangle$, then $tv\in{\mathbb K}$ and it is mixed therein.

Let $\Pi=\{\alpha_1,\,\alpha_3,\,\alpha_4,\,\alpha_5,\,\alpha_6,\,-\alpha_0\}$. A direct computation shows that, 
$t=\alpha_1^\vee(\zeta)\alpha_3^\vee(\zeta^2)\alpha_4^\vee(-1)\alpha_5^\vee(\zeta^4)\alpha_6^\vee(\zeta^{5})$ for 
$\zeta\in\kk$ such that  $\zeta^6=1$, $\zeta^3=-1\neq1$. Hence $q$ is odd and \eqref{cond:Jv=1} implies that $v\in\langle\U_{\pm\alpha_0}\rangle$.  
Also, the condition $t^q=t^{-1}$ gives $F(t)=\vartheta t^q=\vartheta t^{-1}=t$, so $t\in \G^F$ and we may take 
$\x_s=t$ and replace $w=w_0$ by $w=1$. By Remark \ref{obs:representativev} we assume $v=x_{\alpha_0}(\xi)$ 
for some $\xi\in\kk^\times$. Lemma \ref{lem:A3} applies taking $\beta_1=\alpha_0$, $\beta_2=\alpha_2$ and $\beta_3=\alpha_4$.
\epf

Putting together Lemma \ref{lem:specifico},  Proposition \ref{prop:center} 
and Lemmata \ref{lem:d2k-steinberg} up to \ref{lem:e6-steinberg}, we prove the main Theorem of this Section:
\begin{theorem}
Let $\Gb$ be a Steinberg group and $x=x_{s}x_{u} \in \Gb$ with $x_s,\,x_u\neq1$.
Then $\Oc_x^{\Gb}$ is not kthulhu, and consequently $\Oc_x^{\Gb}$ collapses.\qed
\end{theorem}

\section{Nichols algebras over Chevalley and Steinberg groups}\label{sec:nichols}
In order to prove Theorems \ref{thm:collapse} and \ref{thm:collapse-gp}, we have to consider now Nichols algebras attached to Yetter-Drinfeld 
modules over some specific finite groups, in other words, the cocycle $\bq$ is not arbitrary, but determined by a representation of the centraliser
of a fixed element in the conjugacy class.

\subsection{Unipotent orbits in Chevalley and Steinberg groups}
Let $\Gb$ be a Chevalley or Steinberg group isomorphic to neither $\PSL_{3}(2)$ nor $\PSL_{2}(3)\simeq \A_{4}$, the latter being non-simple.
Note that, as  $\PSL_{3}(2)\simeq \PSL_{2}(7)$,
the only  unipotent conjugacy class in $\PSL_{3}(2)$ which is kthulhu is isomorphic to a semisimple
class in $\PSL_{2}(7)$.

We now deal with Nichols algebras associated with unipotent classes  $\Oc$ in $\Gb$ and representations of the centraliser of $x\in\Oc$. 
In most of the cases we will  prove that they are infinite-dimensional by reducing to a subgroup of $\Gb$ for which a similar statement is known.
We begin with a case for which this strategy cannot be implemented, see Remark \ref{obs:PSL23} below.





\begin{lema}\label{lem:uno}Let $G=\Sp_{4}(3)$, $\Gb=\PSp_{4}(3)$ and  $\pi\colon G\to\Gb$ the natural projection.
Let $\Oc$ be the unipotent class of an element $g_0\in\Gb$ whose Jordan form has blocks of size $(2,1,1)$. 
Then, $\dim\toba(M(\Oc,\rho))=\infty$ for every $\rho\in\Irr {C_{G}(g_0)}$.  
The same statement holds replacing $g_0$ by $\pi(g_0)$ and $C_G(g_0)$ by $C_{\Gb}(\pi(g_0))$.
\end{lema}
\pf

There are two such unipotent classes in $G$ and they are represented by $g_0=\id\eta e_{1,4}$, with $\eta^2=1$. 
We  find a suitable abelian subrack of $\Oc$
and apply the results in \cite{H} on Nichols algebras of diagonal type.

For both choices of $g_{0}$, 
we consider the abelian subgroup $H$ of $C_G(g_0)$ consisting of matrices of block form
$g(X):=\left(\begin{smallmatrix}\id_2&X\\
0&\id_2
\end{smallmatrix}\right)$ where $X=\left(\begin{smallmatrix}a&b\\
c&a
\end{smallmatrix}\right)$ for $a,b,c\in\F_3$. Then 
\begin{align*}
\Oc\cap H=\left\{\id+\eta e_{1,4}, \id+\eta e_{2,3}, g\left(\begin{smallmatrix}\eta&\eta\\
 \eta&\eta\end{smallmatrix}\right),g\left(\begin{smallmatrix}-\eta&\eta\\
 \eta&-\eta\end{smallmatrix}\right) \right\}=\{g_i=x_i\trid g_0,\,i\in\I_{0,3}\}                                                                                                                     
 \end{align*} with: 
\begin{align*}
x_0=1,&&x_1={\rm diag}(\Jf_2,\Jf_2),&&u=\left(\begin{smallmatrix}
1&1\\
&1\\
&&1&2\\
&&&1
\end{smallmatrix}\right),
&&x_2=ux_1,&&x_3=u^2x_1.
\end{align*}
Observe that 
\begin{align*}
u\trid g_0=g_0,&& x_1\trid g_2=g_2,&& x_1\trid g_3=g_3, && x_1\trid g_1=g_0, && g_0g_1g_2g_3=1.
\end{align*}
As $H$ is abelian, the restriction of $\rho$ to $H$ decomposes as a direct sum of $1$-dimensional representations; let ${\mathbb C}v$ 
be one of these and let $\rho(g_i)v=\zeta_iv$ for every 
$i\in\I_{0,3}$.
We necessarily have $\zeta_i^3=1$ and $\zeta_0\zeta_1\zeta_2\zeta_3=1$. 
Then ${\rm span}_{\mathbb C}\{x_i\otimes v,\,i\in\I_{0,3}\}$ is a braided vector subspace of $M(\Oc,\rho)$, where
\begin{align*}c((x_i\otimes v)\otimes (x_j\otimes v))=q_{ij}(x_j\otimes v)\otimes (x_i\otimes v),&& q_{ij}v=\rho(x_j^{-1}x_i\trid g_0)v.
\end{align*}
By direct computation we obtain
\begin{center}
\begin{tabular}{|c|c|c|c|c|}
\hline
$q_{ij}$&0&1&2&3\\
\hline
0&$\zeta_0$&$\zeta_1$&$\zeta_1$&$\zeta_1$\\
\hline
1&$\zeta_1$&$\zeta_0$&$\zeta_0^2\zeta_1^2\zeta_2^2$&$\zeta_2$\\
\hline
2&$\zeta_2$&$\zeta_2$&$\zeta_0$&$\zeta_0^2\zeta_1^2\zeta_2^2$\\
\hline
3&$\zeta_0^2\zeta_1^2\zeta_2^2$&$\zeta_0^2\zeta_1^2\zeta_2^2$&$\zeta_2$&$\zeta_0$\\
\hline
\end{tabular}
\end{center}
so 
\begin{align*}
q_{01}q_{10}=\zeta_1^2,&&q_{02}q_{20}=\zeta_1\zeta_2,&&q_{03}q_{30}=\zeta_0^2\zeta_2^2,\\
q_{12}q_{21}=\zeta_0^2\zeta_1^2,&&q_{13}q_{31}=\zeta_0^2\zeta_1^2,&&q_{23}q_{32}=\zeta_0^2\zeta_1^2.
\end{align*}
Let $W={\rm span}_{\mathbb C}\{x_i\otimes v,\,i\in\I_{1,3}\}$. If either $\zeta_0=1$ or  $\dim\toba(W)=\infty$, then  $\dim\toba(M(\Oc,\rho))=\infty$, 
so  we assume that $\zeta_0$ is a primitive third root of $1$ and that  $\dim\toba(W)<\infty$. 
Since the generalized Dynkin diagram of $W$ does not occur in \cite[Table 2]{H}, it must be 
disconnected. This forces $\zeta_1=\zeta_0^2$, but then the generalized Dynkin diagram of 
$W'={\rm span}_{\mathbb C}\{x_0\otimes v,x_1\otimes v\}$ is connected and  does not occur in 
\cite[Table 1]{H}, so $\dim\toba(W')=\infty$, and a fortiori  $\dim\toba(M(\Oc,\rho))=\infty$. 
The statement for $\Gb$ follows similarly because the restriction of $\pi$ to $\Oc\cap H$ is injective.
\epf

\begin{prop}\label{prop:psp}
Let $n\geq 2$, let $M$ be either $\Sp_{2n}(3)$ or $\PSp_{2n}(3)$, and let $\Oc$ be the unipotent class of an element $x_u$ in $M$ 
whose Jordan form has blocks of size $(2,1^{2n-2})$. 
Then, $\dim\toba(M(\Oc,\rho))=\infty$ for every $\rho\in\Irr {C_{M}(x_u)}$. 
\end{prop}
 \pf If $n=2$ this is Lemma \ref{lem:uno}. If $n>2$ we use the embeddings of $\Sp_{4}(3)$ into $\Sp_{2n}(3)$ and $\PSp_{2n}(3)$  given by 
 $\left(\begin{smallmatrix}
 A&B\\
 C&D\end{smallmatrix}\right)\mapsto \left(\begin{smallmatrix}
 A&&B\\
 &\id_{2n-4}\\
 C&&D\end{smallmatrix}\right)$ and \cite[Lemma 3.2]{AFGV-ampa}.
 \epf
 
 We are now in a position to prove our next theorem.

\begin{theorem}\label{thm:unip_collapse}
Let $x$ be a unipotent element in a Chevalley or Steinberg group $\Gb \not \simeq \PSL_{3}(2)$, $\PSL_{2}(3)$. 
Then $\dim \toba(\Oc_{x},\rho) = \infty$ for all $\rho \in \Irr C_{\Gb}(x)$.
\end{theorem}

\pf We may assume that $x$ is non-trivial. 
By Theorem \ref{thm:slspsu} it is enough to prove the statement for the conjugacy classes in Table \ref{tab:unip-kthulhu}. 
If $q=3$, then we invoke Proposition \ref{prop:psp}.  Let $q\neq 3$.  These conjugacy classes are represented by an element 
$x_{\beta}(\xi)$ with $\xi \in \F_{q}^{\times}$ and $\beta$ a positive root by Remark \ref{obs:properties-kthulhu}. 
Since $x_{\beta}(\xi) \in \langle \U_{\beta}^{F}, \U_{-\beta}^{F}\rangle \leq \SL_{2}(q) $ or $\leq \PSL_{2}(q)$, 
the statement follows from \cite[Lemma 3.2]{AFGV-ampa} together with: \cite[Proposition 3.1]{FGV1} for $q$ even; 
\cite[Lemma 2.2]{FGV2} for $p>3$, 
and the proof of \cite[Lemma 3.7]{FGV1} for $q=3^{2h+1}$, $h>0$. 
\epf

\begin{obs}\label{obs:PSL23}
$(a)$ We do not know whether $\dim\toba(M(\Oc),\rho)=\infty$ for $\Oc_u$ a non-trivial unipotent conjugacy class in $\Gb=\PSL_{2}(3)$ 
and $\rho$ an irreducible representation of $C_{\Gb}(u)$. Indeed, the proof of \cite[Proposition 4.3]{FGV2}  
does not cover the case of  non-trivial unipotent conjugacy classes in $\PSL_{2}(3)$ because they are not real. 
These conjugacy classes correspond to the tetrahedral rack associated with a class of 3-cycles in $\A_{4}$. 

$(b)$ There are examples of finite-dimensional Nichols algebras associated with the rack $\Oc_u$ as in $(a)$ 
and a cocycle that does not come from a representation of $C_{\Gb}(u)$, see \cite[Proposition 36]{HLV}.  
\end{obs}

\subsection{Proofs of Theorems \ref{thm:collapse} and \ref{thm:collapse-gp}}

\

\medbreak
\noindent{\em Proof of Theorem \ref{thm:collapse}.} By \cite{HS}, 
$V$ should be simple, say $V \simeq M(\Oc, \rho)$.
By Theorem \ref{thm:mixed-Chev-Steinberg}, we know that $\oc$ is either semisimple or unipotent, 
but the latter is discarded by Theorem \ref{thm:unip_collapse}.
 \hfill$\Box$

\

\smallbreak

\noindent{\em Proof of Theorem \ref{thm:collapse-gp}.}
Let $V = M(\oc, \rho) \in \ydgb$. Assume that $\dim \toba(V) < \infty$.
By Theorem \ref{thm:collapse}, $\oc = \oc_x$ is semisimple, hence $\ord x$ is odd since $q$ is even.
Since $\Gb$ is one of the groups in \eqref{eq:gps-1inW}, $-\id$ belongs to the Weyl group $W$.
Thus $\oc$ is real by \S \ref{subsec:w0}. 
This contradicts \cite{AZ}, hence  $\dim \toba(V) = \infty$. We conclude by \cite[Lemma 1.4]{AFGV-ampa}.
 \hfill$\Box$

\

\smallbreak
Notice that \cite{HS} is not needed for the last proof.

\end{document}